# Non-contractible periodic orbits, Gromov invariants, and Floer-theoretic torsions

Yi-Jen Lee

Dept. of Mathematics, Princeton University

Princeton NJ 08544, U.S.A.

ylee@math.princeton.edu

This version: August 2003

## Abstract

In a previous paper [L], the author introduced a Floer-theoretic torsion invariant $I_F$, which roughly takes the form of a product of a power series counting perturbed pseudo-holomorphic tori, and the Reidemeister torsion of the symplectic Floer complex. We pointed out the formal resemblance of $I_F$ with a generating function of genus 1 Gromov invariant [IP2]; furthermore, for heuristic reasons one also expects a relation with the 1-loop generating function in the A-model side of mirror symmetry [BCOV, FuOOO], which counts genus 1 holomorphic curves.

The present article makes this expected relation precise in the simplest cases, in two variants of the $I_F$ defined in [L]: the lagrangian intersection version, $I_F(L, L')$, and an $S^1$-equivariant version, $I_F^{S^1}$.

As a by-product, we obtain some existence results of *noncontractible* periodic orbits in symplectic dynamics. For example, the results of Gatien-Lalonde [GL] are extended to a much wider class of manifolds.

The two versions $I_F(L, L')$ and $I_F^{S^1}$ are only minimally developed in this paper, leaving fuller accounts to future work (e.g. [LS]). The lagrangian intersection version, $I_F(L, L')$, should be viewed as a simplest example of a rigorous definition of the higher-loop "open Gromov-Witten *invariants*" proposed by physicists [OV, Wi].

## Contents

# 1 Introduction

The existence of periodic orbits in a Hamiltonian dynamical system, or more generally, in a flow generated by a symplectic vector field, is a central problem in dynamical systems and classical mechanics. The abundance of periodic orbits in symplectic dynamical systems is a manifestation of the "rigidity" phenomena in symplectic topology. A famous example is the various versions of the Arnold conjecture. Very roughly, it states that the number of periodic orbits of a Hamiltonian system is at least the total betti number of the underlying symplectic manifold; in contrast, via differential topological methods one only expects it to be greater than or equal to the Lefschetz number of the symplectomorphism generated by the Hamiltonian vector field.

Basing on an infinite-dimensional analog of the Morse homology, Floer introduced the powerful tool of Floer homology to attack this type of problems [F89b, F88b]. Floer's technique has since been refined/modified to deal with various variations/generalizations of the Arnold-type problem; see e.g. [LeO, LiT, FuO]. However, the results thus obtained are typically limited to *contractible* periodic orbits, or equivalently, symplectic fixed points in a standard Nielsen class. The reason is that this method depends on the computation of the Floer homology, which is reduced to the case of small symplectic vector fields by the invariance of Floer homology. Obviously, when the symplectic vector field is small, there are no periodic orbits other than the contractible ones [BuH]. In other words, one may define a Floer homology using the space of loops in a nontrivial homotopy class, but it would be trivial.

Due to the lack of available tools, literature on the existence of *noncontractible* periodic orbits has been very rare. However, very recently some isolated results begin to appear (e.g. [GL, BPS]). These results use the holomorphic-curve techniques introduced by Gromov [Gr], and the difficulty in computing the relevant counting invariant of curves often places severe limitation on the applicability of such results. For example, the results of Gatien-Lalonde only work for $T^*V$, where $V$ is the mapping torus of an involution of a flat manifold.

In this paper, we provide some more existence results of noncontractible periodic orbits in symplectic dynamics, via certain Floer-theoretic *torsions* $I_F$. These torsions are variants of the torsion invariant $I_F$ introduced by the author in the foundational paper [L]. More specifically, we shall be concerned with two versions of $I_F$: the first, $I_F(L, L')$, is defined from the Floer theory of lagrangian intersections, and the second, $I_F^{S^1}$, is based on an $S^1$-equivariant version of Floer theory, where $S^1$ is the rotation of the loops. (See §2.2, 2.3 below respectively.)

By construction, $I_F$ is a centaur that combines a Floer-theoretic part and a curve-counting part:

$$I_F = \tau_F \zeta_F,$$

where $\tau_F$ is the Reidemeister torsion of (a twisted version of) the Floer complex, and $\zeta_F$ is a "zeta function" that counts (perturbed) pseudo-holomorphic tori or annuli, depending on the version of Floer complex used. They take values in certain gener-alization of power-series rings. (see §2.1). The point of the definition is that neither $\tau_F$ nor $\zeta_F$ alone is invariant, but their product is. (See Theorems 2.2.2, 2.3.1 below).

One may on one hand view $I_F$ as a refinement of the Floer homology, in the same way that the Reidemeister torsion refines the homology; on the other hand, $I_F$ may also be regarded as a "corrected" 1-loop generating function of Gromov invariants. These dual perspectives of $I_F$ yield different types of applications.

A typical application from the first perspective is: If the leading term of $I_F$ is nontrivial ($\neq 1$), then critical points exist in the relevant Floer theory, which can be symplectic fixed points, periodic orbits, or lagrangian intersection points. This type of application is discussed in [L, LS]; see also Corollary 2.2.4 below for a sample result.

In the present article, we shall however focus on the second perspective, and obtain applications that apply in situations complementary to the previous case: here the leading term of $I_F$ is trivial, and we rely on the nontriviality of higher-order terms. At the heart of these applications is a relation between $I_F$ and Gromov-type invariants that count tori or annuli. This relation was briefly mentioned in [L], and is made precise here under simplifying assumptions in Theorem 2.3.3 and Lemma 5.1.3 below. To understand the relation in general remains a fascinating challenge.

We shall wield this relation as a two-edge sword: In the lagrangian intersection version described in §2.2, it is the Floer-theoretic torsion $I_F(L, L')$ that is easier to compute; the relation mentioned above therefore yields a computation of the relevant Gromov invariant. One may then follow the arguments of Gromov-Gatien-Lalonde to arrive at significant strengthening of the results of [GL]. (See Theorems A, B, Corollary C in §2.4). On the other hand, in the $S^1$-equivariant version described in §2.3, we shall use the computation of the genus 1 Gromov invariants from the existent literature to establish the nontriviality of $I_F^{S^1}$, which implies existence of periodic orbits. (Corollary D in §2.4).

We next briefly comment on some other motivations for defining these two particular versions of $I_F$.

Expanding on the second perspective of $I_F$ mentioned above, we believe that the definition of $I_F(L, L)$ here showcases a potential approach to define general "open Gromov-Witten invariants". Though a mathematical foundation has been lacking, physicists (see e.g. [OV, Wi]) have proposed an "open Gromov-Witten theory" which counts holomorphic curves with boundaries on lagrangian submanifolds, and produced fascinating conjectures on its relation with the Chern-Simons theory. The main difficulty for a rigorous definition of such open Gromov invariants is the appearance of boundary nodes on the holomorphic curve while varying the parameters in the theory (i.e. almost complex structures, perturbations, lagrangian boundary conditions), which is a codimension 1 phenomenon. Thus, in general the count of such curves can not be an *invariant*. The definition of $I_F(L, L)$ here shows that by adding Floer-theoretic "corrections" $\tau_F$ (which vanishes in special cases), one may obtain an invariant for the 1-loop case (i.e. counting annuli). One expects a generalization of this picture to higher loop cases: there should be an interpretation of Gromov invariants in terms of generalized Morse theory (*a la* Cohen, Fukaya, see references in [Fu]) on the (relative) loop space. The corrected invariants should count, in addition to open holomorphic curves, contributions from a (deformed) Floer complex [L2]. (See [Fu] for a finite-dimensional, 2-loop analog). The Floer-theoretic correction should be understood as a compensation to the effects of a *new* type of boundary nodes arising in the higher-loop situation, in contrast to the more familiar bubbling-off-disks, which occurs already in the 0-loop (disk) case. In the ribbon-tree description of bordered Riemann surfaces, this new mechanism corresponds to "breaking an (internal) edge".

It is also interesting to note that in Chern-Simons theory, the perturbative invariant originally defined from path integrals need to be regularized by combining with certain Ray-Singer analytic torsions, analogous to the correction by Reidemeister torsion of Floer complex in our story.

On the other hand, we arrived at the $S^1$-equivariant version of $I_F$ from the following considerations: first, it ties in better with the usual Gromov invariants that has been computed in some cases, thus providing many interesting examples. Second, when one works with the space of loops in a primitive homotopy class, rotation generates a free $S^1$-action on this loop space. In this case, the usual, non-equivariant version of $I_F$ defined in [L] is trivial, and it is more interesting to consider the equivariant version, modeling on Morse theory on the quotient of the loop space by the $S^1$-action.

Historically, an $S^1$-equivariant Floer theory was first proposed by Givental [Gi], to give a heuristic explanation of the predictions of mirror symmetry. The foundation has however been hitherto missing. Givental considered the space of contractible loops, on which the $S^1$-action is not free, and one needs a full-blown equivariant theory. In contrast, in this paper we use the free $S^1$-action to reduce everything to the quotient space. Viterbo [V] described a different approach of defining an $S^1$-equivariant Floer theory via a Floer theory for families. (Namely, consider the Floer theory of the homotopy quotient). The details of this approach are yet to appear, but it does not seem to suit our purposes.

Lastly, we emphasize that even restricted to just the 1-loop version, the story presented here is only a vastly simplified version of a general theory. Two major omissions are: first, we avoid the possibility of bubbling by making several restrictive assumptions. In general, one shall need the virtual moduli method and the machinery of [FuOOO]. Second, we define $I_F$ only using the moduli spaces of expected dimension 0, resulting in a more restrictive class of manifolds with interesting $I_F$. More generally, one may try to construct generalizations of $I_F$ using higher dimensional moduli spaces. This would involve a "quantum" version of Floer theory.

The plan for the rest of the paper is as follows. We summarize the necessary algebraic framework and the results obtained in this paper in section 2. In sections 3 and 4 respectively, we lay the foundations of the two versions of Floer theories we need. In sections 5 and 6, we establish relations of the two versions of $I_F$ with the corresponding versions of Gromov invariants, which is then used to prove the existence theorems of periodic orbits. The invariance of $I_F(L, L')$ and $I_F^{S^1}$ is proven in section 7, via simple adaptations of the arguments in [L].

## 2 Statements of results

### 2.1 Definition of $I_F$: the algebraic framework

We quickly summarize the general framework for defining the torsion invariant $I_F$. More details may be found in [L, HL2].

#### 2.1.1 The Novikov ring.

Let $G$ be an abelian group, $R$ a ring, and $N : G \to \mathbb{R}$ a homomorphism. The *Novikov ring* $\mathrm{Nov}(G, N; R)$ is the set of formal sums $\sum_{g \in G} a_g \cdot g$, with $a_g \in R$, such that for every $C \in \mathbb{R}$, the set $\{g \in G \mid N(g) < C \text{ and } a_g \neq 0\}$ is finite. $\mathrm{Nov}(G, N; R)$ is a ring with the obvious addition and the convolution product. (See e.g. [MS2].)

Notice that $\mathrm{Nov}(G, 0, R) = R[G]$, and there is an inclusion

$$i_N : R[G] \hookrightarrow \mathrm{Nov}(G, N; R).$$

The Novikov ring should thus be viewed as a completion of the group ring.

The *degree* of $a$, denoted $\deg(a)$, is defined to be the minimum of $N(g)$ among $g$ such that $a_g \neq 0$.

Given $a = \sum_g a_g g \in \mathrm{Nov}(G, N; R)$, the "leading term" of $a$ is defined to be

$$\mathrm{lt}(a) := \sum_{N(g) = \deg(a)} a_g g.$$

$a - \mathrm{lt}(a)$ is called "higher order terms". Notice that lt defines a homomorphism

$$\mathrm{Nov}(G, N; R)/(\pm G) \to R[\ker N]/(\pm \ker N).$$

For our applications, we assume from now on

$$G \text{ is a finitely generated abelian group}; R = \mathbb{Z} \text{ or } \mathbb{Q}. \tag{1}$$

In this case, the Novikov ring is commutative.

*Notation.* Given a commutative ring $R$, $Q(R)$ denotes its total ring of fractions, namely its localization at all non-zero-divisors.

We shall often need to consider the rings of fractions of Novikov rings.

First, observe that a splitting

$$G = \ker N \oplus G/\ker N$$

induces an embedding:
$$Q(\text{Nov}(G, N; R)) \hookrightarrow \text{Nov}(G/\ker N, N; Q(R[\ker N])), \qquad (2)$$
and different embeddings are related by the natural action of the space of splittings ker $N$ on the right hand side.

Furthermore, in the case of (1), both sides of (2) are finite sums of fields. (As a special case, $Q(R[G])$ is a finite sum of fields, see e.g. [Tu2], [H] Lemma A.4). The embedding (2) is compatible with the decompositions on both sides as sums of fields.

**Remark.** In comparison with the first ring, the second ring in (2) has nicer properties (e.g. existence of the notions of degree, order, and limit), which the invariance proofs in [HL2, H, L] made use of. There is a confusion between the two rings in [HL2], which also propagate to later papers such as [H]. In these papers, the notion of order for elements in $Q(\text{Nov}(G, N, R))$ should be understood in terms of the larger ring above through the embedding (2).

Through the embedding (2), we may extend the notion of leading-term to:
$$\text{lt}: Q^{\text{Nov}}(G, N; R)/\pm G \to Q(R[\ker N])/\pm \ker N.$$
Note that the above map is independent of the choice of splitting since we mod out ker $N$. The embedding $i_N$ also extends to the ring of fractions,
$$i_N: Q(R[G]) \hookrightarrow Q(\text{Nov}(G, N; R)).$$
Later we shall also use the same notation $i_N$ to denote the induced map from $Q(R[G])/(\pm G)$ to $Q(\text{Nov}(G, N; R))/(\pm G)$.

Let $\text{Nov}^+(G, N; \mathbb{Q}) \subset \text{Nov}(G, N; \mathbb{Q})$ denote the subset of elements of positive degree.

Let $\text{Nov}^1(G, N; \mathbb{Q}) \subset \text{Nov}(G, N; \mathbb{Q})$ be the subgroup consisting of elements of the form $1 + c$, $c \in \text{Nov}^+(G, N; \mathbb{Q})$. The exponential
$$\exp: \text{Nov}^+(G, N; \mathbb{Q}) \to \text{Nov}^1(G, N; \mathbb{Q}) \hookrightarrow Q(\text{Nov}(G, N; \mathbb{Q}))$$
is well defined via the usual power series. Conversely, the logarithm
$$\ln: \text{Nov}^1(G, N; \mathbb{Q}) \to \text{Nov}^+(G, N; \mathbb{Q})$$
also makes sense formally.

Novikov rings arise naturally in Morse-Novikov theory in the following way: Given a generic closed 1-form $\theta$ on a compact manifold $Y$, let $\tilde{Y}$ be a regular covering of $Y$ with covering transformation group $G$ such that the lift of $\theta$ to $\tilde{Y}$ is exact, then $\theta$ defines a homomorphism from $G$ to $\mathbb{R}$ which we still denote by $\theta$. $\text{Nov}(G, -\theta; R)$ is the coefficient rings of a twisted version of Morse complex $(\text{CN}(Y, G, \theta; R), \partial_N(Y, G, \theta; R))$, where the chain group $\text{CN}(Y, G, \theta; R)$ is a graded free $R$-module generated by lifts of zeros of $\theta$ in $\tilde{Y}$, and the boundary map $\partial_N(Y, G, \theta)$ is defined by counting flow lines of the vector field dual to $-\theta$. We shall call the homology of this complex the *Novikov homology*, denoted $\text{HN}(Y, G, \theta; R)$, or simply $\text{HN}(Y, \theta; R)$ when $G = H_1(Y; \mathbb{Z})$.

## 2.1.2 The flow on loop space and the Floer complex.

We first outline the basic ingredients for the definition of a Floer complex. The precise definitions, constructions, and proofs of the desired properties of the relevant moduli spaces are deferred to sections 3 and 4.

A Floer theory is modeled on Morse theory on an infinite dimensional loop space $\Omega$. A (path of) symplectic vector field(s) $X$ on the symplectic manifold $(M, \omega)$ naturally defines a closed 1-form $\mathcal{Y}_X$ on $\Omega$, which might not be exact. Let $\tilde{\Omega}$ denote the universal abelian covering of $\Omega$, namely the regular covering with covering transformation group $H_1(\Omega; \mathbb{Z})$. $\mathcal{Y}_X$ lifts to an exact 1-form, $d\tilde{\mathcal{A}}_X$, over $\tilde{\Omega}$. We call $\mathcal{Y}_X$ the *action 1-form*, and $\tilde{\mathcal{A}}_X$ the *action functional*.

On the other hand, a (path of) almost complex structure(s) $J$ specifies a metric on $\Omega$. For generic $(J, X)$, the vector field dual to $-\mathcal{Y}_X$ with respect to the above metric generates a (formal) flow on $\Omega$ (and hence on $\tilde{\Omega}$) satisfying the following properties:

**Properties.** (a) *The set of critical points (i.e. zeros of $\mathcal{Y}_X$), denoted $\mathcal{P}(X)$, consists of finitely many nondegenerate points.*

*Let $\tilde{\mathcal{P}}(X) \subset \tilde{\Omega}$ denote the lift of $\mathcal{P}(X)$. There is a map (called* the index map*)*

$$\mathrm{ind} : \tilde{\mathcal{P}}(X) \to \mathbb{Z}$$

*and a homomorphism (called* the SF-homomorphism*)*

$$\psi : H_1(\Omega; \mathbb{Z}) \to 2\mathbb{Z}$$

*satisfying the following: Let $(x, [w])$, $(x, [w']) = A \cdot (x, [w])$ be different lifts of the same $x \in \mathcal{P}(X)$, where $A \in H_1(\Omega; \mathbb{Z})$ acts by deck transformation. Then*

$$\mathrm{ind}(x, [w']) - \mathrm{ind}(x, [w]) = \psi(A). \tag{3}$$

(b) *The set of flow lines beginning and ending at $(x, [w]), (y, [v]) \in \tilde{\mathcal{P}}(X)$ respectively, denoted $\mathcal{M}_P((x, [w]), (y, [v]))$, is an oriented smooth manifold of dimension $\mathrm{ind}(x, [w]) - \mathrm{ind}(y, [v])$ with a free $\mathbb{R}$ action (translation). Furthermore, for any real constant $\Re$,*

$$\coprod_{(y, [v]) \in S^1_\Re(x, [w])} \mathcal{M}_P((x, [w]), (y, [v]))/\mathbb{R} \tag{4}$$

*consists of finitely many smooth points, where*

$$S^1_\Re(x, [w]) := \Big\{ (y, [v]) \,\Big|\, (y, [v]) \in \tilde{\mathcal{P}}(X),$$
$$\mathrm{ind}(y, [v]) = \mathrm{ind}(x, [w]) - 1, \tilde{\mathcal{A}}_X(y, [v]) - \tilde{\mathcal{A}}_X(x, [w]) > \Re \Big\}.$$

*For any pair $(x, [w]), (y, [v])$ with $\mathrm{ind}(x, [w]) - \mathrm{ind}(y, [v]) = 2$, $\mathcal{M}_P((x, [w]), (y, [v]))/\mathbb{R}$ has a compactification $\overline{\mathcal{M}}_P((x, [w]), (y, [v]))/\mathbb{R}$, which is a compact oriented 1-manifold* with boundary,

$$\begin{aligned}
\partial \overline{\mathcal{M}_P((x,[w]),(y,[v]))/\mathbb{R}} \\
= \overline{\mathcal{M}_P((x,[w]),(y,[v]))/\mathbb{R}} - \mathcal{M}_P((x,[w]),(y,[v]))/\mathbb{R} \\
= \coprod_{(z,[r])\in S^1(x,[w])} \mathcal{M}_P((x,[w]),(z,[r]))/\mathbb{R} \times \mathcal{M}_P((z,[r]),(y,[v]))/\mathbb{R},
\end{aligned}$$

where

$$S^1(x,[w]) := \left\{ (z,[r]) \,\Big|\, (z,[r]) \in \tilde{\mathcal{P}}(X),\, \mathrm{ind}(z,[r]) = \mathrm{ind}(x,[w]) - 1 \right\}.$$

**(c)** *The space of closed orbits in the homology class $A \in H_1(\Omega;\mathbb{Z})$, denoted $\mathcal{M}_O(A)$, is a compact oriented manifold of dimension $\psi(A) + 1$ with a semi-free $S^1$ action (translation). Moreover, for any real constant $\Re$,*

$$\coprod_{A\in \ker\psi \,\mathrm{s.t.}\,-\mathcal{Y}_X(A)<\Re} \mathcal{M}_O(A) \tag{5}$$

*is compact.*

**(d)** *In the $S^1$-equivariant version of Floer theory described in section 4, for $X$ with "$H^2$-induced" flux the homomorphisms $\psi : H_1(\Omega;\mathbb{Z}) \to 2\mathbb{Z}$, $\mathcal{Y}_X : H_1(\Omega;\mathbb{Z}) \to \mathbb{R}$ both factors through a homomorphism*

$$\underline{\mathrm{im}} : H_1(\Omega;\mathbb{Z}) \to H_2(M;\mathbb{Z}).$$

*Furthermore, $\mathcal{P}(X) = \emptyset$ when $X$ is sufficiently small.*

In the above, we follow the convention of calling a point in a moduli space *smooth* or *nondegenerate* if the deformation operator at that point is surjective. A moduli space is said to be *smooth* if it consists of smooth points.

A pair $(J, X)$ for which the smoothness and compactness properties described in (a), (b), (c) above hold is said to be *regular*.

We may now describe the construction of the twisted Floer complex in this context, following the framework of Morse-Novikov theory of closed 1-forms. Start by choosing an appropriate regular covering of $\Omega$, which is typically taken to be the universal abelian covering $\tilde{\Omega}$, namely the covering transformation group $G = H_1(\Omega;\mathbb{Z})$. However in the $S^1$-equivariant version, because of Property (d) above we may work with a smaller covering, with

$$G = \mathrm{Image}(\underline{\mathrm{im}}) \subset H_2(M;\mathbb{Z})$$

in this case. To simplify notation, we shall provisionally denote the factored homomorphisms in (d) by the same notations $\psi : G \to \mathbb{Z}$, $\mathcal{Y}_X : G \to \mathbb{R}$.

Due to the phenomenon of spectral flow in the infinite dimensional setting (i.e. $\psi \neq 0$), the coefficient ring of the (twisted) Floer complex involves *not* the covering transformation group $G$, but the subgroup $\ker \psi \subset G$:

$$\Lambda_F := \text{Nov}(\ker \psi, -[\mathcal{Y}_X]; \mathbb{Z}).$$

Let the chain groups $\tilde{\text{CF}}$ be the free $\Lambda_F$-module generated by elements in $\tilde{\mathcal{P}}(X)$, on which $\ker \psi \subset G$ acts by deck-transformation.

Let the boundary map $\tilde{\partial}_F : \tilde{\text{CF}} \to \tilde{\text{CF}}$ be defined by

$$\tilde{\partial}_F(x, [w]) = \sum_{(y,[v]) \in S^1(x,[w])} \chi(\mathcal{M}_P((x,[w]),(y,[v]))/\mathbb{R}) \, (y,[v]),$$

where the Euler number $\chi$ is just the signed count of $\mathcal{M}_P/\mathbb{R}$ in this case. By Property (b) above, $\tilde{\partial}_F$ is a well-defined $\Lambda_F$-linear transformation, and the last statement of Property (b) ensures that $\tilde{\partial}_F^2 = 0$.

According to Property (a) above, $\tilde{\text{CF}}$ is $\mathbb{Z}$-graded by ind;

$$(\tilde{\text{CF}}, \tilde{\partial}_F) = \bigoplus_{k \in \mathbb{Z}} (\tilde{\text{CF}}^k, \tilde{\partial}_F^k),$$

where each summand $\tilde{\text{CF}}^k$ is of finite rank, and satisfies a periodicity condition

$$(\tilde{\text{CF}}^k, \tilde{\partial}_F^k) = (\tilde{\text{CF}}^{k+2N}, \tilde{\partial}_F^{k+2N}),$$

$2N$ being the gcd of the values of $\psi$. We may thus introduce a reduced version of Floer complex, which is $\mathbb{Z}/2\mathbb{Z}$-graded and is of finite rank:

$$(\text{CF}, \partial_F) = \bigoplus_{k \in \mathbb{Z}/2\mathbb{Z}} (\text{CF}, \partial_F); \quad \text{CF}_k = \bigoplus_{i=1}^{N} \tilde{\text{CF}}_{k'+2i},$$

where $k' \equiv k \mod 2$ is an integer.

### 2.1.3 The Reidemeister torsion.

We now specify the version of torsion used in this paper.

Suppose first for simplicity that the coefficient ring $F$ is a field. Let $(C_i, \partial_i)$, $i \in \mathbb{Z}/2\mathbb{Z}$ be a complex of finite dimensional $F$-vector spaces. The standard short exact sequences $0 \to Z_i \to C_i \to B_{i-1} \to 0$ and $0 \to B_i \to Z_i \to H_i \to 0$ induce a canonical isomorphism

$$\mathcal{T} : \bigotimes_i \det(C_i)^{(-1)^i} \longrightarrow \bigotimes_i \det(H_i)^{(-1)^i}.$$

Let $e$ be an ordered basis for $C_*$, i.e. an ordered basis $e_i$ for each $C_i$. Let $h$ be an ordered basis for $H_*$. Let $[e] \in \bigotimes_i \det(C_i)^{(-1)^i}$ and $[h] \in \bigotimes_i \det(H_i)^{(-1)^i}$ denote the resulting volume forms.

In this simplest case, when the coefficient ring is a field, define the Reidemeister torsion

$$\tau(C_*; e) := \begin{cases} \mathcal{T}([e])/[h] \in F^\times & \text{if } H_* = 0, \\ 0 & \text{otherwise.} \end{cases}$$

For our applications, the coefficient ring of the complex is a Novikov ring of the type specified in (1) (including group rings). We saw that in this case it is in general not a field, but its total ring of fractions is a direct sum of fields. Following Turaev, we define:

**Definition.** [Tu2] Let $R$ be a ring, and assume that its total ring of fractions $Q(R)$ is a finite sum of fields, $Q(R) = \bigoplus_j F_j$. Let $(C_i, \partial_i)$, $i \in \mathbb{Z}/2\mathbb{Z}$ be a complex of finitely generated free $R$-modules with an ordered basis $e$. Then

$$\tau(C_*, e) := \sum_j \tau(C_* \otimes_R F_j, e \otimes 1) \in \bigoplus_j F_j = Q(R).$$

The *Reidemeister torsion of the Floer complex* is defined as

$$\tau_F := \tau(CF, e_X) \in Q(\Lambda_F)/(\pm \ker \psi), \tag{6}$$

where $e_X$ is an ordered basis of $CF$ given by lifts of critical points in $\mathcal{P}(X)$. Different ordering of these elements results in a possible change of sign for $\tau(CF, e_X)$, and different lifts result in a multiplication of $\tau$ by an element in $\ker \psi$; so by modding out $\pm \ker \psi$ in the definition we obtain an invariant independent of these choices.

Another version of torsion that is important in topology is the *Reidemeister torsion of a manifold*, denoted $\tau(Y)$. Let $Y$ be a manifold with a cell-decomposition, such that the cell chain complex $C_*(Y)$ is a finite complex of finite-rank $\mathbb{Z}$-modules. The universal abelian covering $\tilde{Y}$ is endowed with an induced equivariant cell-decomposition, and $C_*(\tilde{Y})$ is a $\mathbb{Z}[H_1(Y; \mathbb{Z})]$-module.

$$\tau(Y) := \tau(C(\tilde{Y}), e_Y) \in Q(\mathbb{Z}[H_1(Y; \mathbb{Z})])/ \pm H_1(Y; \mathbb{Z}),$$

where $e_Y$ is an ordered basis consisting of lifts of cells in $Y$.

### 2.1.4 The zeta function and the counting invariant $I_F$.

Imitating the definition of the dynamical zeta function, we define the the Floer-theoretic *zeta function* as

$$\zeta_F := \exp\left(\sum_{A \in \ker \psi} \chi(\mathcal{M}_O(A)/S^1) A\right)$$

$$= \exp\left(\sum_{A \in \ker \psi} \sum_{u \in \mathcal{M}_O(A)/S^1} \frac{\text{sign}(u)}{m(u)} A\right) \tag{7}$$

$$\in \text{Nov}^1(\ker \psi, -[\mathcal{Y}_X]; \mathbb{Q}) \subset Q(\text{Nov}(\ker \psi, -[\mathcal{Y}_X]; \mathbb{Q})).$$

In the above, $\chi$ is the "orbifold Euler number", $m(u)$ is the "multiplicity", i.e. the order of the stabilizer at $u$. By Property (c), the sum in the exponent in the second expression is in $\mathrm{Nov}^+(\ker\psi, -[\mathcal{Y}_X]; \mathbb{Q})$, and the exponential is well-defined.

Finally, viewing both $\zeta_F$ and $\tau_F$ as elements in $Q(\mathrm{Nov}(\ker\psi, -[\mathcal{Y}_X]; \mathbb{Q}))/\pm \ker\psi$, $I_F$ is simply defined as the product

$$I_F := \zeta_F \tau_F \in Q(\mathrm{Nov}(\ker\psi, -[\mathcal{Y}_X]; \mathbb{Q}))/\pm \ker\psi. \tag{8}$$

**Remark.** Ideally, $I_F$ should be defined in $Q(\Lambda_F)/\pm\ker\psi$ instead of the above larger monoid. However, this would require proving a product formula similar to [HL2] equation (2) to ensure that $\zeta_F \in \mathrm{Nov}^1(\ker\psi, -\mathcal{Y}_X; \mathbb{Z})$.

Of course, the definitions of $\tau_F, \zeta_F, I_F$ all depend on $(J, X)$, and we shall specify this dependence in parenthesis when necessary.

Note that by the last statement in Property (d), in the $S^1$-equivariant version $\tau_F(J, X) = 1$ when $X$ is small. Thus, in this case we may write $I_F(J, X)$ as an element in $\mathrm{Nov}^1(\ker\psi, -[\mathcal{Y}_X]; \mathbb{Q})$, and $\ln I_F(J, X)$ is defined.

## 2.2 Foundational results on $I_F$: the lagrangian intersection version

Let $(M, \omega)$ be a complete symplectic manifold of dimension $2n$, and let $L \subset M$ be a lagrangian submanifold.

In this paper we always assume the lagrangian submanifold $L$ to be compact and oriented.

We first recall some terminologies.

The symplectic form $\omega$ and the Maslov index define two homomorphisms

$$N_\omega^L, N_\mu^L : \pi_2(M, L) \to \mathbb{R}, \mathbb{Z} \quad \text{respectively, as follows.}$$

Let $u : (D, \partial D) \to (M, L)$ be a representative of $[u] \in \pi_2(M, L)$, then $N_\omega^L(u) = \int_D u^*\omega$. With respect to a fixed trivialization of the $\mathrm{Sp}(2n, \mathbb{R})$-bundle $u^*TM$ over $D$, $(u|_{\partial D})^*TL \subset (u|_{\partial D})^*TM$ gives a loop of lagrangian planes in $\mathbb{R}^{2n}$. $N_\mu^L$ is the Maslov index of this loop.

A lagrangian submanifold $L \subset M$ is said to be *monotone* if

$$N_\omega^L = \alpha N_\mu^L \quad \text{for some } \alpha \geq 0.$$

The *minimal Maslov number* of $L$, denoted $\mathbb{N}(L)$, is the gcd of the values of $N_\mu^L$.

Following [FuOOO], we say that $L \subset M$ is *relatively spin* if $L$ is orientable and $w_2(L) \in \mathrm{Image}(\iota_L^*)$, where $\iota_L^*$ is the natural map $H^2(M; \mathbb{Z}_2) \to H^2(L; \mathbb{Z}_2)$ induced by the embedding

$$\iota_L : L \hookrightarrow M.$$

A pair of lagrangian submanifolds $L, L' \subset M$ are *relatively spin* if $L, L'$ are both relatively spin and there is a common $\mathrm{st} \in H^2(M; \mathbb{Z}_2)$ such that $w_2(L) = \iota_L^*(\mathrm{st})$, $w_2(L') = \iota_{L'}^*(\mathrm{st})$.

**2.2.1 Definition.** A lagrangian submanifold $L \subset M$ is *admissible* if $L$ is monotone, $\mathbb{N}(L) \geq 3$, and relatively spin.

**Basic Example.** Let $Y$ be an orientable, closed manifold. Let $T^*Y$ be endowed with the standard exact symplectic form. Then a lagrangian section of $T^*Y$ is admissible. In this paper we shall always assume the standard symplectic structure on $T^*Y$.

An almost complex structure $J_0$ on $M$ is *admissible* if it is compatible with $\omega$ and $(M, \omega, J_0)$ is geometrically bounded, which means $\omega$ tames $J_0$ uniformly, and the sectional curvature of $M$ is bounded above, while the injectivity radius is bounded below. (See [ALP] Definition 2.2.1).

Given a vector bundle $F$ over $M$, let $C_\epsilon = C_\epsilon(M; F) \subset C^\infty(M; F)$ be the Banach space of sections defined by Floer [F88a]. Let $\mathcal{J}$ be the space of $C_\epsilon$-paths $J_t, t \in [0,1])$ of admissible almost complex structures on $M$ which are $t$-independent outside a compact set in $M$. Let

$$\mathcal{H} = \left\{ H \in C_\epsilon([0,1] \times M; \mathbb{R}) \,\middle|\, dH_t := dH(t, \cdot) \text{ is compactly supported} \right\}.$$

(These spaces are in general not complete, but that doesn't matter for our purpose, cf. §3.2.2, last paragraph). For an $H \in \mathcal{H}$, let $\chi_H$ denote the path of hamiltonian vector fields $\chi_{H_t}$, $\omega(\chi_{H_t}, \cdot) = dH_t$.

Similarly, let $\mathcal{X}$ be the space of $C_\epsilon$-paths of compactly supported symplectic vector fields over $M$.

*Notation.* When $H$ is a function or section over $[0,1] \times M$, we often regard it as a path of functions or sections over $M$, and denote $H_t := H(t, \cdot)$.

We denote by $\mathrm{Symp}(M)$ the group of all $C_\epsilon$-symplectomorphisms of $M$, and by $\mathrm{Symp}_c^0(M)$ the group of symplectomorphisms generated by vector fields in $\mathcal{X}$. A *symplectic isotopy* is a path $\phi_t \in \mathrm{Symp}(M)$, $t \in [0,1]$. The *flux* of the symplectic isotopy $\phi_t$ is $[\int_0^1 \iota(d\phi_t/dt)\omega \, dt] \in H^1(M)$.

*Notation.* There is a one-to-one correspondence between symplectic vector fields and closed 1-forms on $M$. In this paper, $\theta_X$ denotes $\iota(X)\omega$; conversely, $X_\theta$ denotes the symplectic vector field corresponding to the 1-form $\theta$.

Given a pair of admissible lagrangian submanifolds $L, L' \subset M$, let $\Omega(M; L, L')$ be the space of $L_1^p$ paths of unit length from $L$ to $L'$. Let $\Omega^{\gamma_0}(M; L, L') \subset \Omega(M; L, L')$ be the path component containing $\gamma_0$. Each $[\gamma_0] \in \pi_0(\Omega(M; L, L'))$ defines a homomorphism $h_{[\gamma_0]} : \pi_1(L) \to \pi_1(M, L')$, sending $a \in \pi_1(L)$ to the concatenation $[\gamma_0^{-1} * a * \gamma_0] \in \pi_1(M, L')$.

For convenience, we assume that

$$h_{[\gamma_0]} = 0. \tag{9}$$

In section 3, we lay the foundation for constructing a Floer theory satisfying the properties listed in §2.1.2, based on a formal flow on $\Omega = \Omega^{\gamma_0}(M; L, L')$ generated by generic $(J, X) \in \mathcal{J} \times \mathcal{X}$. A Floer-theoretic torsion is thereby defined via the recipe in §2.1; we denote this version of $I_F$ by $I_F^{\gamma_0}(L, L'; J, X)$. Note that the cohomology class $[\mathcal{Y}_X]$ depends only on the cohomology class $[\theta_X]$ (in addition to $[\omega]$ and the loop space $\Omega$ itself); see (22). Thus so does the relevant Novikov ring $\Lambda_F = \Lambda_F^{\gamma_0}(L, L'; X)$ and the monoid that $I_F^{\gamma_0}(L, L'; J, X)$ takes value in.

$I_F^{\gamma_0}(L, L'; J, X)$ enjoys the following invariance properties.

**2.2.2 Theorem.** *Let $L, L' \subset M$ be a pair of admissible lagrangian submanifolds with a fixed relative spin structure ([FuOOO] Definition 21.5), and let $[\gamma_0]$ satisfy (9) above.*

(a) *Let $X \in \mathcal{X}$ and let $(J_1, H_1), (J_2, H_2) \in \mathcal{J} \times \mathcal{H}$ be such that $(J_1, X + \chi_{H_1})$, $(J_2, X + \chi_{H_2})$ are regular. Then*

$$I_F^{\gamma_0}(L, L'; J_1, X + \chi_{H_1}) = I_F^{\gamma_0}(L, L'; J_2, X + \chi_{H_2}). \tag{10}$$

(b) *If moreover $(\iota_L)_* : H_1(L) \to H_1(M)$ is injective, then there is a $I_F^{\gamma_0}(L, L') \in Q(\mathbb{Z}[\ker \psi])/(\pm \ker \psi)$ such that*

$$I_F^{\gamma_0}(L, L'; J, X) = i_{-[\mathcal{Y}_X]} I_F^{\gamma_0}(L, L') \quad \text{for all regular } (J, X). \tag{11}$$

When (10) above holds, we say that $I_F$ is *invariant under hamiltonian isotopies*, whereas when (11) holds, we say that $I_F$ is *invariant under symplectic isotopies*.

These terminologies reflect a well-known principle that equates perturbations with symplectic isotopies, which we review in §3.1.3 below. According to this principle, if $\Phi = \{\phi_t \,|\, t \in [0, 1], \phi_0 = \text{id}\}$ is a symplectic isotopy, then it induces a canonical isomorphism

$$\Phi : \Omega^{\gamma_0}(M; L, L') \to \Omega^{\Phi \cdot \gamma_0}(M; L, \phi_1(L')),$$

and for any regular $(J, X)$, there is another regular $(J', X')$ such that $\Phi^* \mathcal{Y}_{X'} = \mathcal{Y}_X$ (hence $(\Phi)_* \Lambda_F^{\gamma_0}(L, L'; X) = \Lambda_F^{\Phi \cdot \gamma_0}(L, \phi_1(L'); X'))$, and

$$(\Phi)_* I_F^{\gamma_0}(L, L'; J, X) = I_F^{\Phi \cdot \gamma_0}(L, \phi_1(L'); J', X').$$

Thus, via this principle the above theorem implies that $I_F$ is actually also invariant under compactly supported hamiltonian isotopies of one of $L, L'$, or more generally, symplectic isotopies under the stronger conditions of statement (b).

In fact, (11) always holds when $L' = L$. However, it is in general not possible to reduce the computation of $I_F^{\gamma_0}(L, L')$ to a neighborhood of $L$, due to the presence of "quantum corrections" from "large" pseudo-holomorphic curves (i.e. those that goes out of a neighborhood of $L$). With the additional assumption that $\pi_2(M, L) = 1$, the possibility of quantum corrections is eliminated, and one has:

**2.2.3 Theorem.** *Let $L \subset M$ be an admissible lagrangian submanifold, and let $\phi = \phi_1 \in \text{Symp}_c^0(M)$, which is connected to $\phi_0 = \text{Id}$ via a symplectic isotopy $\Phi = \{\phi_t \,|\, \phi_t \in$* $\mathrm{Symp}_c^0(M)$, $t \in [0, 1]\}$. Let $\gamma_\phi(t) := \phi_t(p)$, where $p$ is a base point in $M$, and let $\gamma_p(t) = p$ be the constant path. Then

$$I_F^{\gamma_\phi}(L, \phi(L)) = (\Phi)_* I_F^{\gamma_p}(L, L)$$

is invariant under symplectic isotopies.

If in addition $\pi_2(M, L) = 1$, then there is an isomorphism

$$e_{0*} : \ker \psi = H_1(\Omega^{\gamma_\phi}(L, \phi(L)); \mathbb{Z}) \to H_1(L; \mathbb{Z})$$

(see (19)) which induces another isomorphism

$$e_{0*} : Q(\mathbb{Z}[\ker \psi])/ \pm \ker \psi \to Q(\mathbb{Z}[H_1(L; \mathbb{Z})])/ \pm H_1(L; \mathbb{Z}),$$

under which

$$e_{0*} I_F^{\gamma_\phi}(L, \phi(L)) = \tau(L).$$

**Examples.** (a) Let $M = \Sigma_2$, a genus 2 surface, and $L$ be the circle in the middle separating $\Sigma_2$ into two 1-handles. Then the above theorem implies in this case $I_F^{\gamma_p}(L, \phi(L)) = (1-t)^{-1}$, where $t$ is a generator of $H_1(S^1)$.
(b) Let $Y$ be *any* oriented, compact manifold; let $M = T^*Y$, and $L = Y$ be the zero section. Then the last theorem implies that

$$I_F^{\gamma_p}(L, \phi(L)) = \tau(Y).$$

Since $Y$ can be quite arbitrary, we obtain in this way many examples of interesting $I_F$. The Reidemeister torsion $\tau(Y)$ can be computed in many important cases. For instance, if $Y = \Sigma_f$ is a *mapping torus* of a diffeomorphism $f : \Sigma \to \Sigma$ of a closed manifold $\Sigma$,

$$\Sigma_f := \Sigma \times [0, 1]/(x, 0) \sim (f(x), 1), \tag{12}$$

then $\tau(\Sigma_f)$ may be computed from the (twisted) Lefschetz numbers of $f$ and its iterates. In the simplest case, when $H_1(\Sigma_f) = \mathbb{Z}$ of which $t$ is a generator,

$$\tau(\Sigma_f) = \zeta(f) := \exp\left(\frac{\#\operatorname{Fix}(f^k)}{k} t^k\right). \text{ (See [Mil])}.$$

On the other hand, if $Y$ is a 3-manifold, there are surgery formulae to compute $\tau(Y)$ [Tu]. For example, if $Y = S_0^3(K)$ is the 0-surgery of a knot $K \subset S^3$, then

$$\tau(Y) = \operatorname{Alex}(K)/(1-t)^2,$$

$\operatorname{Alex}(K)$ being the Alexander polynomial of $K$. It is well-known that the Alexander polynomial can be any A-polynomial (i.e. polynomials $P(t)$ with symmetric coefficients and $P(1) = 1$).

An immediate corollary of the computation in Theorem 2.2.3 is:

**2.2.4 Corollary.** *Let $Y$ be a closed orientable manifold with $\mathrm{lt}(i_\theta \tau(Y)) \neq 1$, where $\theta: H_1(Y;\mathbb{Z}) \to \mathbb{R}$ is a cohomology class of $Y$. Let $L \in T^*Y$ be the zero section. Then it is impossible to disengage $L$ from itself via a symplectic isotopy of flux $\theta \in H^1(T^*Y) = H^1(Y)$.*

**Example.** $Y = S_0^3(K)$, where $K$ is the $k$-twisted knot with $k > 1$. (See Figure 1). Then $\mathrm{lt}(i_\theta \tau(Y)) = k$ for any $\theta \neq 0$. Thus it is impossible to disengage the zero section of $T^*Y$ from itself via any symplectic isotopy. More generally, $K$ can be any knot whose Alexander polynomial has leading coefficient $\neq 1$. In particular, $K$ is necessarily not fibered. If $K$ is fibered, the zero section of $T^*S_0^3(K)$ *can* be disengaged from itself by symplectic isotopies, since $S_0^3(K)$ carries a nowhere vanishing closed 1-form.

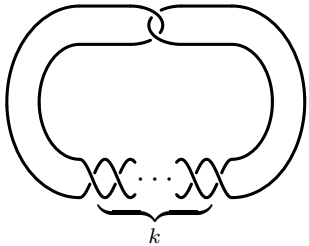

FIGURE 1: the $k$-twisted knot.
Alexander polynomial $= k - (2k+1)t + kt^2$

## 2.3 Foundational results on $I_F$: the $S^1$-equivariant version

In this version, let $(M, \omega)$ be a compact monotone symplectic manifold of dimension $\geq 4$.

*Notation.* We shall use the same notation to denote analogous but often different objects in the two versions of Floer theory (e.g. the conditions on $(M, \omega)$ are different); what it actually means should be clear from the context.

Let $\gamma_0: S_1^1 \to M$ be a loop in $M$, and let $\mathcal{C} = \Omega(M; \gamma_0)$ be the path component containing $\gamma_0$ in the Banach manifold of $L_1^p$-loops of unit-length in $M$.

For simplicity, we assume that the homotopy class of $\gamma_0$ is in the center of $\pi_1(M)$, and that its homology class $[\gamma_0]$ is primitive. A homology class $g \in H_*(M;\mathbb{Z})$ is said to be *primitive* if $g \neq ng'$ for any $n \in \mathbb{Z}\setminus\{\pm 1\}$ and $g' \in H_*(M;\mathbb{Z})$. (In particular, $g$ is nontorsion).

There is a $S^1$ action on $\mathcal{C}$ by rotation:

$$(R_q \gamma)(t) := \gamma(t+q); \quad q \in S^1; \gamma \in \mathcal{C}.$$

This action is free when $[\gamma_0]$ is primitive. Let
$$\underline{\mathcal{C}} = \underline{\Omega(M; \gamma_0)} := \mathcal{C}/S^1.$$

For comparison with the Gromov invariants, it is useful to regard loops in $\underline{\mathcal{C}}$ as tori in $M$: A loop in $\underline{\mathcal{C}}$ corresponds to a "twisted loop" in $\mathcal{C}$, which may be written as a $M$-valued function $u(s,t)$, where $t \in S^1$ parameterizes the element in $\mathcal{C}$, $s \in [0,1]$ parameterizes the twisted loop, and $u(1, \cdot) = R_q u(0, \cdot)$ for some $q$. In other words,

$$u \text{ maps } [0, 1] \times S^1 / \{(0, t+q) \sim (1, t)\} = \mathbb{T}^2 \text{ to } M. \tag{13}$$

It is easy to see that this defines a homomorphism

$$\underline{\text{im}} : H_1(\underline{\mathcal{C}}; \mathbb{Z}) \to H_2(M; \mathbb{Z}).$$

Similarly, one may define a homomorphism $\text{im} : H_1(\mathcal{C}; \mathbb{Z}) \to H_2(M; \mathbb{Z})$ by setting $q = 0$ in (13).

As mentioned in §2.1, in this version of Floer theory, we let $G = \text{Image}(\underline{\text{im}})$. We denote
$$\mathfrak{H} := \ker \psi \Big|_G \subset H_2(M; \mathbb{Z}).$$

Let $\mathcal{J}$ be the Banach manifold of ($t$-independent) $C_\epsilon$ almost complex structures on $M$ that are compatible with $\omega$. In Definition 4.1.2 we introduce the notion of "$H^2$-induced classes" in $H^1(M; \mathbb{R})$ and a homomorphism

$$\mathfrak{h} : \{H^2\text{-induced classes}\} \to \text{Hom}(\mathfrak{H}, \mathbb{R}).$$

Let $\mathcal{X}$ be the Banach space of ($t$-independent) $C_\epsilon$ symplectic vector fields on $M$ with $H^2$-induced flux $[\theta_X]$.

In section 4, we establish the foundation of an $S^1$-equivariant version of Floer theory modeling on a formal flow generated by generic $(J, X) \in \mathcal{J} \times \mathcal{X}$ on

$$\Omega = \underline{\mathcal{C}} = \underline{\Omega(M; \gamma_0)},$$

and verify that it satisfies the properties listed in §2.1.2. The version of $I_F$ thus obtained (via §2.1) is denoted $I_F^{S^1}(M, \gamma_0; J, X)$.

**2.3.1 Theorem.** *Let $(M, \omega)$ be a compact monotone symplectic manifold of dimension $\geq 4$. Suppose $\gamma_0$ represents a central homotopy class and a primitive homology class as above. Then $I_F^{S^1}$ is invariant under symplectic isotopies in the sense of (11), namely, there exists $I_F^{S^1}(M, \gamma_0) \in Q(\mathbb{Z}[\mathfrak{H}])/(\pm \mathfrak{H})$ such that*

$$I_F^{S^1}(M, \gamma_0; J, X) = i_{-\underline{\text{im}}^*[\mathcal{Y}_X]} I_F^{S^1}(M, \gamma_0) \quad \text{for any regular } (J, X).$$

In the above, $-\underline{\text{im}}^*[\mathcal{Y}_X]\Big|_{\mathfrak{H}} = [\omega]\Big|_{\mathfrak{H}} - \mathfrak{h}[\theta_X]$.

**2.3.2 Basic example: the symplectic mapping torus.** Here is an example satisfying all the above assumptions. Let $(\Sigma, \omega_\Sigma)$ be a compact monotone symplectic manifold of dimension $\geq 2$. (For example, $\Sigma$ can be any orientable surface). Let $f : \Sigma \to \Sigma$ be a symplectomorphism and let $\Sigma_f$ be the mapping torus of $f$ (Cf. (12)). The projection $\pi_S : \Sigma_f \to S^1$ defines a closed 1-form $q := d\pi_S$; $[q] \in H^1(\Sigma_f; \mathbb{Z})$ is nontrivial. Let
$$M = M(\Sigma, f) := \Sigma_f \times S^1_\tau,$$
and endow it with the standard symplectic form
$$\omega_{M(\Sigma, f)} = \omega_\Sigma + q \wedge d\tau, \tag{14}$$
where $\tau \in S^1_1$ parameterizes the unit circle $S^1_\tau$, and $\omega_\Sigma$ makes sense as a closed 2-form on $M$, since $f$ preserves the symplectic form on $\Sigma$.

Let $\gamma_0 = \{p\} \times S^1 \subset M$.

We shall explain in §6.3 that the first summand in the decomposition
$$H^1(M; \mathbb{R}) = H^1(\Sigma_f; \mathbb{R}) \oplus H^1(S^1; \mathbb{R})$$
is the space of $H^2$-induced classes. Thus $X$ can be any symplectic vector field with $[\theta_X] \in H^1(\Sigma_f; \mathbb{R}) \subset H^1(M; \mathbb{R})$. On the other hand,
$$\mathfrak{H} = H_1(\Sigma_f; \mathbb{Z}) \otimes H_1(S^1; \mathbb{Z}) \oplus \ker c_1\big|_{\pi_2(\Sigma)} \subset H_2(M; \mathbb{Z}), \tag{15}$$
and the homomorphism $\mathfrak{h} : H^1(\Sigma_f; \mathbb{R}) \to \mathrm{Hom}(\mathfrak{H}, \mathbb{R})$ is such that for $\beta \in H^1(\Sigma_f; \mathbb{R})$, the function $\mathfrak{h}\beta$ is defined by
$$\begin{cases} (\mathfrak{h}\beta)(a \otimes [\gamma_0]) = \beta(a) & \text{for all } a \in H_1(\Sigma_f; \mathbb{Z}); \\ (\mathfrak{h}\beta)(b) = 0 & \text{for all } b \in \ker c_1\big|_{\pi_2(\Sigma)}. \end{cases}$$

Let $\mathrm{Gr}_{A,g,n}(M, \beta_1, \cdots, \beta_n)$ denote the Gromov invariants defined by Ruan-Tian, which counts pseudo-holomorphic maps of genus $g$ curves of $n$ marked points of the class $A \in H_2(M)$, sending the $n$ marked points to representatives of $\beta_1, \ldots \beta_n \in H_{2n-2}(M)$ respectively [RT2]. Usually, $g, n$ are required to be in a "stable range", but though $1, 0$ is not in the stable range, $\mathrm{Gr}_{A,1,0}(M)$ can be defined via stabilization [IP2],
$$\mathrm{Gr}_{A,1,0}(M) := \mathrm{Gr}_{A,1,1}(M, \beta)/A \cdot \beta \quad \text{for any } \beta \in H_{2n-2}(M). \tag{16}$$
Recall now the remarks in the last paragraph of §2.1.4, which imply in particular that $\ln i_{[\omega]} I_F^{S^1}(M, \gamma_0)$ is well-defined.

**2.3.3 Theorem.** *Under the conditions of the previous theorem, for any primitive class $A \in \mathfrak{H}$,*
$$(\ln i_{[\omega]} I_F^{S^1}(M, \gamma_0))(A) = \mathrm{Gr}_{A,1,0}(M).$$

**Example.** Continuing Example 2.3.2, According to Ionel-Parker's computation of the Gromov invariants in this case [IP2], when $b \in H_1(\Sigma_f; \mathbb{Z})$ is a primitive class,

$$\mathrm{Gr}_{b \otimes [\gamma_0], 1, 0}(M) = \# \mathrm{Fix}_b(f),$$

where $\# \mathrm{Fix}_b$ is the number of fixed points "of homology class $b$". This is an abelian version of the Nielsen number. By the *homology class* of a fixed point we mean the following: noticing that each fixed point of $f$ corresponds to a loop in $\Sigma_f$; the homology class of this fixed point is the homology class of this loop.

In fact, inspired by Taubes's definition of a generating function of Gromov invariants of tori, Ionel-Parker defined a "Gromov series" which agrees with the zeta function of $f$ in this case. Formally, the Gromov series much resembles our definition of the Floer-theoretic zeta function. It can be shown directly that $I_F^{S^1}$ agrees with the (twisted) zeta function of $f$ in the case of symplectic mapping tori. However, we do not expect the Gromov series to coincide with $I_F^{S^1}$ in general. [Ln]

**Remarks.** The various artificial assumptions on $M$, $\gamma_0$, and $X$ in this version arise as a compromise under the dichotomy between the two main motivations of this work: to compare $I_F$ with Gr, and to apply to existence problems of periodic orbits.

Using $H_1(\underline{\mathcal{C}})$ instead of Image($\underline{\mathrm{im}}$) $\subset H_2(M; \mathbb{Z})$ for $G$, one may define a more-refined version of $I_F^{S^1}(M, \gamma_0; J, X)$, for $X$ with flux $\theta_X$ in *any* cohomology class satisfying $[\theta_X](\gamma_0) = 0$. This version however does not have a simple expression in Gromov invariants.

One may relax or perhaps entirely remove the assumptions on the homotopy class of $\gamma_0$, and still manage to prove the comparison theorem between Gr and this generalized $I_F^{S^1}$, but it would only be invariant under Hamiltonian isotopies.

Monotonicity is used to simplify the proofs for compactness of moduli spaces, and more essentially, in the proof of invariance under (non-Hamiltonian) symplectic isotopies, and in the comparison theorem with Gr.

The assumption that $\dim(M) \geq 4$ is used in the transversality proof.

## 2.4 Noncontractible periodic orbits from Floer-theoretic torsions

In this paper, being $\lambda$-*periodic* means being periodic with period $\lambda \in \mathbb{R}^+$. $S_T^1$ denotes a circle of length $T$. Given a 1-periodic path of symplectic vector fields $X_t$ on $M$, $t \in S_1^1$, a *periodic orbit* of $X_t$ means a solution of

$$\frac{dx}{dt} = X_t(x(t)), \quad x : S_1^1 \to M.$$

A $\lambda$-*periodic orbit* of a symplectic vector field $X$ is a solution of

$$\frac{dx}{dt} = X(x(t)), \quad x : S_\lambda^1 \to M.$$

PO($X$; $b$) denotes the set of periodic orbits of $X$ in the class $b \in H_1(M; \mathbb{Z})$.

We now state some existence results of such periodic orbits in a homology class $b \neq 0$. Theorems A, B, and Corollary C below use the lagrangian intersection version of $I_F$ described in §2.2; Corollary D uses the $S^1$-equivariant version described in §2.3.

The statements in Theorems A, B, and Corollary C below are all lifted from [GL]. In [GL], $Y$ is assumed to be the mapping torus of an involution of a flat manifold. We extend it to the following more general situation:

### 2.4.1 Manifolds of Type F.

**Definition.** Given $g \in H^1(Y)$, a primitive element $b \in H_1(Y; \mathbb{Z})$ is said to be *g-essential* if $g(b) > 0$, and $\ln i_g \tau(Y) : H_1(Y; \mathbb{Z}) \to \mathbb{Q}$ (when defined) has nontrivial value at $b$. (Cf. remarks following the Examples below).

A closed oriented manifold $Y$ is said to be of *Type F* if it admits a nowhere vanishing closed 1-form $\theta$, and a $[\theta]$-essential class $b \in H_1(Y)$. We also call $(Y, \theta, b)$ a *Type F triple*.

**Examples.** For the simplest example, let $Y = S^1$ and $\theta, b$ be generators of $H^1(Y), H_1(Y; \mathbb{Z})$ respectively. A 3-dimensional example is $Y = S_0^3(K)$, where $K$ is a trefoil knot (See Figure 2), and $\theta, b$ are as in the previous case. More generally, $K$ can be any fibered knot with Alexander polynomial of the form $1 + at +$ higher order terms in $t$ (modulo multiplication by $t$ and $t^{-1}$), where $a \neq 0$. (Note that the Alexander polynomial of any fibered knot has leading coefficient 1).

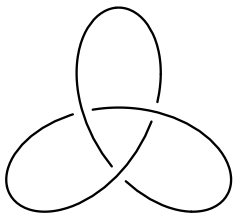

FIGURE 2: the trefoil knot
Alexander polynomial= $1 - t + t^2$

Topologically, a manifold $Y$ admitting a nowhere vanishing closed 1-form fibers over $S^1$. Algebraically, it is easy to see that the Novikov homology $\mathrm{HN}(Y, \theta; \mathbb{Q}) = 0$ and the torsion $i_{-[\theta]} \tau(Y)$ is represented by an element in $\mathrm{Nov}^1(H_1(Y), -\theta; \mathbb{Q})$. Thus, in this case $\ln i_{-[\theta]} \tau(Y)$ is well defined. (It is computable via Lefschetz numbers). In fact, it is possible to strengthen the above observations to obtain an algebraic characterization of manifolds admitting nowhere vanishing closed 1-forms. (This is due to Latour for dimensions $\geq 6$; to Stallings and Thurston for dimension $= 3$).

The condition of Type F is somewhat stronger than the existence of a nowhere vanishing closed 1-form. In particular, a Type F manifold $Y$ must have Reidemeister torsion $\tau(Y) \neq 1$. A torus $T^n$, $n > 1$ is the simplest example of manifolds admitting a nowhere vanishing closed 1-form but is not of Type F.

**Theorem A.** *Let $(Y, \theta, b)$ be a Type F triple, and let $L, L' \subset T^*Y$ be respectively the zero section and the section corresponding to $\theta$. Let $H$ be a 1-periodic hamiltonian on $T^*Y$ separating $L, L'$; namely, $H : S_1^1 \times T^*Y \to \mathbb{R}$; $H_t := H(t, \cdot)$, then*

$$\mu := \int_{S_1^1} \left(\min H_t \Big|_{L'} - \max H_t \Big|_{L}\right) dt > 0.$$

*Furthermore, suppose $\nabla H_t$ is compactly supported $\forall t \in S_1^1$, and let $m := [\theta](b)$. Then there exists a $\lambda \in (0, m/\mu]$ such that there is a periodic orbit of $X_{\lambda H_t}$ in the homology class $b$.*

**Remarks. (a)** Note that if $(Y, \theta, b)$ is a Type F triple, then $(Y, -\theta, -b)$ is another Type F triple due to the Poincaré duality of Reidemeister torsion [Mil]. Thus the above theorem applies to give another periodic orbit in class $-b$.

**(b)** The number $\lambda$ above may be viewed as a generalization of the notion of period. Indeed, when $H$ is $t$-independent, periodic orbits of $X_{\lambda H}$ and $\lambda$-periodic orbits are equivalent.

It is well known that there is a neighborhood of a lagrangian submanifold $L \subset M$ symplectomorphic to the unit-disk bundle $DT^*L = \{(q, p) | q \in L, p \in T_qL, |p| < 1\} \subset T^*L$ (a *Weinstein neighborhood*). The previous theorem thus implies:

**Theorem B.** *Let $(Y, \theta, b)$ be a type F triple. Let $L$ be a lagrangian embedding of $Y$ into $M$, and let $U_L$ be a Weinstein neighborhood of $L$. Let $L' \subset M$ be the section corresponding to $\varepsilon\theta$ in $U_L = DT^*L$, where $\varepsilon \in \mathbb{R}^+$, $\varepsilon \ll 1$. Let $H, \mu$ be as in Theorem A. Then for sufficiently small $\varepsilon$, there exists a $\lambda \in (0, \varepsilon[\theta](b)/\mu]$ such that $X_{\lambda H_t}$ has a periodic orbit in class $b$.*

Recall from [GL] that a lagrangian submanifold $L \subset \mathbb{C}^n$ is not *symplectically self-linked* if it is impossible to separate $L$ from a nearby copy of itself by a Hamiltonian isotopy fixing $L$. To be more precise, there exists a small lagrangian isotopy $L_t$, $t \in [0, 1]$ with $L_0 = L$, $L_t$ disjoint from $L$ $\forall t > 0$, and there is a Hamiltonian isotopy fixing $L$, and moving $L_1$ to either side of any hyperplane in $\mathbb{C}^n$. $L$ is symplectically self-linked otherwise.

**Corollary C.** *If a manifold of Type F lagrangian embeds into $\mathbb{C}^n$, then it must be symplectically self-linked.*

**Remarks.** It is desirable to remove the assumptions of orientability of $Y$ and/or primitivity of $b$ from the Definition of Type F. Note that $Y$ is allowed to be unorientable in [GL]. Oh showed in [Oh97a] that the relevant moduli spaces are oriented for $M = T^*Y$ even for unorientable $Y$; so the Floer theory of lagrangian intersection is well-defined in this case. However, our definition of $I_F(L, L')$ uses the orientability of the lagrangian submanifolds $L, L'$.

On the other hand, if the primitivity of $b$ may be removed, it would be easy to characterize such manifolds algebraically.

Turning now to the $S^1$-equivariant version, here is an almost immediate corollary of Theorem 2.3.3.

**Corollary D.** *Under the assumptions of Theorem 2.3.3, let $A \in \mathfrak{H}$ be a primitive class for which $\mathrm{Gr}_{A,1,0}(M) \neq 0$ and $\omega(A) =: m > 0$. Then for any $X \in \mathcal{X}$, either $\mathrm{PO}(X;[\gamma_0]) \neq \emptyset$, or there exists a $J|X$-torus (i.e. an element of $\mathcal{M}_O$, see (37)). In particular, if $\mathfrak{h}[\theta_X] = \mu/m[\omega]\big|_\mathfrak{H}$ for $\mu \in \mathbb{R}^+$, then $\mathrm{PO}(\lambda X;[\gamma_0]) \neq \emptyset$ for $\lambda = m/\mu$.*

**2.4.2 Example.** Continuing Example 2.3.2 again, Corollary D says in this case that:

Let $b \in H_1(\Sigma_f; \mathbb{Z})$ be a primitive class such that $\#\mathrm{Fix}_b(f) \neq 0$ and $m := [q](b) > 0$. Let $X$ be a symplectic vector field on $M$ with $[\theta_X] = [q]/\lambda$ for $\lambda \in \mathbb{R}^+$. Then there exists a $\lambda$-periodic orbit of $X$ in the homology class of $\{p\} \times S^1 \subset M$.

For instance, this holds when $\Sigma_f = S_0^3(K)$ where $K$ is a trefoil knot, and $b, \eta$ generate $H_1(S_0^3(K); \mathbb{Z}), H^1(S_0^3(K); \mathbb{R})$ respectively.

It is interesting to compare Theorems A, D with the following recent result of Biran-Polterovich-Salamon:

**2.4.3 Theorem.** ([BPS] Theorem B) *Let $\mathbb{T}^n = \mathbb{R}^n/\mathbb{Z}^n$ and $M = T^*\mathbb{T}^n$. Let $H : [0,1] \times M \to \mathbb{R}$ be a path of smooth compactly supported functions on $M$, and $e \in \mathbb{Z}^n$ is such that*

$$|e| \leq \mu := \max\nolimits_{[0,1] \times Z} H, \tag{17}$$

*then $X_{dH_t}$ has a periodic orbit in the homology class $e \in \mathbb{Z}^n \simeq H_1(\mathbb{R}^n/\mathbb{Z}^n) \simeq H_1(M)$, with action $\geq \mu$.*

Their proof makes use of another refinement of the Floer homology—a filtered version (by action) of Floer homology called 'relative symplectic homology'. The calculation of the relative symplectic homology is in general difficult; in the case above, the computation was reduced to certain special $H$ for which the periodic orbits can be understood in terms of closed geodesics in $\mathbb{T}^n$.

We have deliberately chosen analogous notations in the statements of Theorems A, D and Theorem 2.4.3 to emphasize the similarity of these results. Very heuristically, in all three $\mu$ is a sort of 'energy', and $m$ or $|e|$ represent a sort of 'capacity'. Theorem D and Theorem 2.4.3 are more precise than Theorem A, reflecting the advantage of Floer theory over Gromov invariants in finding periodic orbits. When computable (which is often hard), the symplectic homology gives the best result. The comparison of the three theorems also raises the following natural questions:

**Questions.** In analogy with (17), is it possible to strengthen Theorems A and D to say that periodic orbits exist for all $\lambda \geq [m/\mu, \infty)$? Is this in any way related to th energy-capacity inequalities in symplectic topology? (See e.g. [MS2, La]). We know that generically, periodic orbits exist for an interval of $\lambda$. On the other hand, the dimensional assumption on $M$ in Theorem D can not be removed, since it is easy to construct 2-dimensional counterexamples.

# 3 Floer theory of lagrangian intersections

In this section we gather some foundational materials on the Floer theory of lagrangian intersections that are necessary for the definition of $I_F(L, L')$. These are just minor modifications of the existent literature such as the standard references [F88b, Oh93] and [L]; we include them here only for the sake of completeness.

## 3.1 Setup and basics

Here we review of the setup and state some basic facts.

### 3.1.1 The relative loop space $\Omega(M; L, L')$.

Let $\Omega := \Omega_{\gamma_0}(M; L, L')$ be as in §2.2. We shall frequently make use of the following alternative description of the universal abelian covering $\tilde{\Omega}$:

$$\tilde{\Omega} = \Big\{ (\gamma, [w]) \Big| \gamma \in \Omega; w : [0,1] \to \Omega; w(0) = \gamma_0, w(1) = \gamma;$$
$$[w_1] = [w_2] \text{ if } [w_1 - w_2] = 0 \in H_1(\Omega; \mathbb{Z}) \Big\}.$$

Let $\Omega(M; L) := \Omega(M; L, L)$, and let $\Omega_0(M, L)$ denote the component containing the constant paths. Under assumption (9), the map

$$e_0 : \Omega(M, L, L') \to L, \quad \gamma \mapsto \gamma(0) \tag{18}$$

induces a fibration

$$\Omega_0(M, L') \xrightarrow{f_0} \Omega_{\gamma_0}(M; L, L') \xrightarrow{e_0} L.$$

Abelianizing the homotopy sequence of this fibration and using (9) again, we have:

$$\pi_2(M, L') \xrightarrow{f_{0*}} H_1(\Omega; \mathbb{Z}) \xrightarrow{e_{0*}} H_1(L; \mathbb{Z}) \to 0. \tag{19}$$

(We used the fact that $\pi_1(\Omega_0(M, L')) = \pi_2(M, L')$; see [Wh]).

Similar to $N_\omega^L, N_\mu^L$ in §2.2, below we define two homomorphisms

$$\psi_\omega, \psi_\mu : H_1(\Omega; \mathbb{Z}) \to \mathbb{R}, \mathbb{Z} \quad \text{respectively}.$$

Let $u : S^1 \times [0,1] \to M$ be a representative of $[u] \in H_1(\Omega)$, then

$$\psi_\omega(u) := \int_{S^1 \times [0,1]} u^*\omega. \tag{20}$$

This is independent of the choice of representatives since $L, L'$ are lagrangian.

Choose a trivialization of $u^*TM$ over $S^1 \times [0,1]$, Then with respect to the trivialization, $(u|_{S^1 \times \{0\}})^*TL \subset (u|_{S^1 \times \{0\}})^*TM$ and $(u|_{S^1 \times \{1\}})^*TL' \subset (u|_{S^1 \times \{1\}})^*TM$ give

two loops of lagrangian planes in $\mathbb{R}^{2n}$. Let $\mu_0(u), \mu_1(u)$ respectively be the Maslov indices of them, then

$$\psi_\mu(u) := \mu_1(u) - \mu_0(u). \tag{21}$$

It is easy to see that $\psi_\mu$ does not depend on the choice of representatives or trivialization. (Cf. e.g. [Oh93].)

Since our $L$, $L'$ are orientable by assumption, $\psi_\mu$ always has even values. We denote by $\mathbb{N}(L, L')$ the gcd of the values of $\psi_\mu$.

In the lagrangian-intersection version of Floer theory, we set

$$\psi = \psi_\mu,$$

$\psi$ being the SF homomorphism described in §2.1.2.

### 3.1.2 A formal flow on $\Omega$ associated to $(J, X)$.

Given a path of symplectic vector fields $X_t, t \in [0, 1]$ over $M$, the associated action 1-form on $\Omega$ is:

$$\mathcal{Y}_X(\gamma)(\xi) := -\int_0^1 \omega(\partial_t\gamma(t), \xi(t))\, dt + \int_0^1 \theta_{X,t}(\xi(t))\, dt \quad \text{for } \xi \in T_\gamma\Omega. \tag{22}$$

An admissible almost complex structure on $M$ defines a metric on $M$ by definition, and hence an element in $\mathcal{J}$ defines a metric on the relative loop space $\Omega$.

The vector field dual to $-\mathcal{Y}_X$ with respect to this metric is $-\mathcal{V}_X$,

$$\mathcal{V}_X(\gamma) := J_t(\gamma)\frac{\partial \gamma}{\partial t} + \check{\theta}_{X,t}^{J_t}(\gamma) \quad \text{for } \gamma \in \Omega, \tag{23}$$

where $\check{\theta}_{X,t}^{J_t}$ is the dual vector of $\theta_{X,t}$ with respect to the metric on $M$ associated to $J_t$. Formally, $\mathcal{V}_X$ generates a flow on $\Omega$ given by a perturbed Cauchy-Riemann equation

$$\bar{\partial}_{JX} u := \frac{\partial u}{\partial s} + J_t(u)\frac{\partial u}{\partial t} + \check{\theta}_{X,t}^{J_t}(u) = 0, \tag{24}$$

where $s \in \mathbb{R}$ or $S_T^1$ for any $T > 0$, depending on whether the flow line is periodic. $t \in [0, 1]$, and $u$ is a smooth $M$-valued function on $(s, t)$ such that $u(s, \cdot) \in \Omega \; \forall s$.

The *energy* of a solution $u$ to (24) is

$$\mathcal{E}(u) = \int |\partial_s t|^2 ds\, dt.$$

We may now clarify the meaning of the spaces $\mathcal{P}(X), \mathcal{M}_P, \mathcal{M}_O$ that appeared in §2.1.2.

An *s*-independent solution of (24) is called a *critical point*; it is also a solution of the equation of motion

$$\frac{dx(t)}{dt} - X_t(x(t)) = 0, \quad \text{for } x \in \Omega. \tag{25}$$

The *moduli space of critical points*,

$$\mathcal{P}(X) = \mathcal{P}(L, L'; X),$$

is the set of all such solutions.

Let

$$\mathcal{M}_P = \mathcal{M}_P(L, L'; J, X), \mathcal{M}_O := \mathcal{M}_O(L, L'; J, X)$$

be respectively the space of *s-dependent* finite-energy solutions of (24) for real or circle-valued $s$. (We consider circles of all possible length.) They consist of perturbed pseudo-holomorphic strips or annuli respectively. We shall refer to them as the *moduli space of $J|X$-strips* and *moduli space of $J|X$-annuli* respectively.

Given $x, y \in \mathcal{P}(L, L'; X)$, let

$$\mathcal{M}_P(L, L', x, y; J, X) =$$
$$\left\{ u \in \mathcal{M}_P(L, L'; J, X); \lim_{s \to -\infty} u(s, \cdot) = x, \lim_{s \to \infty} u(s, \cdot) = y \text{ in } C^0 \text{ norm.} \right\}$$

By Gromov compactness, $\mathcal{M}_P(L, L'; J, X) = \coprod_{x,y} \mathcal{M}_P(L, L', x, y; J, X)$. Let $(x, [w]), (y, [v]) \subset \tilde{\mathcal{P}}(L, L'; X)$ be lifts of $x, y$ respectively, then

$$\mathcal{M}_P(L, L', (x, [w]), (y, [v]); J, X) \subset \mathcal{M}_P(L, L', x, y; J, X)$$

is the subset of elements $u$ such that $u(s, \cdot)$ lifts to a path in $\tilde{\Omega}$ starting at $(x, [w])$, ending in $(y, [v])$.

For any $A \in H_1(\Omega)$, let $\mathcal{M}_O(L, L', A; J, X) \subset \mathcal{M}_O(L, L'; J, X)$ consist of elements $u$ such that the 1-cycle $S^1 \to \Omega : s \mapsto u(s, \cdot)$ represents $A$. Let

$$\mathcal{M}_O^{k+1}(L, L'; J, X) := \coprod_{A, \psi_\mu(A)=k} \mathcal{M}_O(L, L', A; J, X) \subset \mathcal{M}_O(L, L'; J, X);$$

$$\mathcal{M}_P^k(L, L', x, y; J, X) := \coprod_{(y,[v]), \text{ind}(y,[v])-\text{ind}(x,[w])=k} \mathcal{M}_P(L, L', (x, [w]), (y, [v]); J, X)$$
$$\subset \mathcal{M}_P(L, L', x, y; J, X),$$

where $(y, [v])$ is an arbitrary lift of $y$; $(x, [w])$ is a fixed lift of $x$ in $\tilde{\Omega}$, ind is the Maslov-Viterbo index defined in §3.2.1 below. Note that this definition does not depend on the choice of $(x, [w])$.

There is a free $\mathbb{R}$-action and a semi-free $S^1$-action on $\mathcal{M}_P^k(L, L', x, y; J, X), \mathcal{M}_O(L, L'; J, X)$ respectively by translation in $s$. Let $\hat{\mathcal{M}}_P^{k-1}(L, L', x, y; J, X), \hat{\mathcal{M}}_O(L, L'; J, X)$ be the respective quotients. An element $u \in \mathcal{M}_O$ or $\hat{\mathcal{M}}_O$ is said to be *simple* if its multiplicity $m(u) = 1$; otherwise it is a *multiple cover*.

A direct computation using (24) shows that the energy of a $u \in \mathcal{M}_O(L, L', A; J, X)$ is given by

$$\begin{aligned} \mathcal{E}(u) &= \int |\partial_s u|^2 \, ds \, dt \\ &= -\int \langle \partial_s u, J(u) \partial_t u + \check{\theta}_{X,t}^{J_t} \rangle \, ds \, dt \\ &= -[\mathcal{Y}_X](A) \\ &= \psi_\omega(A) - \int e_0^* \theta_{X,t}(A) \, dt. \end{aligned} \qquad (26)$$

Similarly, the energy of a $\mu \in \mathcal{M}_P(L, L', (x, [w]), (y, [v]); J, X)$ is

$$\mathcal{E}(\mu) = \tilde{\mathcal{A}}_X([y, v]) - \tilde{\mathcal{A}}_X([x, w]).$$

Let $\mathcal{M}_O(L, L'; J, X)^\Re \subset \mathcal{M}_O(L, L'; J, X)$, $\hat{\mathcal{M}}_P(L, L', x, y; J, X)^\Re \subset \hat{\mathcal{M}}_P(L, L', x, y; J, X)$ etc. be the subsets consisting of elements with energy $\leq \Re$.

The energy formulae above show that the expressions (4), (5) may be written simply in this notation as $\hat{\mathcal{M}}_P^0(L, L', x, y; J, X)^\Re$, $\hat{\mathcal{M}}_O^1(L, L'; J, X)^\Re$ respectively.

*Notation.* The moduli spaces described above all depend implicitly on the choice of the base point $\gamma_0$; we shall include the superscript $\gamma_0$ in the notation when we wish to emphasize the dependence. E.g. $\mathcal{M}_P^{\gamma_0}(L, L'; J, X) = \mathcal{M}_P(L, L'; J, X)$.

### 3.1.3 Moving $L'$ by a symplectic isotopy and changing $X$.

Below is a well-known and very useful observation.

Let $\Phi = \{\phi_t \mid t \in [0, 1]\}$ be a symplectic isotopy connecting the identity map to $\phi_1$, and let $\vartheta_t$ be the path of closed 1-forms so that $X_{\vartheta_t}$ generates $\phi_t$. For $\gamma \in \Omega_{\gamma_0}(M; L, L')$, let

$$\Phi \cdot \gamma(t) := \phi_t(\gamma(t)).$$

This defines an isomorphism

$$\Phi : \Omega_{\gamma_0}(M; L, L') \to \Omega_{\Phi \cdot \gamma_0}(M; L, \phi_1(L')).$$

It is easy to see that for any $X \in \mathcal{X}$, there is another $X' \in \mathcal{X}$ such that

$$\Phi^* \mathcal{Y}_{X'} = \mathcal{Y}_X \quad \text{and}$$
$$\phi_t^* \theta_{X',t} := \theta_{X,t} + \vartheta_t.$$

Furthermore, if $u(s, t)$ solves (24), then

$$w(s, \cdot) := \Phi \cdot u(s, \cdot)$$

also satisfies the flow equation (24), but with $L'$ replaced by $\phi_1(L')$, $X$ replaced by $X'$, and $J_t$ replaced by $J'_t$, where

$$J'_t(\phi_t(x)) := (D\phi_t)J_t(x)(D\phi_t)^{-1}.$$

This not only defines isomorphisms between the moduli spaces $\mathcal{P}^{\gamma_0}(L, L'; X)$, $\mathcal{M}_P^{\gamma_0}(L, L'; J, X)$, $\mathcal{M}_O^{\gamma_0}(L, L'; J, X)$ and $\mathcal{P}^{\Phi \cdot \gamma_0}(L, \phi(L'); X')$, $\mathcal{M}_P^{\Phi \cdot \gamma_0}(L, \phi(L'); J', X')$, $\mathcal{M}_O^{\Phi \cdot \gamma_0}(L, \phi(L'); J', X')$ respectively, but also equivalences of the relevant deformation operators ($A_x$, $E_u$, $\tilde{D}_u$ in §3.2) by similarity transformations. Thus the Floer theories associated to $(L, L', \gamma_0; J, X)$ and $(L, \phi_1(L'), \Phi \cdot \gamma_0; J', X')$ are completely equivalent.

Because of this equivalence, in this paper we shall fix the lagrangian submanifolds $L, L'$ and vary the almost complex structure $J$ and the symplectic vector field $X$.

## 3.2 Structure of the moduli spaces

Let $M, L, L'$ be as in §2.2. The purpose of this subsection is to verify the Properties listed in §2.1.2 for this version of Floer theory.

### 3.2.1 The Maslov-Viterbo index of a critical point.

Here we define the index map ind (cf. §2.1.2 (a)). It is a version of the Maslov-Viterbo index.

*Notation.* Our convention is to let $\|\xi\|_{p,k}$ denote the sum of $L^p$ norms of $\xi$ and its derivatives up to the $k$-th order, and $L_k^p(M, P)$ denotes the corresponding Banach space, where $\xi$ is a section of the bundle $P \to M$.

Let $A_x : K_1^p(x^*TM) \to L^p(x^*TM)$ be the linearization of $\mathcal{V}_X$ at $x$, where

$$K_k^p(x^*TM) := \left\{ \xi \,\middle|\, \xi \in L_k^p(x^*TM), \xi(0) \in T_{x(0)}L, \xi(1) \in T_{x(1)}L' \right\}.$$

An $x \in \mathcal{P}(L, L'; X)$ is said to be *nondegenerate* if $A_x$ is surjective. We shall assume below that $\mathcal{P}(L, L'; X)$ consists of finite elements and they are all nondegenerate. In this case we see that $x(t) \neq y(t)$ $\forall t$ for any two distinct elements $x, y \in \mathcal{P}(L, L'; X)$, since they both satisfy the ODE (25). For each $x \in \mathcal{P}(L, L'; X)$, one may now choose a trivialization of the symplectic vector bundle $\Phi_x : x^*TM \to [0, 1] \times \mathbb{C}^n$ such that $\Phi_x(T_{x(0)}L) = \{0\} \times \mathbb{R}^n$; $\Phi_x(T_{x(1)}L') = \{1\} \times i\mathbb{R}^n$.

Similarly, fix a trivialization of the symplectic vector bundle $\Phi_{\gamma_0} : \gamma_0^*TM \to [0, 1] \times \mathbb{C}^n$ satisfying the same boundary conditions.

Given a $(x, [w]) \in \tilde{\mathcal{P}}(L, L'; X)$ for nondegenerate $x$, one may now find a trivialization

$$\Phi_w : w^*TM \to [0, 1] \times [0, 1] \times \mathbb{C}^n$$

of the symplectic vector bundle $w^*TM$ satisfying the following conditions. Regard $w$ as a map $w : [0, 1] \times [0, 1] \to M$ such that $w(\cdot, 0) \in L$; $w(\cdot, 1) \in L'$, $w(0, t) = \gamma_0(t)$, and $w(1, t) = x(t)$. Then we require that

$$\Phi_w\Big|_{w(0,\cdot)^*TM} = \Phi_{\gamma_0};$$
$$\Phi_w\Big|_{w(1,\cdot)^*TM} = \Phi_x;$$
$$\Phi_w(w(s,1)^*TL) = \{s\} \times \mathbb{R}^n.$$

It is easy to see that such a trivialization exists, and the space of such trivializations is contractible. Now

$$\Phi_w(w(s,1)^*TL')), \quad s \in [0,1]$$

gives a loop of lagrangian subspaces in $\mathbb{C}^n$; we define the *Maslov-Viterbo index of* $(x, [w]) \in \tilde{\mathcal{P}}(L, L'; X)$, denoted $\mathrm{ind}(x, [w])$, to be the Maslov index of this loop. This is an $\mathbb{Z}$-valued index, and clearly it is independent of the choice of $\Phi_w$. For $x \in \mathcal{P}(L, L'; X)$, let $\mathrm{ind}(x) \in \mathbb{Z}/\mathbb{N}(L, L')\mathbb{Z}$ be the value of $\mathrm{ind}(x, [w])$ for any lift $(x, [w])$ of $x$.

Thus defined, it follows from standard arguments (see e.g. [F88c]) that for a generic representative $w(s)$ ($s \in [0, 1]$) of $[w]$,

$$\mathrm{ind}(x, [w]) \text{ is the spectral flow of the path of self-adjoint operators } A_{w(s, \cdot)}. \tag{27}$$

ind *does* depend on the choices of the trivializations $\Phi_{\gamma_0}$, $\Phi_x$ $\forall x \in \mathcal{P}(L, L'; X)$. However, if $L, L'$ are oriented and we require in addition that these trivializations are orientation-preserving, then different choices of such trivializations only change the definition of ind by even numbers. On the other hand, in this case $\mathbb{N}(L, L')$ is even; thus the orientations of $L, L'$ alone determine an

$$\text{absolute } \mathbb{Z}/2\mathbb{Z}\text{-index of } \mathcal{P}(L, L'; X). \tag{28}$$

It is easy to see that when $X = 0$, $(-1)^{\mathrm{ind}(p)}$ agrees with the sign of intersection of $L, L'$ at $p$, up to an overall sign.

**Remarks.** The usual definition of Maslov-Viterbo index only gives a *relative* $\mathbb{Z}/\mathbb{N}(L, L')\mathbb{Z}$-valued index (see e.g. [F88c]) instead of an absolute index. The properties (27) and (28) of the index ind are important to us: we need (28) for the definition of $\tau_F$, and (27) in the proof of the invariance theorem 2.2.2.

Similar ideas were used in [Se] to show that for more general $N$, certain more refined structures on the lagrangian manifolds $L, L'$ determine absolute $\mathbb{Z}/N\mathbb{Z}$ gradings for the Floer complex. However, the definitions of such an absolute grading in [Se] (cf. also [FuOOO] §2) are slightly different from ours. We have chosen this description for a more straightforward connection with spectral flow.

### 3.2.2 Transversality and compactness results for the moduli spaces.

We verify the smoothness and compactness properties of the moduli spaces listed in §2.1.2 Properties (b), (c). The orientation issue is postponed to the next subsubsection.

Like [L], we achieve transversality by perturbing the symplectic vector field $X$ with hamiltonian vector fields, fixing $J$ and the lagrangian submanifolds. (This of course implies that transversality may be achieved via perturbation by general symplectic vector fields). It is also possible to achieve transversality by fixing $X$ and $L, L'$, and perturbing $J$ instead.

Recall the definitions of $\mathcal{H}$ and $\chi_H$ in §2.2. We now make precise what "generic" means in the present context.

**Proposition.** *Let $M$, $L, L'$, $J$, $X$ be as in Theorem 2.2.2. Then:*

*(1) There is a Baire set $U_\mathcal{H} \subset \mathcal{H}$ such that for any $H \in \mathcal{H}$, $\mathcal{P}(L, L'; X + \chi_H)$ consists of finitely many nondegenerate critical points.*

*(2) Let $H_0 \in U_\mathcal{H}$. Then for any small positive number $\delta$ and any positive integer $k \geq 2$, there is a Baire set $\mathcal{H}^k_{\text{reg}}(H_0, \delta)$ in*

$$\mathcal{H}^k(H_0, \delta) := \Big\{ H_0 + h \,\Big|\, H_0 + h \in \mathcal{H}, \|h\|_{C_\epsilon} < \delta; \\ \nabla^i h_t(x(t)) = 0 \;\forall x \in \mathcal{P}(L, L'; X + \chi_H), i \in \{1, 2, \ldots, k\} \Big\}, \tag{29}$$

*such that for any $H \in \mathcal{H}^k_{\text{reg}}(H_0, \delta)$, $(J, X + \chi_H)$ is regular in the sense of §2.1.2.*

The issue of orientations will be dealt with in §3.2.3. The rest of the proof follow the standard outline: very roughly, the compactness results are based on the Gromov compactness theorem; to prove smoothness, one describes the moduli space as the zero locus of a Fredholm section of a Banach bundle over the "configuration space" ($\mathcal{B}_P$ for $\mathcal{M}_P$, and $\mathcal{B}_O$ for $\mathcal{M}_O$), then try to show that for generic perturbations, the Fredholm section is transverse to the zero section. Locally, the section is described by the Kuranishi model: a map between the kernel and cokernel of the "deformation operator". Since very similar statements were proven in [L], below we shall only indicate how the proofs should be modified for part (2) above. See [L] section 3 for more details.

We first recall some preliminaries from [L] and clarify certain terminologies. For a $u \in \mathcal{M}_P$ or $(1/T, u) \in \mathcal{M}_O$ ($T$ specifies the period), the deformation operators are the linearization of $\bar{\partial}_{JX}$ at $u$ or $(1/T, u)$, denoted $E_u$ or $\tilde{D}_{(1/T,u)}$ respectively. $\tilde{D}_{(1/T,u)} : \mathbb{R} \oplus K^p_1(u^*TM) \to L^p(u^*TM)$ has the form

$$\tilde{D}_{(1/T,u)}(\varrho, \xi) = D_u \xi + \varrho T \partial_s u,$$

where $u : S^1_T \times [0, 1] \to M$,

$$K^p_k(u^*TM) := \Big\{ \xi \,\Big|\, \xi \in L^p_k(u^*TM), \xi(s,0) \in T_{u(s,0)}L, \xi(s,1) \in T_{u(s,1)}L' \Big\},$$

$p > 2$, and $D_u$ is an elliptic differential operator between $K_1^p(u^*TM)$ and $L^p(u^*TM)$.

$(J, X)$ is said to be $\Re$-*regular* if $\hat{\mathcal{M}}_P(L, L'; J, X)^\Re$, $\mathcal{M}_O(L, L'; J, X)^\Re$ are both smooth and compact.

We now are ready for the:

*Proof.* **(2)** Let's first assume that $M$ is compact. When $u \in \mathcal{M}_P$, the above claim about smoothness follows from arguments in the proof of [FHS] Theorem 5.1 (ii). Although [FHS] (as well as most literature) works with the free loop space instead of the relative loop space $\Omega(M; L, L')$, the arguments may still be copied to our situation, basically because the fact that $L, L'$ are lagrangian makes the boundary contributions vanish in integration-by-parts arguments. For example, to show that an element $\eta$ in the cokernel of the deformation operator $E_u$ vanishes, it suffices to show that it vanishes in the interior $(0, 1) \times \mathbb{R} \subset [0, 1] \times \mathbb{R}$. Let $E_u^*$ be the formal $L^2$-adjoint of $E_u$. Note that $\eta \in L_1^q(u^*TM)$ with $1/q + 1/p = 1$, and $E_u^*\eta = 0$: for any $\xi \in T_u\mathcal{B}_P$,

$$\int_{\mathbb{R} \times [0,1]} \langle \eta, E_u\xi \rangle \, ds \, dt = \int_{\mathbb{R} \times [0,1]} \langle E_u^*\eta, \xi \rangle \, ds \, dt + \text{boundary contribution}$$

of the form $\int \langle \eta(s, 1), J_1(u(s, 1))\xi(s, 1) \rangle ds - \int \langle \eta(s, 0), J_0(u(s, 0))\xi(s, 0) \rangle ds$,

where $\langle \cdot, \cdot \rangle$ denotes the pointwise inner-product given by the metric induced by $J_t$. The boundary contribution above vanishes precisely when $\xi, \eta$ satisfy the specified boundary conditions. Now one may apply the unique continuation theorem to show that $\eta$ vanishes. The needed injectivity results for interior points of $\mathbb{R} \times [0, 1]$ is provided by the arguments in [FHS] section 4.

For $\mathcal{M}_O$, smoothness may be obtained by modifying [L] in the same manner.

The claims about compactness follow from the admissibility assumption on $L, L'$ and the transversality established above, since they exclude the possibility of bubbling (cf. [Oh93]).

When $M$ is noncompact, one needs the following additional observations.

First, note that by the monotonicity lemma ([Law] 3.15; [MS] p.44) all elements $u$ in $\mathcal{M}_P(L, L'; J, X + \chi_{H_0})$, $\mathcal{M}_O(L, L'; J, X + \chi_{H_0})$ take values in a compact region in $M$ for fixed $H_0$. Thus compactness of $\hat{\mathcal{M}}_P(L, L'; J, X + \chi_{H_0})^\Re$ and $\mathcal{M}_O(L, L'; J, X + \chi_{H_0})^\Re$ follows from the same argument as in the case when $M$ is compact. In fact, for all sufficiently small $\delta > 0$ and any $h$ such that $H_0 + h \in \mathcal{H}^k(H_0, \delta)$, there is a compact neighborhood $U(H_0, \delta, \Re) \subset M$ containing the image of all $u \in \mathcal{M}_P(L, L'; J, X + \chi_{H_0+h})^\Re$ and $\mathcal{M}_O(L, L'; J, X + \chi_{H_0+h})^\Re$. We shall choose $U(H_0, \delta, \Re)$ large enough so that it also contains $L, L'$, and the support of $X + \nabla H_0$.

We noticed before that when $M$ is noncompact, $\mathcal{H}$ is not complete and thus may not be used directly in the Sard-Smale theorem. Instead, given $H_0 \in \mathcal{H}$, and any compact neighborhood $U' \supset U(H_0, \delta, \Re)$, we use the Banach manifold

$$\mathcal{H}^k(H_0; \delta, U') := \left\{ H_0 + h \mid H_0 + h \in \mathcal{H}^k(H_0, \delta), h \text{ is supported in } U' \right\}$$

in the Sard-Smale theorem to achieve transversality. In other words, in view of the compactness of $\mathcal{M}_P^\Re(L, L'; J, X + \chi_{H_0+h})$ and $\mathcal{M}_O^\Re(L, L'; J, X + \chi_{H_0+h})$, for each $\Re$ and $U'$ containing $U(H_0, \delta, \Re)$, the subset $\mathcal{H}_{\Re-\mathrm{reg}}^k(H_0; \delta, U') \subset \mathcal{H}^k(H_0; \delta, U')$ of $H$ for which $(J, X + \chi_H)$ is $\Re$-regular is open and dense. Note that

$$\mathcal{H}_{\Re-\mathrm{reg}}^k(H_0; \delta, U') \subset \mathcal{H}_{\Re-\mathrm{reg}}^k(H_0; \delta, U'') \text{ for } U' \subset U'',$$

and the subset of $\mathcal{H}^k(H_0; \delta)$ consisting of all $H$ such that $(J, X + \chi_H)$ is $\Re$-regular is

$$\mathcal{H}_{\Re-\mathrm{reg}}^k(H_0; \delta) = \bigcup_{U' \supset U(H_0, \delta, \Re)} \mathcal{H}_{\Re-\mathrm{reg}}^k(H_0; \delta, U').$$

Hence $\mathcal{H}_{\Re-\mathrm{reg}}^k(H_0; \delta) \subset \mathcal{H}^k(H_0; \delta)$ is also open and dense. Thus, the subset of $H$ such that $(J, X + \chi_H)$ is regular is

$$\mathcal{H}_{\mathrm{reg}}^k(H_0; \delta) = \bigcap_\Re \mathcal{H}_{\Re-\mathrm{reg}}^k(H_0; \delta) \subset \mathcal{H}^k(H_0; \delta)$$

is Baire. □

### 3.2.3 Orientations.

We now specify the orientations on the relevant moduli spaces. Fix a pair of regular $(J, X)$ throughout this subsection.

**(a) Orienting $\mathcal{M}_O^1$.** Since our definition of $I_F$ only uses the 1-dimensional moduli space $\mathcal{M}_O^1$, we shall restrict our discussion to this case. Higher dimensional $\mathcal{M}_O^k$ may be oriented from the chosen orientation of $\mathcal{M}_O^1$ by an excision argument; for an example see the proof of Proposition 21.3 in [FuOOO].

Let $\mathcal{B}_O^1 \subset \mathcal{B}_O$ be the path component where $\mathrm{ind}\, \tilde{D}_{(1/T,u)} = 1$. There are determinant line bundles $\det \tilde{D}$, $\det D$ over $\mathcal{B}_O^1$, with fiber

$$\det \tilde{D}_{(1/T,u)} := \det \ker \tilde{D}_{(1/T,u)} \otimes (\det \mathrm{coker}\, \tilde{D}_{(1/T,u)})^*, \quad \det D_u$$

respectively over $(1/T, u) \in \mathcal{B}_O^1$, where det denotes either determinant or top exterior product. The two determinant bundles are isomorphic via "stabilization" (see e.g. [L] §5.1.1). To orient $\mathcal{M}_O^1$, it suffices to orient $\det \tilde{D}$, which in turn may be induced from an orientation of $\det D$ ([L] equation (5.2)). We orient $\det D_u$ for $(1/T, u) \in \mathcal{B}_O^1$ as follows.

The orientability of $L$ and $L'$ means that $TL\big|_{u(\cdot,0)}$ and $TL'\big|_{u(\cdot,1)}$ are trivializable. Fix a pair of such trivializations.

*Claim.* The homotopy classes of stable trivializations of $TL\big|_{u(\cdot,0)}$ and $TL'\big|_{u(\cdot,1)}$, together with the homotopy class of $\Phi_{\gamma_0}$ in §3.2.1, assigns an orientation of $\det D_u$.

To prove the Claim, we may assume without loss of generality that we have trivializations of $TL\big|_{u(\cdot,0)}, TL'\big|_{u(\cdot,1)}$. (In general, we consider stabilizations of $D_u$ given stable trivializations of $TL\big|_{u(\cdot,0)}, TL'\big|_{u(\cdot,1)}$, but orienting determinant line bundles of stabilizations of $D_u$ is equivalent to orienting $\det D_u$). Note that the assumption that $(1/T, u) \in \mathcal{B}_O^1$ implies $\psi_\mu([u]) = 0$, where $[u] \in H_1(\Omega)$ is the homology of $u$ as a loop in $\Omega$. Thus, there is a trivialization $\Phi_u : u^*TM \to S_T^1 \times [0,1] \times \mathbb{C}^n$ of $u^*TM$ such that

$$\Phi_u(T_{u(s,0)}L) = \{(s,0)\} \times \mathbb{R}^n, \quad \Phi_u(T_{u(s,1)}L') = \{(s,1)\} \times i\mathbb{R}^n \quad \forall s$$

agrees precisely with the chosen trivializations over $\partial(S_T^1 \times [0,1])$. [1]

Now via this trivialization, $D_u$ is identified with a differential operator on $\Gamma(S^1 \times [0,1], \mathbb{C}^n)$, namely $\Phi_u D_u \Phi_u^{-1}$. We may thus assign an orientation on $\det D_u$ (independent of the choice of $\Phi_u$) by orienting the determinant line bundle $\mathcal{L}_O$ on the space of such differential operators, $\Sigma_O$. More precisely,

$$\Sigma_O := \Big\{ \partial_s + \alpha j(s,t) \partial_t + \nu(s,t) : K_1^p(\mathbb{C}^n) \to L^p(\mathbb{C}^n) \,\Big|\, \alpha \in \mathbb{R}^+,$$
$$j \text{ is a } C_\epsilon \text{ complex structure}, \nu \text{ is a } C_\epsilon \text{ matrix-valued function},$$
$$(s,t) \in S_T^1 \times [0,1] \text{ for some } T > 0, \Big\},$$

where $\mathbb{C}^n$ denotes the trivial bundle over $S_T^1 \times [0,1]$;

$$K_k^p(\mathbb{C}^n) = \Big\{ \xi(s,t) \,\Big|\, \xi \in L_k^p(\mathbb{C}^n), \xi(\cdot,0) \in \mathbb{R}^n, \xi(\cdot,1) \in i\mathbb{R}^n \Big\}.$$

Clearly $\Sigma_O$ is contractible, since the space of possible $\nu$ can be contracted to 0, and the space of almost complex structures is well-known to be contractible. Thus $\mathcal{L}_O$ is orientable, and we may choose the orientation by fixing the orientation of $\det D_0$ at any $D_0 \in \Sigma_O$. We choose

$$D_0 = \partial_s + \Phi_{\gamma_0} A_{\gamma_0} \Phi_{\gamma_0}^{-1}.$$

Thus chosen, $\det D_0$ has a canonical orientation, since by the $s$-independence of $D_0$, $\ker D_0$ and $\operatorname{coker} D_0$ may respectively be identified with $\ker A_{\gamma_0}$ and $\operatorname{coker} A_{\gamma_0}$, which are isomorphic by the self-adjointness of $A_{\gamma_0}$.

Of course, the choice of $\Phi_u$ is not unique; we now show that this orientation is independent of the choice. The space of the choices of $\Phi_u$'s has $\mathbb{Z}$ path-components, related by transformations of the form

$$\phi_k(s,t) = e^{ik\pi\beta(t)} \oplus \operatorname{Id}^{n-1} \in \operatorname{Sp}_n \text{ for } (s,t) \in S^1 \times [0,1], k \in \mathbb{Z},$$

where $\beta$ is a smooth real cutoff function with values in $[0,1]$, which is 0 at $t = 0$ and is 1 at $t = 1$. The independence thus follows from the following lemma:

---

[1] One may also remove the second constraint, so that space of $\Phi_u$ is connected. But then instead of $\Sigma_O$ below, one would have to work with the larger space $\tilde{\Sigma}_O$–the space of all possible differential operators of the same form, but with possibly different boundary conditions. $\tilde{\Sigma}_O$ fibers over $\Sigma_O$ with noncontractible fibers.

**Lemma.** Let $\mathrm{Or}(\mathcal{L}_O) = \mathbb{Z}/2\mathbb{Z}$ be the space of possible orientations of the determinant line bundle $\mathcal{L}_O \to \Sigma_O$. Then the $\mathbb{Z}$-action of $\{\phi_k\}$ by conjugation on $\mathrm{Or}(\mathcal{L}_O)$ is trivial.

*Proof.* This is essentially [FH] Lemma 13, with the roles of $s, t$ permuted. The idea is to consider the action of $\phi_k$ on an operator with small $\alpha$. However, conjugating by $\phi_k$ on such an operator only changes it by a term small in the operator norm. □

End of proof for the Claim. We now show how the relative spin structure of $L, L'$ determines the homotopy classes of stable trivializations of $TL\big|_{u(\cdot,0)}$, $TL'\big|_{u(\cdot,1)}$ consistently. For simplicity, assume first that $L, L'$ are spin. In this case, an orientation and a spin structure of $L$ determines a homotopy class of (stable) trivializations of $TL$ over the 1-skeletons of $L$; similarly for $L'$. Our trivializations should be chosen in these classes. Furthermore, such a choice is consistent for different $u$ in $\mathcal{B}_O^1$. To see this, it suffices to show that given a loop $\{u_\lambda \in \mathcal{B}_O^1 \,|\, \lambda \in S^1\}$, $TL\big|_{\bigcup_{\lambda,s} u_\lambda(s,0)}$, $TL'\big|_{\bigcup_{\lambda,s} u_\lambda(s,1)}$ may be trivialized. Since $\bigcup_{\lambda,s} u_\lambda(s,0)$, $\bigcup_{\lambda,s} u_\lambda(s,1)$ trace out two-cycles in $L$ and $L'$ respectively, their (stable) trivializability is guaranteed by the spinness of $L, L'$ respectively, which implies the (stable) trivializability of $TL, TL'$ over the two-skeletons of $L, L'$. (Stabilization is necessary when $n = \dim L < 3$).

We refer the reader to [FuOOO] §21 for the generalization to the relative spin case. The idea is that the $w_2$ condition makes possible a stabilization argument, reducing the proof basically to the spin case.

**(b) Orienting $\mathcal{M}_P$.** For each pair $x, y \in \mathcal{P}(L, L; X)$, and any $k \in \{0\} \cup \mathbb{Z}^+$, one may orient $\mathcal{M}_P^k(L, L', x, y; J, X)$ in a similar manner. Let $\mathcal{B}_P^k(L, L', x, y; J, X) \subset \mathcal{B}_P$ denote the component containing $\mathcal{M}_P^k(L, L', x, y; J, X)$. Then each $u \in \mathcal{B}_P^k(L, L', x, y; J, X)$ extends to define a map $\bar{u} : \bar{\mathbb{R}} \times [0, 1] \to M$, where $\bar{\mathbb{R}} := \mathbb{R} \cup \{-\infty, \infty\} \simeq [0, 1]$.

Similarly to **(a)** above, we may find a trivialization $\Phi_u : \bar{u}^* TM : [0,1] \times [0,1] \times \mathbb{C}^n$ of $\bar{u}^* TM$ that satisfies:

$$\Phi_u\big|_{\bar{u}(0,\cdot)^* TM} = \Phi_x;$$

$$\Phi_u\big|_{\bar{u}(1,\cdot)^* TM} = \Phi_y;$$

$$\Phi_u(T_{\bar{u}(s,0)}L) = \{(s,0)\} \times \mathbb{R}^n;$$

$$\Phi_u(T_{\bar{u}(s,1)}L') = \{(s,1)\} \times \psi_k(s)(i\mathbb{R}^n) \quad \forall s \in [0,1],$$

where $\psi_k(s) = e^{i\pi k(1-s)} \oplus \mathrm{Id}_{n-1} \in \mathrm{Sp}_n$, and that $\Phi_u\big|_{\bar{u}(\cdot,0)^* TL}$ and $\psi_{-k}\Phi_u\big|_{\bar{u}(\cdot,1)^* TL'}$ agree respectively with the stable trivializations induced from the relative spin structure of $L, L'$. $\Phi_x, \Phi_y$ are specified in §3.2.1. One may then use the argument described in **(a)** to show that the determinant line bundle $\det E$ is orientable over $\mathcal{B}_P^k(L, L', x, y; J, X)$. However, in this case we do not have a canonical orientation. Instead, one chooses a coherent orientation for each pair of $(x, y)$ according to the scheme of [FH].

# 4 An $S^1$-equivariant Floer theory

This section contains foundational materials for an $S^1$-equivariant version of Floer theory. We shall describe an $S^1$-invariant flow on the loop space $\mathcal{C}$, which induces a flow on the quotient space $\underline{\mathcal{C}} = \mathcal{C}/S^1$. $I_F^{S^1}$ is defined by setting $\Omega = \underline{\mathcal{C}}$. All the moduli spaces involved in §2.1.2 will be the quotient spaces in this version: $\mathcal{P}$ should be replaced by $\underline{\mathcal{P}} = \mathcal{P}/S^1$, $\mathcal{M}_P$ by $\mathcal{M}_P/S^1$ etc.

We reiterate that this is not a full-blown $S^1$-equivariant version since the $S^1$ action is free under our assumptions, and that since this section is parallel to the last one, we shall often adopt the same or similar notations without further explanation. We shall however provide more details, since there is no adequate reference for this version of Floer theory in the literature.

## 4.1 Topological preliminaries and basics

Recall the definitions and notations in §2.3. As in §2.3, we assume for simplicity that $(M, \omega)$ is a compact monotone symplectic manifold.

Regarding $c_1 = c_1(TM), \omega$ as maps from $H_2(M; \mathbb{Z})$ to $\mathbb{Z}, \mathbb{R}$ respectively, we define

$$\underline{\psi}_c := c_1 \circ \underline{\text{im}}; \quad \underline{\psi}_\omega := \omega \circ \underline{\text{im}}.$$

$\underline{\psi}_c$ will be the SF-homomorphism $\psi$ in this version of Floer theory (cf. §4.2.1). Let $\underline{\mathfrak{H}} := \underline{\text{im}}(\ker \underline{\psi}_c) \subset H_2(M; \mathbb{Z})$.

Let $\mathcal{L}M$ be space of (based) contractible loops. Again let $e_0 : \mathcal{C} \to M$ be the evaluation map $e_0(\gamma) = \gamma(0)$. From the fibrations

$$\mathcal{L}M \xrightarrow{f_0} \mathcal{C} \xrightarrow{e_0} M, \quad S^1 \to \mathcal{C} \xrightarrow{\text{qt}} \underline{\mathcal{C}},$$

we have the following exact sequences by abelianizing the associated homotopy sequences:

$$\pi_2(M) \xrightarrow{f_{0*}} H_1(\mathcal{C}; \mathbb{Z}) \xrightarrow{e_{0*}} H_1(M; \mathbb{Z}) \to 0; \tag{30}$$

$$H_1(S^1; \mathbb{Z}) \xrightarrow{f_S} H_1(\mathcal{C}; \mathbb{Z}) \xrightarrow{e_S} H_1(\underline{\mathcal{C}}; \mathbb{Z}) \to 0. \tag{31}$$

We used the assumption that $[\gamma_0]$ is a central element above. In general, the fiber of $e_0$ might have several path components, and instead of (31), one would have:

$$\pi_2'(M) \to H_1(\mathcal{C}; \mathbb{Z}) \to (\pi_1^{[\gamma_0]}(M))^{\text{ab}},$$

where $\pi_1^{[\gamma_0]}(M)$ is the stabilizer of $[\gamma_0] \in \pi_1(M)$ under the conjugation action of $\pi_1(M)$ on itself; $\pi_2'(M)$ is the fixed-point-set under the action of $\pi_1(M)$ on $\pi_2(M)$.

Since $e_{0*} \circ f_S$ maps a generator of $H_1(S^1)$ to $[\gamma_0] \in H_1(M)$, we have

$$\pi_2(M) \xrightarrow{f} H_1(\underline{\mathcal{C}}; \mathbb{Z}) \xrightarrow{e} H_1(M; \mathbb{Z})/\mathbb{Z}[\gamma_0] \to 0. \tag{32}$$

Note also that the maps $\mathrm{im}, \underline{\mathrm{im}}$ and (31) fit into a commutative diagram:

$$H_1(S^1;\mathbb{Z}) \xrightarrow{f_S} H_1(\mathcal{C};\mathbb{Z}) \longrightarrow H_1(\underline{\mathcal{C}};\mathbb{Z}) \longrightarrow 0. \tag{33}$$

with $\mathrm{im}$ from $H_1(\mathcal{C};\mathbb{Z})$ and $\underline{\mathrm{im}}$ from $H_1(\underline{\mathcal{C}};\mathbb{Z})$ into $H_2(M;\mathbb{Z}) \xrightarrow{c_1} \mathbb{Z}$.

On the other hand, (32) fits into the following useful commutative diagram:

$$\begin{array}{c} \ker \underline{\psi}_c \\ \downarrow \text{rest} \\ \pi_2(M) \xrightarrow{\underline{f}} H_1(\underline{\mathcal{C}};\mathbb{Z}) \xrightarrow{\underline{e}} H_1(M;\mathbb{Z})/\mathbb{Z}[\gamma_0] \end{array} \tag{34}$$

with Hurewicz, $\underline{\mathrm{im}}$ into $H_2(M;\mathbb{Z})$, then $c_1$ to $\mathbb{Z}$ and $\omega$ to $\mathbb{R}$,

where "Hurewicz" denotes the Hurewicz homomorphism, and "rest" denotes the restriction.

Simple diagram chasing gives us:

**4.1.1 Lemma.** *If $M$ is monotone, then:*
*(a) There is a $\theta_\omega \in H^1(M;\mathbb{R})$ such that $\theta_\omega[\gamma_0] = 0$ and*

$$\underline{\mathrm{im}}^*\omega\Big|_{\ker \underline{\psi}_c} = \underline{e}^*\theta_\omega\Big|_{\ker \underline{\psi}_c};$$

*(b) If $a \in \ker \underline{\psi}_c$ such that $\underline{\mathrm{im}}^*\omega(a) \neq 0$, then $\underline{e}(a)$ is nontorsion.*

*Proof.* Both are straightforward. For (a), notice that the monotonicity assumption and the diagram (34) imply:

$$\underline{\mathrm{im}}^*\omega\Big|_{\underline{f}(\pi_2(M))\cap \ker \underline{\psi}_c} = 0. \tag{35}$$

The exact sequence part of the diagram then implies

$$\mathrm{rest}^* \underline{\mathrm{im}}^*\omega \in \mathrm{rest}^* \underline{e}^*(H_1(M;\mathbb{R})/\mathbb{R}[\gamma_0]),$$

hence the claim.

For (b), it suffices to show that $\underline{e}(a) \neq 0$, since if $m\underline{e}(a) = 0$ for $m \in \mathbb{Z}^+$, then $\underline{e}(ma) = 0$, reducing to the same proof.

By (34), $\underline{e}(a) \neq 0$ means $a \notin \mathrm{Image}(\underline{f})$, which follows from (35). □

**4.1.2 Definition.** A class $\theta \in H^1(M; \mathbb{R})$ is $\gamma_0$-*exact* iff $\theta([\gamma_0]) = 0$.

A $\gamma_0$-exact class $\theta \in H^1(M; \mathbb{R})$ is $H^2$-*induced* iff $\mathrm{rest}^* \, \underline{e}^* \theta = \mathrm{rest}^* \, \underline{\mathrm{im}}^* \vartheta$ for some $\vartheta \in H^2(M; \mathbb{R})$.

A symplectic vector field $X$ is said to be $\gamma_0$-exact or $H^2$-induced if $\theta_X$ is in a $\gamma_0$-exact or $H^2$-induced class.

Note that $0$ is by definition an $H^2$-induced class, and so is $\theta_\omega$ under the monotonicity assumption by Lemma 4.1.1 (b) above. Furthermore, the definition 4.1.2 introduces a homomorphism

$$\mathfrak{h} : \{H^2\text{-induced classes}\} \to \mathrm{Hom}(\mathfrak{H}, \mathbb{R}).$$

Let $\mathcal{X}(\theta) \subset \mathcal{X}$ be the subspace of symplectic vector fields $X$ of flux $[\theta_X] = \theta$.

Given a pair $(J, X) \in \mathcal{J} \times \mathcal{X}$, one may define similarly to §3.1.2 a $L^2$-metric and a closed 1-form $\mathcal{Y}_X$ on $\mathcal{C}$:

$$\mathcal{Y}_X(\xi) := -\int_{S^1} \omega(\partial_t \gamma, \xi) \, dt + \int_{S^1} \theta_X(\xi) \, dt \quad \text{for } \xi \in T_\gamma \mathcal{C} = \Gamma(\gamma^* TM).$$

It lifts to $\tilde{\mathcal{Y}}_X = d\tilde{\mathcal{A}}_X$ for an $S^1$-invariant $\tilde{\mathcal{A}}_X$ due to the $t$-independence and $\gamma_0$-exactness of $X$. Thus $\mathcal{Y}_X, \tilde{\mathcal{A}}_X$ descend to be defined on $\underline{\mathcal{C}}, \underline{\tilde{\mathcal{C}}}$ respectively, which we shall denote by the same notation.

It follows immediately from the definitions that

$$\underline{\mathrm{im}}^*[\mathcal{Y}_X]\Big|_{\mathfrak{H}} = [\omega]\Big|_{\mathfrak{H}} - \mathfrak{h}[\theta_X].$$

On the other hand, the $L^2$-metric on $\mathcal{C}$ induced by $J$ is also $S^1$-invariant due to the $t$-independence of $J$. The (formal) flow generated by $\mathcal{Y}_X$ with respect to this $L^2$-metric on $\mathcal{C}$ is defined by the perturbed Cauchy-Riemann equation (24) (without the subscript $t$). We may then define

$$\mathcal{P} = \mathcal{P}(M, \gamma_0; X) = \mathrm{PO}(X, [\gamma_0]),$$

and the *moduli space of $J|X$-cylinders* $\mathcal{M}_P = \mathcal{M}_P(M, \gamma_0; J, X)$ and its many variants as in §3.1. The definition of the moduli space of $J|X$-tori, $\mathcal{M}_O = \mathcal{M}_O(M, \gamma_0; J, X)$, however needs some explanation. As mentioned in §2.3, we need to consider not only closed orbits under the formal flow on $\mathcal{C}$, but more generally, also "twisted orbits" $u : [0, T] \to \mathcal{C}, u(T) = R_q u(0)$, for some $q \in S^1$. To be more explicit, let $\Theta_{T,Q}$ be the complex torus

$$\Theta_{T,Q} := \mathbb{C}/\Gamma_{T,Q}, \tag{36}$$

where $\Gamma_{T,Q}$ is the lattice spanned by $\{i, T + iQ\}$.

$$\mathcal{M}_O = \mathcal{M}_O(M, \gamma_0; J, X) :=$$
$$\left\{ (1/T, Q, u) \,\Big|\, (T, Q) \in \mathbb{R}_+ \times S^1_1; \, u : \Theta_{T,Q} \to M; \, \bar{\partial}_{JX} u = 0 \right\}. \tag{37}$$

We shall sometimes drop the parameters $1/T, Q$ and simply say that $u \in \mathcal{M}_O$ when there is no danger of confusion.

In other words, by taking twisted orbits into account, we consider (perturbed) pseudo-holomorphic maps from tori of *all* conformal structure.

By the above construction, the moduli spaces $\mathcal{P}, \mathcal{M}_P, \mathcal{M}_O$ all inherit a free $S^1$ action from that on $\mathcal{C}$; in addition, there is a $\mathbb{R}$-action on $\mathcal{M}_P$ and another $S^1$ action on $\mathcal{M}_O$ from translation in $s$. To distinguish the two $S^1$ actions, we shall denote them by $S_t^1, S_s^1$ respectively. Let $\underline{\mathcal{P}} = \mathcal{P}/S_t^1$; $\hat{\mathcal{M}}_P = \mathcal{M}_P/(\mathbb{R} \times S_t^1)$; $\hat{\mathcal{M}}_O = \mathcal{M}_O/(S_t^1 \times S_s^1)$, and similarly for variants of these moduli spaces.

## 4.2 The critical $S^1$-orbits and their neighborhoods

A *critical $S^1$-orbit* refers to either an element in $\underline{\mathcal{P}}$, or a corresponding orbit of the $S^1$ action on $\mathcal{P}$. Our convention is to denote the critical $S^1$-orbit containing $x \in \mathcal{P}$ by $\underline{x}$. We denote an element in $\underline{\tilde{\mathcal{P}}}$ as $(\underline{x}, [w])$: this refers to the $\mathbb{R}$-orbit in $\tilde{\mathcal{C}}$ containing $(x, [w]) \in \tilde{\mathcal{P}}$. ($\mathbb{R}$ is the universal abelian covering of the $S^1$-orbit $\underline{x}$, and $[w]$ is the image of $[w]$ under the map $e_S$ in (31). $[w]$ is the class of paths connecting $\gamma_0$ and $x$ as in §3.1).

A critical $S^1$-orbit $\underline{x}$ is *nondegenerate* if for any of its representative $x \in \mathcal{P}$, $\ker A_x$ is the 1-dimensional space spanned by $\dot{x} = dx/dt$. ($\ker A_x$ always contains $\dot{x}$ in this theory due to $S^1$-invariance.) If $\ker A_x$ is 2-dimensional, $\underline{x}$ is said to be *minimally degenerate*.

A symplectic vector field $X \in \mathcal{X}$ is said to be *nondegenerate* if all $x \in \underline{P}(M, \gamma_0; X)$ are nondegenerate.

*Notation.* Our convention is to denote $\dot{u} = \partial_t u$, $u' = \partial_s u$.

**4.2.1 An index for critical $S^1$-orbits.** We now specify the index map ind in this version of Floer theory. It is a variant of the Conley-Zehnder index for nondegenerate or minimally-degenerate[2] critical $S^1$-orbits. Suppose $X \in \mathcal{X}$ is such that $\underline{\mathcal{P}} = \underline{\mathcal{P}}(M, \gamma_0; X)$ contains only degenerate or minimally-degenerate critical points. Recall that the Conley-Zehnder index assigns an integer $CZ(A)$ to each surjective operator $A \in \Sigma_C$, where $\Sigma_C$ is the space of operators $A : L_1^p(S^1, \mathbb{R}^{2n}) \to L^p(S^1, \mathbb{R}^{2n})$ of the form $A = J_0 \partial_t + \nu(t)$, where $\mathbb{R}^{2n}$ is equipped with the standard symplectic structure and complex structure $J_0$, and $\nu$ is a $C_\epsilon$-function taking values in the space of self-adjoint matrices. Furthermore, the spectral flow of a path of such operators $A(s)$, $s \in [0, 1]$ is the difference in the Conley-Zehnder indices of $A(1)$ and $A(0)$. (See e.g. [Sa97]). We fix a unitary trivialization of $\gamma_0^* TM$. This trivialization of $\gamma_0^* TM$ induces a (homotopy class of) trivialization of $x^* TM$, $\Phi_{x,[w]}$, for each $(x, [w]) \in \tilde{\mathcal{P}}$, and $\mathbb{A}_{(x,[w])} := \Phi_{x,[w]} A_x \Phi_{x,[w]}^{-1}$ is a (non-surjective) operator in $\Sigma_C$. We define

$$\text{ind}(x, [w]) = \text{CZ}(\mathbb{A}_{(x,[w])} - \delta) \in \mathbb{Z},$$

where $\delta$ is a positive number much smaller than the magnitude of any nonzero eigenvalue of $A_x$. It is clear that $\text{ind}(x, [w])$ does not depend on the choice of $\Phi_{x,[w]}$ and $\delta$. Since $\mathcal{P}$ is compact (see the following Lemma), we can and shall fix once and for all an $\delta$, which is much smaller than the magnitude of any nonzero eigenvalue of $A_x$, *for all $x \in \mathcal{P}$*.

---

[2]In the case of minimally degenerate critical points, the index ind defined below is the analog to $\text{ind}_+$ in [L].

For two different liftings $(x, [w]), (x, [v]) \in \tilde{\mathcal{P}}$ of $x$,

$$\text{ind}(x, [w]) - \text{ind}(x, [v]) = c_1 \circ \text{im}([w - v]).$$

From the commutative diagram (33) we see that the index ind defined above descends to define an index map

$$\text{ind} : \tilde{\mathcal{P}} \to \mathbb{Z}$$

satisfying the property (3) with $\Omega = \mathcal{C}$ and the SF-homomorphism $\psi = \underline{\psi}_c$.

Hence it also defines an index of critical $S^1$-orbits with values in $\mathbb{Z}/(2N\mathbb{Z})$, where $N$ is the g.c.d of the values of $\underline{\psi}_c$.

By now it should have become clear that §2.1.2 (d) is just a simple consequence of our setup.

Next we observe that nondegenerate symplectic vector fields are generic.

**4.2.2 Lemma.** *For any $\gamma_0$-exact class $\theta$, the set of nondegenerate symplectic vector fields $\mathcal{X}_{nondeg}(\theta) \subset \mathcal{X}(\theta)$ is open and dense. In fact, for any $X \in \mathcal{X}_{nondeg}(\theta)$, $\mathcal{P}(M, \gamma_0; X)$ consists of a finite number of points.*

*Proof.* Obviously $\mathcal{P}$ is compact, which can be seen, for example, from its interpretation as the fixed-point-set of a symplectomorphism of the compact manifold $M$. Since $\text{ind}(A_x) = 0$, we only need to demonstrate that $\mathcal{X}_{nondeg}(\theta)$ is Baire.

To see this, first observe that for all $x \in \mathcal{P}$, the map $x : S_1^1 \to M$ is injective. In fact, $\dot{x}(t) \neq 0$ for any $t$ and $x \in \mathcal{P}$. This is because $\dot{x}$ satisfies an ODE; whenever it vanishes at a point, it vanishes identically, which is impossible because by assumption $x$ traces out a loop in $M$ with nontrivial homotopy class $\gamma_0$. Similarly, $x$ can not have self-intersection; otherwise it is a multiple cover of another curve, which can not happen either, since $\gamma_0$ is primitive.

The rest of the proof is basically the standard fare via the Sard-Smale theorem: one need only show that the image of $A_x$, $\dot{x}$, and $x^*\chi$ $\forall \chi \in T_X \mathcal{X}(\theta) = \mathcal{X}(\theta)$ span $L^p(S_1^1, x^*TM)$ for any $X \in \mathcal{X}(\theta)$ and $x \in \mathcal{P}(M, \gamma_0; X)$. If this were not true, there exists a $\eta \in L^q(S_1^1, x^*TM)$ such that $A_x\eta = 0$, $\langle \dot{x}, \eta \rangle = 0$, and $\langle x^*\chi, \eta \rangle = 0$ $\forall \chi$. Given the injectivity of $x$, it is easy to construct a smooth vector field $\chi \in T_X \mathcal{X}(\theta)$ violating the last condition. □

### 4.2.3 Asymptotics near critical $S^1$-orbits.
We shall need to understand the behavior of flows near critical $S^1$-orbits.

It is standard that when $y \in \mathcal{P}$ is in a nondegenerate critical $S^1$-orbit, then the flows decay exponentially to $y$, with an exponent greater than $\delta$. (This is the Bott-Morse case).

The minimally degenerate case only needs the following minor modification of [L] section 4. Let $y \in \mathcal{P}$ be in a minimally degenerate critical $S^1$-orbit. Via the standard arguments (see [L] §4.1 and [MMR]), one may again reduce the problem to describing the behavior of flows in the center manifold to $y$. In this case $A_y$ has a 2-dimensional zero eigenspace, spanned by the two orthonormal vectors $\dot{y}/\|\dot{y}\|_{2,t}$ and $\mathbf{e}_y$. Thus, the center manifold is now a 2-manifold with a free $S^1$-action. The vector field $\mathcal{V}_X$ is tangent to the the center manifold, and is orthogonal to this $S^1$-action everywhere. We thus have a description of the center manifold near the critical $S^1$-orbit $\underline{y}$ as a cylinder $\bigcup_{q \in S^1} R_q \mathbb{R}_y^+$, where $\mathbb{R}_y^+$ is an integral curve of the flow to $y$. By the orthogonality of $\mathcal{V}_X$ to the $S^1$-action, $T_y \mathbb{R}_y^+ = \mathbb{R}\mathbf{e}_y$; furthermore, we may further reduce the description of the flow on the central manifold to a description of the flow along $\mathbb{R}_y^+$.

In a local coordinate system centered at $y$, $\mathbb{R}_y^+$ is described as

$$\{\beta \mathbf{e}_y + \zeta(\beta) \,|\, \beta \in \mathbb{R}^+\},$$

where $\zeta : \mathbb{R} \to \ker A_y^\perp \oplus \mathbb{R}\dot{y}$ is a smooth function vanishing up to the first order at 0. The flow on $\mathbb{R}_y^+$ is again described by the following equation (cf. [L] (4.3)):

$$d\beta(s)/ds + \langle \mathbf{e}_y, n_y(\beta \mathbf{e}_y + \zeta(\beta)) \rangle_{2,t} = 0,$$

where $n_y$ is the nonlinear part of $\mathcal{V}_X$ in local coordinates. We choose the orientation of $\mathbf{e}_y$ such that $\langle \mathbf{e}_y, \nabla_{\mathbf{e}_y} \nabla_{\mathbf{e}_y} n_y(0) \rangle_{2,t} =: C_y \geq 0$. In fact, as in [L] we assume that $C_y > 0$. Then the rest of the argument in [L] §4.1 goes through, and we obtain the same polynomially decay estimates for the flow near $y$.

To sum up, the asymptotic behaviors of the flows in this $S^1$-equivariant version is entirely identical with the non-equivariant version in [L].

## 4.3 The Fredholm framework

We now setup the Fredholm framework for describing the moduli spaces $\mathcal{M}_P, \mathcal{M}_O$.

### 4.3.1 The Fredholm framework for $\mathcal{M}_P$.

We shall need various versions of weighted Sobolev norms. Let's start by recalling some definitions. Let $\Theta := \mathbb{R} \times S_1^1$; and $\mathrm{pr}_1, \mathrm{pr}_2$ be the projection to the first and second component respectively. Let $\beta : \mathbb{R} \to [0,1]$ be a smooth cutoff function supported on $\mathbb{R}^+$ such that $\beta(s) = 1$ as $s \geq 1$. Given $\delta_-, \delta_+ \in \mathbb{R}$, we define the *weight function* $\varsigma^{(\delta_-,\delta_+)} : \mathbb{R} \to \mathbb{R}$ by

$$\varsigma^{(\delta_-,\delta_+)}(s) := e^{(s\beta(s)\delta_+ - s\beta(-s)\delta_-)}.$$

We'll often denote $\mathrm{pr}_1^* \varsigma^{(\delta_-,\delta_+)} : \Theta \to \mathbb{R}$ by $\varsigma^{(\delta_-,\delta_+)}$ as well. Suppose $P$ is an euclidean or hermitian vector bundle over $\Theta$. Given any section $\mu \in C_0^\infty(\Theta, P)$, Let

$$\|\mu\|_{p,k:(\delta_-,\delta_+)} := \|\varsigma^{(\delta_-,\delta_+)}\mu\|_{p,k};$$

let $L^p_{k:(\delta_-,\delta_+)}(P)$ denote the completion of $C_0^\infty$ with respect to the above norm.

Given two nondegenerate critical $S^1$-orbits $x, y \in \mathcal{P}$, the *configuration space* $\mathcal{B}_P(x,y)$ is defined as

$$\begin{aligned}
\mathcal{B}_P(x,y) := \Big\{ u \,\big|\, & u \in L^p_{1,loc}(\Theta, M), \\
& u(s,t) = \exp(y(t), \xi_+(s,t)) \text{ for some } \xi_+ \in L^p_{1:(0,\delta)}(\mathrm{pr}_2^*(y^*TM)) \text{ when } s > \rho_+(u); \\
& u(s,t) = \exp(x(t), \xi_-(s,t)) \text{ for some } \xi_- \in L^p_{1:(\delta,0)}(\mathrm{pr}_2^*(x^*TM)) \text{ when } s < \rho_-(u) \Big\},
\end{aligned} \tag{38}$$

where $\delta$ is as in §4.2.1; $\pm\rho_\pm(u)$ are large positive numbers depending on $u$. According to the discussion in §4.2.3, $\mathcal{M}_P(M, \gamma_0, x, y; J, X)$ is contained in $\mathcal{B}_P(x,y)$. $\mathcal{B}_P(x,y)$ is a Banach manifold modeled on $T_u \mathcal{B}_P(x,y) = L^p_{1:(\delta,\delta)}(\Theta, u^*TM)$.

For a $u \in \mathcal{M}_P(M, \gamma_0, x, y; J, X)$, let $E_u$ denote the linearization of $\bar{\partial}_{JX}$ at $u$ between $L^p_{1:(\delta,\delta)}(\Theta, u^*TM)$ and $L^p_{:(\delta,\delta)}(\Theta, u^*TM)$ It is by now standard that this is a Fredholm operator whose index is computed by the spectral flow from $\mathbb{A}_{x,[w]} + \delta$ to $\mathbb{A}_{y,[w\#u]} - \delta$, namely,

$$\mathrm{ind}\, E_u = \mathrm{ind}(x,[w]) - \mathrm{ind}(y,[w\#u]) - 1,$$

where $w\#u$ denotes the concatenation of $w$ and $u$. Also, note that the cokernel of $E_u$ always contains at least the 1-space spanned by $(\varsigma^{(\delta,\delta)})^{-2}\dot{u}$.

We also define the extended configuration spaces $\mathcal{B}_P(x,\underline{y}) = \bigcup_{y' \in \underline{y}} \mathcal{B}_P(x,y')$, $\mathcal{B}_P(\underline{x},y) = \bigcup_{x' \in \underline{x}} \mathcal{B}_P(x',y)$, $\mathcal{B}_P(\underline{x},\underline{y}) = \bigcup_{x' \in \underline{x}, y' \in \underline{y}} \mathcal{B}_P(x',y')$. They all inherit Banach manifold structures from $\mathcal{B}_P(x,y)$ via the fibrations

$$\mathcal{B}_P(x,y) \longrightarrow \mathcal{B}_P(\underline{x},y) \xrightarrow{\partial_-} S^1;$$
$$\mathcal{B}_P(x,y) \longrightarrow \mathcal{B}_P(x,\underline{y}) \xrightarrow{\partial_+} S^1;$$
$$\mathcal{B}_P(x,y) \longrightarrow \mathcal{B}_P(\underline{x},\underline{y}) \xrightarrow{\partial_- \times \partial_+} S^1 \times S^1,$$

$\partial_-, \partial_+$ being the end point maps.

The Banach spaces they model on are, respectively,

$$T_u\mathcal{B}_P(\underline{x},y) = L_{1:(\delta,\delta)}^{-;,p}(\Theta; u^*TM) = \mathbb{R} \oplus L_{1:(\delta,\delta)}^p(\Theta, u^*TM)$$
$$= \left\{C\beta(-s)\dot{u}(s,t) + \xi \,\middle|\, C \in \mathbb{R}, \xi \in L_{1:(\delta,\delta)}^p(\Theta, u^*TM)\right\};$$
$$T_u\mathcal{B}_P(x,\underline{y}) = L_{1:(\delta,\delta)}^{+;,p}(\Theta; u^*TM) = \mathbb{R} \oplus L_{1:(\delta,\delta)}^p(\Theta, u^*TM)$$
$$= \left\{C\beta(s)\dot{u}(s,t) + \xi \,\middle|\, C \in \mathbb{R}, \xi \in L_{1:(\delta,\delta)}^p(\Theta, u^*TM)\right\};$$
$$T_u\mathcal{B}_P(\underline{x},\underline{y}) = L_{1:(\delta,\delta)}^{-+,p}(\Theta; u^*TM) = \mathbb{R} \oplus \mathbb{R} \oplus L_{1:(\delta,\delta)}^p(\Theta, u^*TM)$$
$$= \left\{(C_-\beta(-s) + C_+\beta(s))\dot{u}(s,t) + \xi \,\middle|\, C_\pm \in \mathbb{R}, \xi \in L_{1:(\delta,\delta)}^p(\Theta, u^*TM)\right\}.$$

Let $E_u^{-;}, E_u^{;+}, E_u^{-+}$ be respectively the Fredholm operators that extend the domain of $E_u$ to these spaces. Note that since $E_u^{-+}(\chi(s)\dot{u}) = \chi'(s)\dot{u}$ for any real-valued function $\chi(s)$, these operators all have cokernels of dimension $\dim \operatorname{coker} E_u - 1$.

Let

$$\mathcal{M}_P(x,y) = \mathcal{M}_P(M, \gamma_0, x, y; J, X) := \mathcal{M}_P(M, \gamma_0; J, X) \cap \mathcal{B}_P(x,y),$$
$$\mathcal{M}_P(\underline{x},\underline{y}) = \mathcal{M}_P(M, \gamma_0, \underline{x}, \underline{y}; J, X) := \bigcup_{x \in \underline{x}, y \in \underline{y}} \mathcal{M}_P(x,y),$$

and similarly for other variants of $\mathcal{M}_P$.

$E_u, E_u^{;+}, E_u^{-;}, E_u^{-+}$ are respectively the deformation operators of $\mathcal{M}_P(x,y), \mathcal{M}_P(x,\underline{y})$, $\mathcal{M}_P(\underline{x},y)$, and $\mathcal{M}_P(\underline{x},\underline{y})$.

Note that the orbits of the $S_t^1$ action map to the diagonal under $\partial_- \times \partial_+ : \mathcal{M}_P(M, \gamma_0, \underline{x}, \underline{y}; J, X) \to S^1 \times S^1$, and hence

$$\underline{\mathcal{M}}_P(\underline{x},\underline{y}) := \mathcal{M}_P(\underline{x},\underline{y})/S_t^1 = \mathcal{M}_P(\underline{x},y) = \mathcal{M}_P(x,\underline{y}).$$

### 4.3.2 The Fredholm framework for $\mathcal{M}_O$.

Similarly one can define the *configuration space* $\mathcal{B}_O$ in which $\mathcal{M}_O(M, \gamma_0; J, X)$ embeds. Let $\Theta_{T,Q}$ be as in (36).

$$\mathcal{B}_O := \left\{(1/T, Q, u) \,\middle|\, T \in \mathbb{R}^+, Q \in S_1^1; u \in L_1^p(\Theta_{T,Q}, M)\right\}.$$

It is endowed with a Banach manifold structure via the fibration

$$\mathcal{B}_O \to \mathbb{R}^+ \times S^1$$

to the parameter space of $(1/T, Q)$; namely, $T_{(1/T,Q,u)}\mathcal{B}_O = \mathbb{R}_\varrho \oplus \mathbb{R}_q \oplus L_1^p(\Theta_{T,Q}, u^*TM)$, where $\mathbb{R}_\varrho, \mathbb{R}_q$ parameterize the variations in $1/T, Q$ respectively. We shall denote by

$$\tilde{D}_{(1/T,Q,u)} : \mathbb{R}_\varrho \oplus \mathbb{R}_q \oplus L_1^p(\Theta_{T,Q}, u^*TM) \to L^p(\Theta_{T,Q}, u^*TM)$$

the linearization of $\bar{\partial}_{JX}$ at $(1/T, Q, u) \in \mathcal{B}_O$, and by

$$D_u : L_1^p(\Theta_{T,Q}, u^*TM) \to L^p(\Theta_{T,Q}, u^*TM)$$

its restriction. One may re-express $(s, t)$ in terms of new variables $(s', t')$ which take values in $\Theta_{1,1}$:

$$(s, t) = (0, t') + (s'T, s'Q).$$

Substituting in the new variables, the flow equation (24) becomes:

$$\partial_s u = 1/T \partial_{s'} u - Q/T \partial_{t'} u; \quad \partial_t u = \partial_{t'} u.$$

Thus,

$$\tilde{D}_{(1/T,Q,u)}(\varrho, q, \xi) = D_u \xi + \varrho T \partial_s u - q/T \partial_t u.$$

Note that $\dot{u}, u'$ are in the cokernel of $D_u$, but not in $\tilde{D}_{(1/T,Q,u)}$.

$D_u$ and $\tilde{D}_{(1/T,Q,u)}$ are both Fredholm operators, since the first is an elliptic differential operator on a compact space $\Theta_{T,Q}$, and the second is just a finite dimensional extension of the first. By Riemann-Roch,

$$\operatorname{ind} D_u = c_1(TM)[u]; \quad \operatorname{ind} \tilde{D}_{(1/T,Q,u)} = 2 + c_1(TM)[u],$$

where $[u] \in H_2(M; \mathbb{Z})$ is the homology class of $u : \Theta_{T,Q} \to M$.

It is also convenient to define $\mathcal{D}_{(1/T,Q,u)} : \mathbb{R}_\varrho \oplus \mathbb{R}_q \oplus L_1^p(u^*TM) \longrightarrow \mathbb{R} \oplus \mathbb{R} \oplus L^p(u^*TM)$,

$$\mathcal{D}_{(1/T,Q,u)} = d_u^{s*} + d_u^{t*} + \tilde{D}_{(1/T,Q,u)} = \begin{pmatrix} 0 & 0 & \Pi_{u'} \\ 0 & 0 & \Pi_{\dot{u}} \\ Tu' & -1/T\dot{u} & D_u \end{pmatrix},$$

where $\Pi_{u'}, \Pi_{\dot{u}}$ denote $L^2$-orthogonal projection to $u', \dot{u}$ respectively, and $d_u^s, d_u^t$ are respectively the linearization of the $S_s^1$, $S_t^1$ action.

### 4.3.3 Fredholm context for $\mathcal{M}_P(x,y)$, when one of $x, y$ is minimally degenerate.

In the invariance proof in section 7 below, we shall also need to consider the case when some of the critical $S^1$-orbits are minimally degenerate. What follows is a modification of the work in [L] §4.2. Without loss of generality, we assume that $x \in \mathcal{P}$ is in a nondegenerate critical $S^1$ orbit, and $y \in \mathcal{P}$ is in a minimally degenerate critical $S^1$-orbit. In this case we need polynomially weighted norms instead of exponentially weighted ones.

We first recall the following definitions from [L]. Let $u \in \mathcal{M}_P(x, y)$, and $\xi \in u^*TM$. The "longitudinal" component of $\xi$ is

$$\xi_L(s) := \|u'(s)\|_{2,t}^{-2} \langle u'(s), \xi(s) \rangle_{2,t} u'(s).$$

*Notation.* The subscript $t$ in $\|\cdot\|_{2,t}$ or $L_t^2$ emphasizes that the $L^2$-norm is taken with respect to the variable $t$.

Let $\xi_{\mathfrak{f}}(s)$ denote the $L_t^2$-orthogonal projection to $\dot{u}(s)$, and define

$$\xi_T := \xi - \xi_L - \xi_{\mathfrak{f}}.$$

Note that when $u \in \mathcal{M}_P(x,y)$, $u', \dot{u}$ are $L_t^2$-orthogonal by the $S^1$-invariance of $\tilde{\mathcal{A}}_X$.

Define the *weight*

$$\sigma_u(s) = \begin{cases} \|u'(0)\|_{2,t}^{-1} & \text{for } s \leq 0, \\ \|u'(s)\|_{2,t}^{-1} & \text{for } s \geq 0. \end{cases}$$

From §4.2.3 and [L] §4.1 we know that there are positive constants $C', C$ such that $C'/s^2 \geq \sigma_u \geq C/s^2$ for large $s$.

**Definition.** Let $u, x, y$ be as above, and $\delta > 0$ be as in (38). We define the following norms on $C_0^\infty(u^*TM)$:

$$\|\xi\|_{W_u} := \|\xi\|_{p,1:(\delta,0)} + \|\sigma_u^{1/2}\xi\|_{p,1} + \|\sigma_u(\xi_L' + \xi_{\mathfrak{f}}')\|_p;$$
$$\|\xi\|_{L_u} := \|\xi\|_{p:(\delta,0)} + \|\sigma_u^{1/2}\xi\|_p + \|\sigma_u(\xi_L + \xi_{\mathfrak{f}})\|_p.$$

Let $W_u, L_u$ denote the completion of $C_0^\infty$ with respect to the above two norms respectively.

**Lemma.** $E_u$ *is a bounded Fredholm operator between $W_u$ and $L_u$ of index* $\text{ind}(x, [w]) - \text{ind}(y, [w\#u])$.

*Proof.* The boundedness follows from routine computations.

The Fredholmness basically follows from the nice fact that

$$E_u \text{ is 'almost diagonal' with respect to the decomposition } \xi = \xi_{\mathfrak{f}} + \xi_{\mathfrak{n}} \text{ when } s \to \pm\infty. \tag{39}$$

To see this, note that the $\mathfrak{fn}$ component of $E_u$ vanishes since by the $S^1$-invariance of $\tilde{\mathcal{A}}_X$,

$$E_u \dot{u} = 0. \tag{40}$$

As for the $\mathfrak{nf}$ component, using the fact that $A_{u(s)}$ is self-adjoint (this can be checked by direct computation), (40), and the fact that $\langle \dot{u}(s), \xi_{\mathfrak{n}}(s)\rangle_{2,t} = 0 \ \forall s$, we have

$$\langle \dot{u}, E_u \xi_{\mathfrak{n}}\rangle_{2,t} = 2\langle \dot{u}, A_u \xi_{\mathfrak{n}}\rangle_{2,t} = 2\langle \nabla_t \mathcal{V}_X(u), \xi_{\mathfrak{n}}\rangle_{2,t}, \tag{41}$$

which is small ($\sim s^{-2}\|\xi_{\mathfrak{n}}\|_{2,t}$ by §4.2.3 and [L] §4.1) when $u(s)$ is close to the critical points. The second equality in (41) follows from linearizing $\langle \dot{u}, \mathcal{V}_X(u)\rangle_{2,t} = 0$, which in turn follows from the $\gamma_0$-exactness of $\theta_X$. Thus the $\mathfrak{nf}$ component is small when $s \to \pm\infty$. By a typical excision argument (see e.g. [F88a]), one may reduce the proof for Fredholmness to the case when $u$ is close to $y$. Now by (39) we may ignore the off-diagonal terms and separate the variables; the Fredholmness of $E_u$ then follows from the Fredholmness of $E_{\mathfrak{nn}}$ and $E_{\mathfrak{ff}}$: The Fredholmness of the former is proved in [L]; the Fredholmness of the latter follows from the fact that $d/ds$ is an isomorphism between $N_1$, $N_2$, which are respectively the completion of $C_0^\infty(\mathbb{R}^+)$ with respect to the norms:

$$\|f\|_{N_1} := \|\sigma_u^{1/2} f\|_p + \|\sigma_u f'\|_p; \quad \|f\|_{N_2} := \|\sigma_u f\|_p.$$

To compute the index, one needs a decay estimate for $\xi_{\mathfrak{f}} \,\forall \xi \in \ker E_u$, in addition to the work of [L] §4.2. This can be done by observing that writing $\xi_{\mathfrak{f}} = g(s)\dot{u}$, the function $g(s)$ satisfies the following ODE:

$$g'\|\dot{u}\|_{2,t}^2 + \langle \dot{u}, E_u \xi_T \rangle_{2,t} = 0.$$

$\|\xi_T\|_{2,t}$ is estimated in [L] to decay exponentially, with an exponent larger than $\delta$. On the other hand, $\|\dot{u}(s)\|_{2,t} \geq C'' s^{-1}$ according to [L]. Thus $g(s)$ also decays exponentially. Combining with the work of [L], one sees that $\mathrm{ind}\, E_u = \mathrm{CZ}(\mathbb{A}_{x,[w]} + \delta) - \mathrm{CZ}(\mathbb{A}_{y,[w\#u]} - \delta\Pi_{\mathfrak{f}} + \delta\Pi_L)$, where $\Pi_{\mathfrak{f}}, \Pi_L$ denote $L_t^2$-orthogonal projection to $\dot{u}$ and $u'$ respectively. □

Of course, all the Fredholm operators introduced in this paper depend implicitly on the pair $(J, X)$. We shall add $J, X$ in the subscripts when we wish to emphasize this dependence.

### 4.4 Structure of the moduli spaces

The purpose of this subsection is to establish Proposition 4.4.2 below about the smoothness and compactness of the relevant moduli spaces.

**4.4.1 The spaces $\mathcal{X}_{ad}$ and $\mathcal{J}_{ad}$.** The argument in this subsection is built on the work of [FHS]; we shall therefore begin by recalling some definitions in [FHS].

According to [FHS] Definition 7.1, a $X \in \mathcal{X}(\theta)$ is *admissible* if the zeros of $X$ are all nondegenerate, and if there is an almost complex structure $J \in \mathcal{J}$ such that for each $p \in X^{-1}(0)$ and each unitary frame $\Phi_p : \mathbb{R}^{2n} \to T_p M$ taking $J(p)$ and $\omega(p)$ to the standard $J_0, \omega_0$, the symmetric matrix $J_0 \Phi_p^{-1} dX(p) \Phi_p$ is "regular" in the sense of [FHS] section 6. Such an almost complex structure is said to be *admissible* (with respect to $X$). (Note that regular symmetric matrices are generic). The set of admissible symplectic vector fields in $\mathcal{X}(\theta)$ is denoted $\mathcal{X}_{ad}(\theta)$. For a fixed $X \in \mathcal{X}_{ad}(\theta)$, $\mathcal{J}_{ad} = \mathcal{J}_{ad}(X) \subset \mathcal{J}$ denotes the corresponding space of admissible almost complex structures.

**Lemma.** ([FHS] Lemma 7.2) $\mathcal{X}_{ad}(\theta)$ *is an open dense set in* $\mathcal{X}(\theta)$, *and for every* $X \in \mathcal{X}_{ad}(\theta)$, $\mathcal{J}_{ad}(X)$ *is open and dense in* $\mathcal{J}$.

The salient point about admissible symplectic vector fields and almost complex structures is the next Lemma. First we need to introduce more definitions.

Given $(J, X) \in \mathcal{J} \times \mathcal{X}$, we denote $V = -JX$, and let $M'_X \subset M$ be the complement of $X^{-1}(0)$. Let the metric on $M$ be that induced by $J$, so that $X, V$ are orthogonal to each other. Over $M'_X$, $X, V$ together span a 2-plane field, which we denote by $F_{J|X}$. Let $F_{J|X}^\perp$ denote the $(2n-2)$-plane field of the orthogonal complement to $F_{J|X}$. We shall also use the same notation to denote the corresponding rank 2 and rank $(2n-2)$ subbundles of $TM$. Let $\Pi_{F_{J|X}^\perp} : TM_{X'} \to F_{J|X}^\perp$ denote the orthogonal projection to $F_{J|X}^\perp$.

We define two sections of $F_{J|X}^\perp$, $\mathfrak{c}_1(J, X)$ and $\mathfrak{c}_2(J, X)$:

$$\mathfrak{c}_1(J, X) := \Pi_{F_{J|X}^\perp}[X, V],$$

$$\mathfrak{c}_2(J, X) := \Pi_{F_{J|X}^\perp}(\nabla_V[X, V] - \mathbf{e}_X \cdot [X, V]\nabla_V X - \mathbf{e}_V \cdot [X, V]\nabla_V V),$$

where $\mathbf{e}_V = V/|V|$, $\mathbf{e}_X = X/|X|$ are unit vector fields, and $\nabla$ is the Levi-Civita connection. Let $\mathfrak{c}(J, X) : M_{X'} \to F_{J|X}^\perp \otimes \mathbb{R}^2$ be

$$\mathfrak{c}(J, X) = \mathfrak{c}_1(J, X) \oplus \mathfrak{c}_2(J, X).$$

The following is a rephrasing of [FHS] Lemma 7.5:

**Lemma.** ([FHS] Lemma 7.5) *Let $X \in \mathcal{X}_{ad}$, and $J \in \mathcal{J}_{ad}(X)$. Then there exists a neighborhood $U(J, X)$ of $X^{-1}(0)$, such that $\mathfrak{c}(J, X)^{-1}(0) \cap U(J, X) = \emptyset$.*

Now we are almost ready to state the

**4.4.2 Structure Theorem for $\mathcal{M}_O$ and $\mathcal{M}_P$.** Fix $X \in \mathcal{X}_{ad}(\theta) \cap \mathcal{X}_{nondeg}(\theta)$, and $J_0 \in \mathcal{J}_{ad}(X)$. (Note that $\mathcal{X}_{ad}(\theta) \cap \mathcal{X}_{nondeg}(\theta)$ is open dense in $\mathcal{X}(\theta)$ according to Lemma 4.2.2 and [FHS] Lemma 7.2). Let $U(J_0, X)$ be as in the previous Lemma. Let $U_0 = U_0(J_0, X) \subset U(J_0, X)$ be a very small tubular neighborhood of $X^{-1}(0)$ with $|c(J_0, X)|(x) \geq \varepsilon > 0 \; \forall x \in \partial U_0(J_0, X)$. (Such $U_0$ exists by the proof of [FHS] Lemma 7.5). Let $\mathcal{J}_{J_0}(U_0) \subset \mathcal{J}$ be the Banach manifold of all $J$ which agrees with $J_0$ over $\bar{U}_0$.

**Proposition.** *Let $M, \gamma_0, \theta$ be as in §2.3 and $X \in \mathcal{X}_{nondeg}(\theta) \cap \mathcal{X}_{ad}(\theta)$, $J_0 \in \mathcal{J}_{ad}(X)$, $U_0$ be as above. Then there is a Baire set $\mathcal{J}_{reg} \subset \mathcal{J}_{J_0}(U_0)$ such that for any $J \in \mathcal{J}_{reg}$, $(J, X)$ is regular in the sense of §2.1.2 (with all the moduli spaces there replaced by $S_t^1$-quotients).*

*Proof.* Compactness is standard due to the assumption that $M$ is monotone; the only difference from the standard literature is that to obtain the description of $\partial \hat{\mathcal{M}}_P^2$ in terms of broken flow lines, one needs a fibered version of gluing theory instead. With this, one shows that

$$\partial \overline{\mathcal{M}_P^2(\underline{x}, \underline{y})}/\mathbb{R}_s = \coprod_{\underline{z} \in \mathcal{P}} \mathcal{M}_P^1(\underline{x}, \underline{z})/\mathbb{R}_s \times_{S_t^1} \mathcal{M}_P^1(\underline{z}, \underline{y})/\mathbb{R}_s, \quad (42)$$

where the fiber produce is taken with respect to the fibrations

$$\mathcal{M}_P^1(\underline{x}, \underline{z})/\mathbb{R}_s \xrightarrow{\partial_+} S_t^1; \quad \mathcal{M}_P^1(\underline{z}, \underline{y})/\mathbb{R}_s \xrightarrow{\partial_-} S_t^1.$$

This implies the desired

$$\partial \overline{\hat{\mathcal{M}}_P^1(\underline{x}, \underline{y})} = \coprod_{\underline{z}\in\mathcal{P}} \hat{\mathcal{M}}_P^0(\underline{x}, \underline{z}) \times \hat{\mathcal{M}}_P^0(\underline{z}, \underline{y}).$$

We shall henceforth concentrate on transversality.

Transversality is also standard except for an injectivity result (Lemma 4.4.6 below), which we postpone until the next subsubsection. Such an injectivity result holds for an open dense set $\mathcal{J}_{good} \subset \mathcal{J}_{J_0}(U_0)$. Using the Banach manifold $\mathcal{J}_{good}$ instead of $\mathcal{J}_{J_0}(U_0)$ in the usual Sard-Smale argument, Lemma 4.4.6 then implies that one may perturb $J$ away from $U_0$ to achieve transversality. (Note that if $U_0$ is sufficiently small, all $u \in \mathcal{M}_P$ or $\mathcal{M}_O^s \subset \mathcal{M}_O$ (the subset of simple elements) must pass through $M\backslash U_0$, since the minimal energy of a $u \in \mathcal{M}_O$ or $\mathcal{M}_P(x,x)$ is determined by the periods of the 1-form $\mathcal{Y}_X$.) For details, see [FHS] Theorem 7.4 and [L] sections 3.2–3.3. Here we only remark that for the case of multiple-covered tori, the arguments in [L] §3.3 need the following minor modifications. Let $u^m$ be the $m$-cover of the simple curve $(1/T, Q, u) \in \mathcal{M}_O^s$. Following the notations in [L] §3.3, in our case $(\ker {}^m D_{u,J|X})^\perp$, $(\text{coker } {}^m D_{u,J|X})^\perp$ should be the orthogonal complements of the 2-space spanned by $\{i_m(\partial_s u^m), i_m(\partial_t u^m)\}$ in $\ker {}^m D_{(1/T,Q,u),J|X}$, $\text{coker } {}^m D_{(1/T,Q,u),J|X}$ respectively; then the same arguments show that

$$(\ker {}^m D_{(u,J|X)})^\perp = \ker {}^m \mathcal{D}_{(1/T,Q,u),J|X},$$
$$(\text{coker } {}^m D_{(u,J|X)})^\perp = \text{coker } {}^m \mathcal{D}_{(1/T,Q,u),J|X}.$$

The other difference is that in this paper we perturb $J$ instead of $H$; thus $\nabla\nabla h$ in the proof of [L] Lemma 3.11 should be replaced by $\nabla j \dot{u}$. It is just as easy to come up with a $j$ satisfying the required properties. □

**Remark.** In this version of Floer theory we fix $X$ and perturb $J$. Since varying $J$ also changes the metric on the loop space, the variation of $\mathcal{V}_X$ might not always be perpendicular to the cokernel of the deformation operator.

We now turn to the missing part of the above proof, namely the relevant injectivity result. First we introduce the notion of good almost complex structures, which refines the notion of admissibility in [FHS].

**4.4.3 Definition.** Given $X \in \mathcal{X}_{ad}$, $J_0 \in \mathcal{J}_{ad}$, $U_0 = U_0(J_0, X)$ as in §4.4.2. An almost complex structure $J \in \mathcal{J}_{J_0}(U_0)$ is said to be *good* if the section $\mathfrak{c}(J, X) : M_X' \to F_{J|X}^\perp \otimes \mathbb{R}^2$ intersects the zero section transversely. The set of all good almost complex structures in $\mathcal{J}_{J_0}(U_0)$ is denoted by $\mathcal{J}_{good}$.

Note that by our definitions and choice of $U_0$, the zero locus $\mathfrak{c}(J,X)^{-1}(0)$ is a compact submanifold in $M'_X$ of dimension $4-2n$.

**4.4.4 Lemma.** *Let $X, J_0, U_0$ be as in §4.4.2. Then $\mathcal{J}_{good} \subset \mathcal{J}_{J_0}(U_0)$ is open and dense.*

*Proof.* This follows easily from a standard transversality argument. Let
$$\tilde{M}'_X(\mathcal{J}_{J_0}(U_0)) := \mathcal{J}_{J_0}(U_0) \times M'_X,$$
and let $F^\perp(\mathcal{J}_{J_0}(U_0), X)$ be the rank $2n-2$ euclidean bundle over $\tilde{M}'_X(\mathcal{J}_{J_0}(U_0))$ which restricts to $F^\perp_{J|X}$ over $\{J\} \times M'_X$ $\forall J \in \mathcal{J}_{J_0}(U_0)$. Similarly, let $\mathfrak{c}(\mathcal{J}_{J_0}(U_0), X)$ be the section of $F^\perp(\mathcal{J}_{J_0}(U_0), X)$ which restricts to $\mathfrak{c}(J,X)$ over $\{J\} \times M'_X$. It is easy to see that $\mathfrak{c}(\mathcal{J}_{J_0}(U_0), X)$ is transverse to the zero section: For example, let $(J,x) \in \tilde{M}'_X(\mathcal{J}_{J_0}(U_0))$ be a zero of $\mathfrak{c}(\mathcal{J}_{J_0}(U_0), X)$, and let $V := -JX$. Our choices of $X$ and $U_0$ imply that $x \in M\setminus\bar{U}_0$ and $X(x) \neq 0$. Consider a variation $\delta J$ of $J$ such that $\delta J(x) = 0$, and let $X$ be fixed. This guarantees that $X(x), V(x), F_{J|X}(x), F^\perp_{J|X}(x)$ remain unchanged, and
$$\delta[X,V](x) = \nabla_X(\delta J)X(x), \quad \delta\nabla_V V(x) = (\nabla_V \delta J)X(x).$$

By choosing the first order term in the Taylor expansion of $\delta J$ at $x$ appropriately, we may arrange such that:
$$\delta\mathfrak{c}_1(x) = \delta[X,V](x) = w,$$
$$\Pi_{F^\perp_{J|X}} \delta\nabla_V V(x) = 0,$$
where $w \neq 0$ is an arbitrary vector orthogonal to $X(x), V(x)$. Notice that the first condition is possible because of our assumption that $\dim M \geq 4$. With such choice of $\delta J$,
$$\delta\mathfrak{c}_2(x) = \Pi_{F^\perp_{J|X}} \nabla_V(\nabla_X(\delta J)X)(x),$$
which may be an arbitrary vector in $F^\perp_{J|X}$ by a suitable choice of the second order term of $\delta J$.

We have now verified the transversality of $\mathfrak{c}(\mathcal{J}_{J_0}(U_0), X)$ to the zero section, and thus $\mathfrak{c}(\mathcal{J}_{J_0}(U_0), X)^{-1}(0)$ is a codimension $4n-4$ submanifold of $\tilde{M}'_X(\mathcal{J}_{J_0}(U_0))$. By the Sard-Smale theorem, the set of regular values of the projection $\mathfrak{c}(\mathcal{J}_{J_0}(U_0), X)^{-1}(0)$ to $\mathcal{J}_{J_0}(U_0)$ is Baire. At such a regular value $J$, $\mathfrak{c}(J,X)^{-1}(0) \subset M'_X$ is a submanifold of dimension $2n - (4n-4) = 4 - 2n$. In fact, it is either empty or consists of finitely many points, since $M$ is compact of dimension $\geq 4$, and by assumption the zero locus avoids $U_0$. Thus in fact the set of regular values is not only Baire, but open and dense. □

The next lemma illuminates the significance of good almost complex structures.

**4.4.5 Lemma.** *Let $X, J_0, U_0$ be as in Definition 4.4.3, and let $J \in \mathcal{J}_{good}$. Then for any $u \in \mathcal{M}_P(M, \gamma_0; J, X)$ or $\mathcal{M}_O^s(M, \gamma_0; J, X)$, the four vectors $\partial_s u, \partial_t u, X, JX$ are linearly independent on an open and dense subset of $\Theta$, where $\Theta$ is the domain of $u$; namely, $\Theta = \mathbb{R} \times S_1^1$ when $u \in \mathcal{M}_P$, and $\Theta = \Theta_{T,Q}$ when $(1/T, Q, u) \in \mathcal{M}_O^s$.*

*Proof.* The openness is obvious. We shall therefore concentrate on denseness.

We saw in the proof of Lemma 4.4.4 that

$$\mathfrak{z} := X^{-1}(0) \cup \mathfrak{c}(J, X)^{-1}(0) \subset M$$

consists of at most finitely many points. Since $u$ is nonconstant, by the unique continuation theorem (cf. e.g. [FHS] Proposition 3.1):

$$\Theta \setminus u^{-1}(\mathfrak{z}) \text{ is open and dense.} \tag{43}$$

On the other hand, we have the following standard fact:

**Sublemma.** *Let $J, X, u, \Theta$ be as in Lemma 4.4.5. Then $\partial_s u, \partial_t u$ are linearly independent on an open and dense set $R(u) \subset \Theta$.*

*Proof.* This follows, for instance, from a simple modification of the proof of [FHS] Lemma 7.6: in our case, $u$ can not be $t$-independent because of homological constraint $[\gamma_0] \neq 0$; in fact, because of the primitivity of $[\gamma_0]$,

$$\forall m \in \mathbb{Z}^+, \exists (s,t) \in \Theta \text{ such that } u(s,t) \neq u(s, t+1/m).$$

[FHS] actually called this condition simple, which is different from our definition of simple curves. Instead, we shall refer to this as "FHS-simple". On the other hand, $u$ can not be $s$-independent either, since by assumption $u$ is nonconstant. Finally, $u$ can not be of the form $u(s,t) = \nu(s + \lambda t)$: in the case when $u \in \mathcal{M}_O$, this is clear from the homological constraint $[u] \neq 0 \in H_2(M)$; in the case when $u \in \mathcal{M}_P$, it follows from the fact that $u$ is FHS-simple and nonconstant as in [FHS]. The rest of the proof is identical to that of [FHS]. $\square$

To continue the proof of Lemma 4.4.5, suppose to the contrary that $\partial_s u, \partial_t u, X, V$ are linearly dependent on a neighborhood $N \subset \Theta$. By (43) and the above sublemma, we may assume without loss of generality that $\partial_s u$ and $\partial_t u$ are linearly independent over $N$, and that $N \cap u^{-1}(\mathfrak{z}) = \emptyset$. This means that $\partial_s u, \partial_t u$ must span the same 2-plane in $T_{u(s,t)}M$ as $X(u(s,t)), V(u(s,t))$ over $N$. In particular, the two-plane field spanned by $X, V$ must be integrable at $u(s,t)$; i.e. $\Pi_{F_{J|X}^\perp}[X, V](u(s,t)) = 0$. Now taking the covariant derivative of this equation with respect to $V$, and plugging in the flow equation (24), one sees that $J, X$ also have to satisfy

$$\Pi_{F_{J|X}^\perp}(\nabla_V [X,V] - \mathbf{e}_X \cdot [X,V]\nabla_V X - \mathbf{e}_V \cdot [X,V]\nabla_V V)(x)) = 0 \quad \text{over } u(N).$$

(See the proof of [FHS] Lemma 7.7 for the explicit computation.) In other words, $u(N) \in \mathfrak{z}$. This contradicts the assumption on $N$ that we started out with. $\square$

The following lemma is an immediate corollary of the above lemma and the argument of [FHS] Lemma 7.8:

**4.4.6 Lemma.** *Under the assumptions of Lemma 4.4.5, the set of injective points*

$$\left\{(s,t) \,|\, (s,t) \in R(u),\, u^{-1}\{u(s,t)\} = \{(s,t)\}\right\} \subset \Theta$$

*is open and dense in $\Theta$.*

Very roughly, the argument shows that if this were not true, the linear independence of $\partial_s u, \partial_t u, X, V$ on an open set implies a periodicity condition $u(s,t) = u(s+a, t+b)$. However, either when $u \in \mathcal{M}_P$ or $\mathcal{M}_O^s$, the simplicity and non-constancy of $u$ means $a + ib \in \Gamma$, where $\Gamma$ is the lattice such that $\Theta = \mathbb{C}/\Gamma$.

## 4.5 Orientations

Orientations of $\mathcal{M}_O, \mathcal{M}_P$ induce orientations on their $S_t^1$ quotients; we now specify the orientations of these moduli spaces. In other words, we need to specify the orientations of $\det \tilde{D}_u, \det E_u^{\cdot-}, \det E_u^{\cdot+}, \det E_u^{-+}$ for $u \in \mathcal{M}_O$ or $\mathcal{M}_P$. We shall orient the determinant lines of $D_u, E_u$ instead; by regarding $\tilde{D}_u, E_u^{\cdot-}, E_u^{\cdot+}, E_u^{-+}$ as stabilizations (i.e. finite dimensional extensions) of $D_u$ and $E_u$, these induce orientations on $\det \tilde{D}_u$, $\det E_u^{\cdot-}, \det E_u^{\cdot+}, \det E_u^{-+}$.

We follow the sign and orientation conventions in [L] §5.1.1; in addition, we shall use the following convention for the fiber product: Given two oriented real vector spaces fibering over $\mathbb{R}$, $V \xrightarrow{\partial_+} \mathbb{R}$ and $V' \xrightarrow{\partial_-} \mathbb{R}$, $V \times_\mathbb{R} V'$ is oriented by regarding it as the inverse image of the diagonal under the map $\partial_+ \times \partial_-$:

$$\begin{array}{ccccccccc}
0 & \longrightarrow & V' & \longrightarrow & V \times V' & \longrightarrow & V & \longrightarrow & 0 \\
& & \downarrow & & \downarrow{\partial_+ \times \partial_-} & & \downarrow & & \\
0 & \longrightarrow & \mathbb{R} & \longrightarrow & \mathbb{R} \times \mathbb{R} & \longrightarrow & \mathbb{R} & \longrightarrow & 0.
\end{array}$$

We fix a pair of regular $(J, X)$ throughout this subsection.

**(a) Orienting $\mathcal{M}_O$.**

In this version, these orientations may be defined basically in the same way as [L] §5.2. We again restrict our discussions to $\mathcal{M}_O^2$ since this is the case we need. The argument is however easily generalizable to higher dimensions. In this case, the relevant space of operators is

$$\Sigma_O := \Big\{ (1/T, Q, D) \,\Big|\, (1/T, Q) \in \mathbb{R}^+ \times S_1^1,$$
$$D = \partial_s + j(s,t)\partial_t + \nu(s,t) : L_1^p(\mathbb{C}^n) \to L^p(\mathbb{C}^n),\, (s,t) \in \Theta_{T,Q} \Big\},$$

where $\mathbb{C}^n$ denotes the trivial bundle over $\Theta_{T,Q}$; $j, \nu$ are as in §4.2.3 (a). ($\mathbb{C}^n$ should be replaced by a nontrivial bundle for higher dimensional moduli spaces).

$\Sigma_O$ contracts to the subset

$$\Sigma_O^c = \{(1/T, Q, D) \mid D \text{ is complex-linear}\} \subset \Sigma_O,$$

over which the determinant line bundle $\mathcal{L}_O$ is canonically oriented by the complex linearity; which also orients $\mathcal{L}_O$ over $\Sigma_O$.

Given $(1/T, Q, u) \in \mathcal{M}_O^2$, since $u^*TM$ is trivial, $D_u$ may be identified with an element in $\Sigma_O$, and $\det D_u$ is accordingly oriented. (This does not depend on the choice of trivialization of $u^*TM$, since the notion of complex linearity is independent of trivializations).

Since the deformation operator $\tilde{D}_{(1/T, Q, u)}$ is a stabilization of $D_u$, $\mathcal{M}_O^2$ is thus oriented. This induces orientations of the quotients $\mathcal{M}_O^2/S_t^1$, $\hat{\mathcal{M}}_O^0$.

### (b) Orienting $\mathcal{M}_P$.

We shall orient the quotient spaces

$$\mathcal{M}_P(\underline{x}, \underline{y})/S_t^1 = \mathcal{M}_P(x, \underline{y}) = \mathcal{M}_P(\underline{x}, y),$$

and endow $\mathcal{M}_P(x, y)$ with the induced orientation.

Following the framework of [FH], we need to address two main components of the argument: First in §4.5.1 we shall establish the orientability of the determinant line bundles $\mathcal{L}_P(\underline{x}, [w]; \underline{y}, [v])$ over the relevant space of operators $\Sigma_P(\underline{x}, [w]; \underline{y}, [v])$ for each pair of $(\underline{x}, [w]), \underline{y}, [v]) \in \tilde{\mathcal{P}}$. Thus the space of possible orientations of $\mathcal{L}_P(\underline{x}, [w]; \underline{y}, [v])$, $\text{or}(\underline{x}, [w]; \underline{y}, [v]) = \mathbb{Z}/2\mathbb{Z}$, is well-defined. Then in §4.5.2 we shall prove a (linear) gluing lemma, which is useful for defining a splicing homomorphism:

$$\# : \text{or}(\underline{x}, [w]; \underline{y}, [v]) \times \text{or}(\underline{y}, [v]; \underline{z}, [r]) \to \text{or}(\underline{x}, [w]; \underline{z}, [r]).$$

A coherent orientation is a systematic way of choosing orientations of all such $\mathcal{L}_P$'s in a way compatible with splicing.

**4.5.1 The spaces of operators $\Sigma_P$.** Recall from §4.2.1 that a trivialization of $\gamma_0^* TM$ induces a (class of) trivialization $\Phi_{x,[w]}$ of $x^*TM$, for each $(x, [w]) \in \tilde{\mathcal{P}}(X)$. (In fact, if $(x, [w]), (x', [w']) \in (\underline{x}, [w])$ is in the same $\mathbb{R}$-orbit in $\tilde{\mathcal{P}}$, then $\Phi_{x,[w]}$, $\Phi_{x',[w']}$ are related by rotation). With these trivializations chosen, each $(x, [w]) \in \tilde{\mathcal{P}}(X)$ gives a self-adjoint differential operator $\mathbb{A}_{(x,[w])} := \Phi_{x,[w]} A_x \Phi_{x,[w]}^{-1}$. Given $(x, [w]), (y, [v]) \in \tilde{\mathcal{P}}$, let

$$\Sigma_P(x, [w]; y, [v]) := \Big\{ E = \partial_s + j(s,t)\partial_t + \nu(s,t) : L_{1:(\delta,\delta)}^p(\mathbb{R}^{2n}) \to L_{:(\delta,\delta)}^p(\mathbb{R}^{2n}) \Big|$$
$$j \in C_{\epsilon:(\delta,\delta)} \ \nu \in C_{\epsilon:(\delta,\delta)}, \ (s,t) \in \mathbb{R} \times S_1^1,$$
$$E_l = \mathbb{A}_{(x,[w])}; \ E_r = \mathbb{A}_{(y,[v])} \Big\},$$

where $\mathbb{R}^{2n}$ denotes the trivial bundle, and $j$, $\nu$ are almost complex structures and matrix-valued functions as in §3.2.3, while $C_{\epsilon:(\delta,\delta)}$ denotes the exponentially weighted space:

$$C_{\epsilon:(\delta,\delta)} := \{\xi| \varsigma^{(\delta,\delta)}\xi \in C_\epsilon\};$$

$\delta > 0$ is as in §4.3.1. The limiting operators $E_l, E_r : L_1^p(S_1^1; \mathbb{R}^{2n}) \to L^p(S_1^1; \mathbb{R}^{2n})$ above are:

$$E_l\eta(t) := j(-\infty, t)d\eta/dt + \nu(-\infty, t)\eta(t); \quad E_r\eta(t) := j(\infty, t)d\eta/dt + \nu(\infty, t)\eta(t).$$

We know that if $u \in \mathcal{M}_P(x, [w]; y, [v])$, then $E_u \in \Sigma_P(x, [w]; y, [v])$ by the exponential decay of $u$.

Given $(\underline{x}, [\underline{w}]), (\underline{y}, [\underline{v}]) \in \tilde{\mathcal{P}}$, let

$$\Sigma_P(\underline{x}, [\underline{w}]; \underline{y}, [\underline{v}]) := \bigcup_{(x,[w])\in(\underline{x},[\underline{w}]);(y,[v])\in(\underline{y},[\underline{v}])} \Sigma_P(x, [w]; y, [v])$$

endowed with the product topology. This is a contractible space because of the fibration:

$$\Sigma_P(x, [w]; y, [v]) \to \Sigma_P(\underline{x}, [\underline{w}]; \underline{y}, [\underline{v}]) \to \tilde{O}_x \times \tilde{O}_y = \mathbb{R} \times \mathbb{R}.$$

The space $\Sigma_P(x, [w]; y, [v])$ is contractible as usual (see e.g. the proof of [FHS] Proposition 7); thus so is $\Sigma_P(\underline{x}, [\underline{w}]; \underline{y}, [\underline{v}])$. The determinant line bundle over $\Sigma_P(\underline{x}, [\underline{w}]; \underline{y}, [\underline{v}])$ is thus trivial, therefore orientable.

Note that for different liftings, the determinant line bundles

$$\mathcal{L}_P(\underline{x}, [\underline{w}]; \underline{y}, [\underline{v}]) \to \Sigma_P(\underline{x}, [\underline{w}]; \underline{y}, [\underline{v}]), \quad \mathcal{L}_P(\underline{x}, [\underline{w}']; \underline{y}, [\underline{v}']) \to \Sigma_P(\underline{x}, [\underline{w}']; \underline{y}, [\underline{v}'])$$

are isomorphic as long as

$$\mathrm{ind}(\underline{x}, [\underline{w}]) - \mathrm{ind}(\underline{y}, [\underline{v}]) = \mathrm{ind}(\underline{x}, [\underline{w}']) - \mathrm{ind}(\underline{y}, [\underline{v}']) = k.$$

We shall therefore identify their respective spaces of possible orientations and denote it as $\mathrm{or}(x, y; k)$.

**4.5.2 The Splicing construction.** The linear gluing theorem in the next lemma is an example of the general principle of "replacing the ordinary gluing theorem by the obvious adaptation to the fibered version" for this $S^1$-equivariant version of Floer theory. Since we have not done so before, we shall provide more details here, to demonstrate this type of modification.

First, we introduce some preliminary definitions. Similarly to the case of $E_u$ in §4.3.1, given $E \in \Sigma_P(x, [w]; y, [v])$, we may define its stabilizations $E^{+}, E^{--}, E^{-+}$ by extending its domain of definition respectively to

$$L_{1:(\delta,\delta)}^{\cdot+,p}(\Theta,\mathbb{R}^{2n};x,[w];y,[v]) := \left\{ C\beta(s)\Phi_{y,[v]}\dot{y} + \xi \,\Big|\, C \in \mathbb{R}, \xi \in L_{1:(\delta,\delta)}^{p}(\Theta,\mathbb{R}^{2n}) \right\},$$

$$L_{1:(\delta,\delta)}^{-,p}(\Theta;\mathbb{R}^{2n};x,[w];y,[v]) := \left\{ C\beta(-s)\Phi_{x,[w]}\dot{x} + \xi \,\Big|\, C \in \mathbb{R}, \xi \in L_{1:(\delta,\delta)}^{p}(\Theta,\mathbb{R}^{2n}) \right\},$$

$$L_{1:(\delta,\delta)}^{-+,p}(\Theta;\mathbb{R}^{2n};x,[w];y,[v]) := \Big\{ (C_-\beta(-s)\Phi_{x,[w]}\dot{x} + C_+\beta(s)\Phi_{y,[v]}\dot{y} + \xi \,\Big|$$
$$C_\pm \in \mathbb{R}, \xi \in L_{1:(\delta,\delta)}^{p}(\Theta,\mathbb{R}^{2n}) \Big\}.$$

Note that by taking left and right limits; these spaces naturally fiber over $\mathbb{R}, \mathbb{R}, \mathbb{R} \times \mathbb{R}$ respectively:

$$L_{1:(\delta,\delta)}^{\cdot+,p}(\Theta,\mathbb{R}^{2n};x,[w];y,[v]) \xrightarrow{\pi_+} \mathbb{R};$$

$$L_{1:(\delta,\delta)}^{-,\cdot,p}(\Theta;\mathbb{R}^{2n};x,[w];y,[v]) \xrightarrow{\pi_-} \mathbb{R};$$

$$L_{1:(\delta,\delta)}^{-+,p}(\Theta;\mathbb{R}^{2n};x,[w];y,[v]) \xrightarrow{\pi_-\times\pi_+} \mathbb{R} \times \mathbb{R}.$$

Next, recall from [FH] that an $E \in \Sigma_P(x,[w];y,[v])$ is *asymptotically constant* if $j,\nu$ are $s$-independent for $s$ large. Let $E_1 = \partial_s + j_1\partial_t + \nu_1 \in \Sigma_P(x,[w];y,[v])$, $E_2 = \partial_s + j_2\partial_t + \nu_2 \in \Sigma_P(y,[v];z,[q])$, then for all sufficiently large $R$, the glued operator

$$E_1 \#_R E_2 := \begin{cases} \partial_s + j_1(s+R,t)\partial_t + \nu_1(s+R,t) & \text{when } s \in [-\infty,0] \times S_1^1; \\ \partial_s + j_2(s-R,t)\partial_t + \nu_2(s-R,t) & \text{when } s \in [0,\infty] \times S_1^1. \end{cases} \quad (44)$$

Also, given $\xi_1 \in L_{:(\delta,\delta)}^{-+,p}(\Theta;\mathbb{R}^{2n};x,[w];y,[v]), \xi_2 \in L_{:(\delta,\delta)}^{-+,p}(\Theta;\mathbb{R}^{2n};y,[v];z,[q])$ with $\pi_+(\xi_1) = \pi_-(\xi_2)$, let

$$\#_R(\xi_1,\xi_2) = \beta(s)\xi_2(s-R,t) + (1-\beta(s))\xi_1(s+R,t),$$

where $\beta(s)$ is a smooth cutoff function that is zero for $s < 0$, and 1 for $s > 1$.

In [L] §5.1.1, we introduced the convenient language of "K-models" : a K-model for a Fredholm operator $D$ is a pair of finite dimensional oriented vector spaces $[\ker_g D; \operatorname{coker}_g D]$ ("generalized kernel" and "generalized cokernel"), such that the domain and range of $D$ split (not necessarily orthogonally) as $B \oplus \ker_g D$, $\operatorname{coker}_g D \oplus D(B)$ respectively, and $D$ is an isomorphism between $B$ and $D(B)$. A K-model for $D$ orients $\det D$ by a canonical isomorphism

$$\det D \simeq \det \ker_g D \otimes \det(\operatorname{coker}_g D)^*.$$

Different K-models which induce the same orientation on $\det D$ are said to be *co-oriented*.

A K-model $[\ker_g E; \operatorname{coker}_g E]$ for $E \in \Sigma_P(x,[w];y,[v])$ induces K-models for the various extensions of $E$ according to the recipe of [L] (5.2). For instance, the induced

K-models for $E^{-\cdot}, E^{-+}, E^{\cdot+}$ are respectively

$$\left[\ker_g E^{-\cdot} = \ker_g E \oplus \mathbb{R}\left(\beta(-s)\Phi_{x,[w]}\dot{x}\right) ; \mathrm{coker}_g E^{-\cdot} = \mathrm{coker}_g E\right].$$

$$\left[\ker_g E^{-+} = \ker_g E \oplus \mathbb{R}\left(\beta(-s)\Phi_{x,[w]}\dot{x}\right) \oplus \mathbb{R}\left(\beta(s)\Phi_{y,[v]}\dot{y}\right) ; \mathrm{coker}_g E^{-+} = \mathrm{coker}_g E\right].$$

$$\left[\ker_g E^{\cdot+} = \ker_g E \oplus \mathbb{R}\left(\beta(s)\Phi_{y,[v]}\dot{y}\right) ; \mathrm{coker}_g E^{\cdot+} = \mathrm{coker}_g E\right].$$

(Note particularly the ordering of the summands above). In terms of this language:

**Lemma.** *Given $E_1$, $E_2$ as in (44), and K-models $[\ker_g E_1, \mathrm{coker}_g E_1]$, $[\ker_g E_2, \mathrm{coker}_g E_2]$ for them, then for all sufficiently large $R$, the following forms a K-model for $(E_1 \#_R E_2)^{\cdot+}$:*

$$\left[\begin{array}{l}\ker_g (E_1 \#_R E_2)^{\cdot+} = (-1)^{\dim(\mathrm{coker}_g E_2) \operatorname{ind} E_1 + \operatorname{ind} E_2} \#_R(\ker_g E_1^{\cdot+} \times_\mathbb{R} \ker_g E_2^{-+}); \\ \mathrm{coker}_g(E_1 \#_R E_2)^{\cdot+} = \#_R(\mathrm{coker}_g E_2 \times \mathrm{coker}_g E_1)\end{array}\right]; \quad (45)$$

*furthermore, co-oriented choices of the K-models $[\ker_g E_1; \mathrm{coker}_g E_1]$, $[\ker_g E_2; \mathrm{coker}_g E_2]$ of $E_1$ and $E_2$ in (45) yield co-oriented K-models.*

Notice that according to the above formula,

$$\ker_g(E_1 \#_R E_2)^{\cdot+} \simeq (-1)^{\dim(\mathrm{coker}_g E_2) \operatorname{ind} E_1^{\cdot+}} \ker_g E_1^{\cdot+} \times \ker_g E_2^{\cdot+}.$$

We take (45) for the definition of the orientation of $\det(E_1 \#_R E_2)^{\cdot+}$, and its (de)stabilizations.

*Proof.* Equivalently, we shall prove that:

Set $W' := \{\xi | \xi \in L^p_{1:(\delta,\delta)}(\Theta, \mathbb{R}^{2n}; x, [w]; z, [r]), \langle \dot{y}, \xi(0, \cdot)\rangle_{2,t} = 0\}$, then the restriction $(E_1 \#_R E_2)' := (E_1 \#_R E_2)\big|_{W'}$ has as a K-model:

$$\left[(-1)^{\dim(\mathrm{coker}_g E_2) \operatorname{ind} E_1} \#_R(\ker_g E_1 \times \ker_g E_2); \#_R(\mathrm{coker}_g E_2 \times \mathrm{coker}_g E_1)\right],$$

and co-oriented choices of $[\ker_g E_1; \mathrm{coker}_g E_1]$, $[\ker_g E_2; \mathrm{coker}_g E_2]$ above yield co-oriented K-models.

We leave the reader the task of verifying the independence of the orientation of the choices involved, and only indicate how the routine argument should be modified to prove the gluing theorem. An example of such routine argument may be found in e.g. [FH]: First one may stabilize to reduce to the case $\mathrm{coker}_g = \mathrm{coker}$ are trivial. Then one may construct a right inverse either by partition of unity, or by running the proof by contradiction argument in e.g. [F88b, L]. In the first method, let $P_1, P_2$ be respectively the right inverses of the translates $E_1^R(s,t) := \partial_s + j_1(s+R,t)\partial_t + \nu_1(s+R,t)$, $E_2^{-R}(s,t) := \partial_s + j_2(s-R,t)\partial_t + \nu_2(s-R,t)$, and set

$$Q := \varphi_1 P_1 (1-\beta) + \varphi_2 P_2 \beta,$$

where $\varphi_1$, $\varphi_2$ are smooth cutoff functions with small $\|d\varphi_i\|_{C_\epsilon}$ and whose supports include those of $1 - \beta, \beta$ respectively. Then $DQ = 1 + \mathcal{R}$, where $\mathcal{R}$ is small in operator norm because of the boundedness of $P_1, P_2$, and we can take $P = Q(1+\mathcal{R})$. Because of the 1-dimensional constraint $\langle \dot{y}, \xi(0, \cdot)\rangle_{2,t} = 0$ in the definition of $W'$ and the index computation, $\ker_g L_1 \oplus \ker_g L_2$ must be complementary to the image of $P$ in $W'$.

For the second method, note that the condition $\langle \dot{y}, \xi(0, \cdot)\rangle_{2,t} = 0$ implies that $\xi$ is small in the middle region $\Theta_0$ (in the notation of [L]), and thus the usual estimates (see e.g. [L] §7.2) go through. $\square$

The rest of the arguments then follow [FH]: The proceeding Lemma defines a splicing homomorphism $\# : \mathrm{or}(x, y; k) \times \mathrm{or}(y, z; k') \to \mathrm{or}(x, z; k + 1)$, and we might choose a coherent system of orientations, namely an element $\mathfrak{o}_{x,y;k} \in \mathrm{or}(x, y; k)$ for every pair $x, y \in \mathcal{P}$ and every $k \in \mathbb{Z}$, such that $\#(\mathfrak{o}_{x,y;k}, \mathfrak{o}_{y,z;k'}) = \mathfrak{o}_{x,z;k+k'}$ for every triple $x, y, z \in \mathcal{P}$, and every pair of integers $k, k'$.

For each $u \in \mathcal{M}_P(x, [w]; y, [v])$, we may choose a trivialization $\Phi_u$ of $u^*TM$ extending $\Phi_{x,[w]}$, $\Phi_{y,[v]}$, so that $\Phi_u E_u \Phi_u^{-1} \in \Sigma_P(x, [w]; y, [v])$ and hence $\det E_u$ is assigned an orientation from the orientation of the determinant line bundle over $\Sigma_P(x, [w]; y, [v])$. This assignment does not depend on the choice of $\Phi_u$ by the arguments of [FH] Lemmas 13, 15.

# 5 Gromov vs. $I_F$: the lagrangian intersection case

In this section, we prove Theorems A, B, and Corollary C, *assuming* the invariance theorem 2.2.2. The key is the relation of $I_F(L, L')$ with a Gromov-type invariant counting pseudo-holomorphic annuli. (Lemma 5.1.3 below).

Throughout this section, let $(Y, \theta, b)$ be a Type F triple and let $L, L' \subset T^*Y$ be the zero section and the section corresponding to $\theta$.

## 5.1 A Gromov-type invariant for annuli

Though not explicitly spelled-out, the key of paper [GL] is the computation of a Gromov-type invariant, which we define below.

Let $M = T^*Y$, $\Theta' = [0, 1] \times S^1_1$, and let $\mathrm{pr}_2 : \Theta' \times M \to M$ be the projection. Note that under our assumptions, $H_2(M, L \cup L')$ is isomorphic to the anti-diagonal of $H_1(L) \oplus H_1(L')$ via the relative homology sequence of the pair $(M, L \cup L')$, and the anti-diagonal is in turn isomorphic to $H_1(L)$. We denote the isomorphism from $H_2(M, L \cup L')$ to $H_1(L)$ obtained in this way by $i_{rel}$.

Given $\nu \in C_\epsilon(\mathrm{pr}_2^* TM)$ and an admissible almost complex structure $J$, let

$$\mathcal{N}(L, L', b; J, \nu) := \Big\{(\rho, v) \,\Big|\, v : [0, 1] \times S^1_1 \to M, \, v(0, \cdot) \in L, v(1, \cdot) \in L',$$

$$\rho \in (0, \infty), \frac{1}{\rho}\partial_s v + J(v)\partial_t v + \nu(s, t, v) = 0, \, (s, t) \in \Theta', \, i_{rel}([v]) = b\Big\}.$$

By rescaling, $\mathcal{N}(L, L', b; J, \nu)$ may be interpreted as the space of perturbed $J$-holomorphic annuli of length 1 and arbitrary width $\rho \in \mathbb{R}^+$. Let

$$\Upsilon(L, L'; J) := \left\{\nu \,\Big|\, \nexists\, (s, x) \text{ s.t. } J(x(t)) + \nu(s, t, x(t)) = 0\, \forall t,\, s \in [0, 1],\, x : S_1^1 \to M \right\}$$
$$\subset C_\epsilon(\mathrm{pr}_2^* TM).$$

By standard arguments (see e.g. [GL]),

**5.1.1 Lemma.** *Let $M, L, L', b, J$ be as above. Then:*
    *(a) There is a Baire set $\mathcal{R}(L, L'; J) \subset C_\epsilon(\mathrm{pr}_2^* TM)$, such that for any $\nu \in \mathcal{R}(L, L'; J)$, $\mathcal{N}(L, L', b; J, \nu)$ is a smooth 1-manifold.*
    *(b) $\mathcal{N}(L, L', b; J, \nu)$ is compact for $\nu \in \Upsilon(L, L'; J)$.*
    *(c) If $\nu_t$, $t \in [0, 1]$ is a path in $\Upsilon(L, L'; J)$, and $\nu_0, \nu_1 \in \mathcal{R}(L, L'; J)$, then $\mathcal{N}(L, L', b; J, \nu_0)$, $\mathcal{N}(L, L', b; J, \nu_1)$ are cobordant.*

By (a) and (b) of the above Lemma, we may define:

**5.1.2 Definition. (Gromov invariant for annuli)** For $\nu \in \Upsilon(L, L'; J) \cap \mathcal{R}(L, L'; J)$, let

$$\mathrm{Gr}(L, L', b; J, \nu) := \mathrm{ev}^*\, \theta[\mathcal{N}(L, L', b; J, \nu)]/\theta(b), \tag{46}$$

where $\mathrm{ev} : \mathcal{N}(L, L', b; J, \nu) \to L$ is the evaluation map:

$$\mathrm{ev}(v) = v(0, 0).$$

Apparently, this is an analogue of (16), and similarly it does not depend on the choice of $\theta$.

Next we compare Gr with $I_F$ when $\nu = 0$. Note that since $M = T^*Y$, the map $e_{0*}$ in (19) is actually an isomorphism.

**5.1.3 Lemma. (Comparison with $I_F(L, L')$)** *If $0 \in \Upsilon(L, L'; J) \cap \mathcal{R}(L, L'; J)$, then*

$$\begin{aligned}\mathrm{Gr}(L, L', b; J, 0) &= \chi(\hat{\mathcal{M}}_O(L, L', (e_{0*})^{-1}b; J, 0)) \\ &= \left(\ln I_F(L, L'; J, 0)\right)(e_{0*})^{-1}b).\end{aligned} \tag{47}$$

*Proof.* We first observe that there is an $S^1$-equivariant isomorphism

$$\mathcal{N}(L, L', b; J, 0) \simeq \mathcal{M}_O(L, L', (e_{0*})^{-1}b; J, 0)$$

by mapping $(\rho, v) \in \mathcal{N}(L, L', b; J, 0)$ to $(\rho = 1/T, u) \in \mathcal{M}_O(L, L', (e_{0*})^{-1}b; J, 0)$, where

$$u(s, t) = v(t, -\rho s).$$

The $S^1$-action on both moduli spaces is free by the primitivity assumption on $b$. Thus both moduli spaces consist of a finite number of circles mapping to circles of homology class $b$ under the evaluation map ev. So the right hand side of (46) counts the number such circles, i.e. $\chi(\hat{\mathcal{M}}_O(L, L', (e_{0*})^{-1}b; J, 0))$. Hence the first equality in (47).

Lastly, since $\theta$ is nowhere vanishing and hence $L \cap L' = \emptyset$, $\mathcal{P}(L, L'; J, 0) = \emptyset$ and

$$I_F(L, L'; J, 0) = \zeta_F(L, L'; J, 0) \in \text{Nov}^1(\ker \psi_\mu, -[\mathcal{Y}_X]; \mathbb{Q}),$$

and the second equality in (47) follows. □

## 5.2 Transversality by perturbing Lagrangians

We shall need the following alternative version of transversality result, where $J, X$ is fixed, and the lagrangian submanifold $L$ is perturbed.

Let

$$\mathcal{Z}_\epsilon^\varepsilon(L) := \{\kappa \,|\, \kappa \in \Omega^1(L),\, d\kappa = 0,\, \|\kappa\|_{C_\epsilon} < \varepsilon\}.$$

A $\kappa$ represents a lagrangian submanifolds $C^\infty$-close to the zero section $L \subset T^*Y = M$. Let $\phi_\kappa(L)$ denote this corresponding lagrangian submanifold.

**5.2.1 Proposition.** *Let $Y, \theta, b, L, L'$ be as in the statement of Theorem A, and let $J$ be the standard almost complex structure on $M = T^*Y$. Then for all sufficiently small $\varepsilon > 0$, there exists a Baire set $\mathcal{Z}_{reg}^\varepsilon \subset \mathcal{Z}_\epsilon^\varepsilon(L)$, such that $\mathcal{M}_O(\phi_\kappa(L), L', b; J, 0)$ is a compact 1-manifold for all $\kappa \in \mathcal{Z}_{reg}^\varepsilon$.*

The Proposition follows from the combination of the following two lemmas.

We call a curve $(1/T, u) \in \mathcal{M}_O(L, L'; J, 0)$ *left-simple* if there is a $s \in S_T^1$, such that $\partial_s u(s, 0) \neq 0$, and $u(s, 0) \neq u(s', 0)$ for any other $s' \in S_T^1$. This is an analogue of the notion of "somewhere injectivity" in the case of closed curves (cf. [MS]). Left-simple curves should be distinguished from simple curves (i.e. not multiple covers): A left-simple curve is simple, but the converse might not be true (cf. [Oh97b]).

Let $\mathcal{M}_O^{ls}(L, L'; J, 0) \subset \mathcal{M}_O(L, L'; J, 0)$ denote the subset of left-simple curves. The following is an analogue of the standard fact that somewhere-injectivity implies transversality.

**5.2.2 Lemma.** *For a fixed $t$-independent compatible almost complex structure $J$, and a $\phi$ in a Baire set of $C_\epsilon(L)$, $\mathcal{M}_O^{ls}(\phi(L), L', B; J, 0)$ is a smooth 1-manifold.*

This lemma follows from an obvious adaptation of [Oh96a]. E.g. $e^{J\theta}$ in p.516 there should be replaced by $J(u(s, 0))$.

**5.2.3 Lemma.** *For any primitive class $b \in H_1(L; \mathbb{Z})$, any $(1/T, u) \in \mathcal{M}_O(L, L', b; J, 0)$ is left-simple.*

*Proof.* Let $\tau : T^*Y \to [0, \infty)$ denote the distance square to the zero section, i.e. $\tau = |p|^2$ in the local coordinates $(q, p)$. $u^*\tau$ is a sub-harmonic function with minimal value 0 occurring at $(s, 0) \in S_T^1 \times [0, 1]$ for all $s$, and with maximum at $(s', 1)$ for some $s' \in S_T^1$. Let $\Theta_c \subset \Theta = S_T^1 \times [0, 1]$ be the path component of $u^{-1}(\tau^{-1}[0, c])$ containing $S^1 \times \{0\}$, and let $\Theta_c^0, \Theta^0$ denote interiors of $\Theta_c, \Theta$ respectively. Due to the standard fact that $u$ has isolated critical points or self-intersections in $\Theta^0$, there is a small enough $\delta$, such that $u_\delta := u\big|_{\Theta_\delta}$ has no critical points or self-intersection points except in $S_T^1 \times \{0\}$. Note that by construction $u^{-1}(L) \cap \Theta_\delta^0 = \emptyset$, $u^{-1}(L') \cap \Theta_\delta^0 = \emptyset$, and $u^{-1}(u(\partial\Theta_\delta)) = \partial\Theta_\delta$. On the other hand, $u_\delta$ factors over the Riemann surface defined by $\Sigma_\delta := (u(\Theta_\delta), J)$ and is a local isomorphism. The degree of $\Theta_\delta \to \Sigma_\delta$ is locally constant in the interior and finite due to the absence of bubbling. (See also [Laz] §5.3.) Therefore $u_\delta$ is a finite covering map in the interior. By the primitivity of $b$, it can only be an isomorphism over $\Theta_\delta^0$.

Now, by [Oh97b] Theorem 2.1 or [MiW], $u_\delta$ only has finitely many critical points on the boundary. So $u(S_T^1 \times \{0\})$ is a 1-complex, with the 0-cells being critical values. Because $u_\delta$ is injective in the interior, by [Oh97b] Corollary 2.4 or [MiW], the preimage of each 1-cell in $u(S^1 \times \{0\})$ either consists of a 1-cell, or two 1-cells inducing opposite orientations in the image. If $(1/T, u)$ is not left-simple, this would imply that the homology class of $u\big|_{S^1 \times \{0\}}$ is trivial, contradicting the primitivity assumption on $b$. (In fact, in our case one may even argue that $u\big|_{S^1 \times \{0\}}$ is an isomorphism: if the second case above happens, neighborhoods of the two 1-cells match up to form a holomorphic disk in $T^*Y$ including this 1-cell in the interior. now $\tau = 0$ on this 1-cell, so $\tau \equiv 0$ on the entire disk. But $\tau$ can not be 0 in the interior of $\Theta$ unless $u$ is entirely in $L$, which is impossible.) □

To continue the proof of the Proposition, note from the proof of last lemma that when $\varepsilon > 0$ is sufficiently small, for any $\kappa \in \mathcal{Z}_\epsilon^\varepsilon(L)$ and any $u \in \mathcal{M}_O(\phi_\kappa(L), L', b; J, 0)$, one may still find suitable $\Theta_\delta^0$ with the same properties, so that the rest of the arguments go through to show that:

All elements in $\mathcal{M}_O(\phi_\kappa(L), L', b; J, 0)$ are left-simple.

This, together with Lemma 5.2.2, proves Proposition 5.2.1. □

## 5.3 Proofs of Theorems A, B, and Corollary C

The proof of Theorem A basically follows the framework of [GL], replacing the count of holomorphic annuli there by the computation of $I_F(L, L')$ in Theorem 2.2.3. Given that, Theorem B and Corollary C follow directly from the arguments of [GL].

### 5.3.1 Proof of Theorem A.

Let $M = T^*Y$, $L$, $L'$, $H_t$, $m$, $\mu$ be as in the statement of Theorem A.

We shall show that Theorem A holds if $L$ is replaced by $\phi_\kappa L$ for any $\kappa \in \mathcal{Z}^\varepsilon_{reg}$. Apparently this would imply Theorem A. This in turn follows from the three observations below:

**(a)** $0 \in \mathcal{R}(\phi_\kappa L, L', J)$ since $\kappa \in \mathcal{Z}^\varepsilon_{reg}$. Furthermore, by Lemma 5.1.3 and Theorem 2.2.3,
$$\mathrm{Gr}(\phi_\kappa L, L', b; J, 0) = \ln i_{-[\theta]}\tau(Y)(b) \neq 0$$
by the assumption that $b$ is $\theta$-essential. (More precisely, $(J, 0)$ might not be regular; if so $I_F(\phi_\kappa L, L'; J, 0)$ is not well-defined. However, by Proposition 3.2.2, we can find a very small $H$ so that $(J, \chi_H)$ is regular while leaving $\chi(\hat{\mathcal{M}}_O(\phi_\kappa L, L'; J, 0)) = \chi(\hat{\mathcal{M}}_O(\phi_\kappa L, L'; J, \chi_H)))$.

**(b)** $m/\mu \nabla H_t \in \mathcal{R}(\phi_\kappa L, L', J)$ and
$$\mathrm{Gr}(\phi_\kappa L, L', b; J, m/\mu \nabla H_t) = 0$$
due to the following lemma:

**Lemma.** *When $\lambda \geq m/\mu$, $\mathcal{N}(L, L', b; J, \lambda \nabla H_t) = \emptyset$.*

*Proof.* Say $(\rho, v) \in \mathcal{N}(L, L', b; J, \lambda \nabla H_t)$. Then
$$0 < \rho^{-1} \int_{\Theta'} |\partial_s v|^2 \, ds \, dt = \omega([u]) - \lambda \int (H_t(u(1,t)) - H_t(u(0,t))) \, dt \leq m - \lambda \mu.$$
Thus when $\lambda \geq m/\mu$, there can not be such $(\rho, v)$. □

**(c)** Let $\nu_\lambda$, be the line in $C_\epsilon(\mathrm{pr}_2^* TM)$ given by
$$\nu_\lambda = \lambda \nabla H_t, \ \lambda \in [0, m/\mu].$$
Suppose to the contrary of Theorem A that $\mathrm{PO}(\lambda \nabla H, b) = \emptyset \ \forall \lambda \in [0, m/\mu]$. Then the entire path $\nu_\lambda$ is in $\Upsilon(\phi_\kappa L, L', J)$, and according to Lemma 5.1.1 (c) and observations (a), (b) above,
$$\mathrm{Gr}(\phi_\kappa L, L', b; J, \nu_0) = \mathrm{Gr}(\phi_\kappa L, L', b; J, \nu_1),$$
contradicting the computations in (a), (b) above. This proves Theorem A. □

### 5.3.2 Theorem B and Corollary C.

Since both Theorem B and Corollary C follow from Theorem A by the same arguments as those in [GL], here we shall only outline the main ideas of the proofs.

To prove Theorem B, one uses the monotonicity lemma to argue that for small $\delta$, the perturbed pseudo-holomorphic annuli must be contained in the Weinstein neighborhood. This thus reduces the problem to Theorem A.

To prove Corollary C, suppose to the contrary that there is a hyperplane $f^{-1}(0) \subset \mathbb{C}^n$ separating $L_0$ and $\phi(L_1)$, where $\phi$ is a hamiltonian diffeomorphism, and $f$ is a linear function on $\mathbb{C}^n$. Then $f$ is a Hamiltonian function separating $L_0$, $\phi(L_1)$, and obviously $X_{df}$ has no periodic orbits. This contradicts a generalization of Theorem B to the case of Hamiltonians that are linear at infinity.

# 6 Gromov vs. $I_F$: the $S^1$-equivariant case

*Assuming the invariance theorem 2.3.1, we prove Theorem 2.3.3 in §6.1 by comparing $I_F^{S^1}$ with the Gromov invariant for tori. In §6.2 we use this result to establish Theorem D.*

## 6.1 Proof of Theorem 2.3.3

According to Theorem 2.3.1 and Proposition 4.4.2, we may compute $I_F^{S^1}(M, \gamma_0)$ at a regular $(J, X)$ satisfying

$$X \in \mathcal{X}_{ad}(0) \text{ and } \|X\|_{C_\epsilon} \ll 1. \tag{48}$$

Clearly $\mathcal{P}(M, \gamma_0; X) = \emptyset$ in this case, and

$$(\ln i_{[\omega]} I_F^{S^1}(M, \gamma_0))(A) = \chi(\hat{\mathcal{M}}_O(M, \gamma_0, A; J, X)); \tag{49}$$

because of the primitivity assumption on $A$, the right hand side is just a signed count of elements in $\hat{\mathcal{M}}_O(M, \gamma_0, A; J, X))$. To compare this with the Gromov invariant, we need to first go into some details of the definition of the Gromov invariant. The version of Gromov invariant more suitable for our purpose is that of Ruan-Tian [RT2], though in our case it agrees with the version of Taubes [T].

Let $\mathcal{M}_{1,1}$ denote the moduli space of genus 1 Riemann surfaces with 1 marked point, namely,

$$\{(j, p)\}/\operatorname{Diff}(\mathbb{T}^2),$$

where $j$ is a complex structure on the torus $\mathbb{T}^2$, and $p \in \mathbb{T}^2$. It is well-known that $\mathcal{M}_{1,1} = \mathbb{H}/\operatorname{SL}_2(\mathbb{Z})$, $\mathbb{H} \subset \mathbb{C}$ being the upper half plane (i.e. the genus 1 Teichmüller space).

Let [3]

$$\mathcal{M}_{1,1,A}(M; J, X) = \{(j, p, u)\}/\operatorname{Diff}(\mathbb{T}^2),$$

where $(j, p)$ is a pair of a complex structure and a point on $\mathbb{T}^2$ as before, and $u: \mathbb{T}^2 \to M$ satisfies the conditions that

$du \circ j - J \circ du = i_\mathbb{R} X(u)$, $i_\mathbb{R}: T_\mathbb{R} M \to T^{0,1} M$ being the canonical isomorphism, and $[u] = A \in H_2(M; \mathbb{Z})$.

(50)

---

[3]In [RT2], the perturbation $X$ is allowed to depend on $s, t, j$.

Note that there is a $\mathbb{T}^2$-action on $\mathcal{M}_{1,1,A}(M;J,X)$ by translation on $\mathbb{T}^2$, and there is a "forgetful map"

$$\mathrm{ft}: \mathcal{M}_{1,1,A}(M;J,X) \to \mathcal{M}_{1,1}, \quad \mathrm{ft}(j,p,u) = (j,p).$$

Now, the moduli space of $J|X$-tori $\mathcal{M}_O(M,\gamma_0,A;J,X)$ also have a similar description: Let's begin with some definitions.

A $j^\star$-*structure* on $\mathbb{T}^2$ is a pair $(j,\gamma)$, where $j$ is a complex structure of the torus, and $\gamma$ is a primitive homology class in $H_1(\mathbb{T}^2;\mathbb{Z})$. A $j^\star$-*holomorphic map* from $\mathbb{T}^2$ to $\mathbb{T}^2$ is a smooth map preserving the $j^\star$ structure. It is easy to see that the group of $j^\star$-automorphisms of $\mathbb{T}^2$ is $\mathbb{T}^2$ itself (translations), by reducing the problem to lattice automorphisms.

Two $j^\star$-*structures* are equivalent if they are related by a diffeomorphism of $\mathbb{T}^2$. The moduli space of all equivalence classes of $j^\star$-structures, is

$$\mathbf{T}_{1,1}^\gamma := \{(j,\gamma,p)\}/\operatorname{Diff}(\mathbb{T}^2) = \{(j,p)\}/\operatorname{Diff}_\gamma(\mathbb{T}^2),$$

where $(j,\gamma)$ is a $j^\star$-structure on $\mathbb{T}^2$, and $p \in \mathbb{T}^2$ is a point, and $\operatorname{Diff}_\gamma(\mathbb{T}^2)$ above is the group of diffeomorphisms of $\mathbb{T}^2$ preserving $\gamma$. It is easy to see that

$$\mathbf{T}_{1,1}^\gamma = \mathbb{R}_+ \times S^1 = \mathbb{H}/\mathbb{Z},$$

where $\mathbb{Z} \subset \operatorname{SL}_2(\mathbb{Z})$ is the subgroup generated by $\begin{pmatrix} 1 & 1 \\ 0 & 1 \end{pmatrix}$.

Each $j^\star$-structure on $\mathbb{T}^2$ has a standard representative as $\Theta_{T,Q}$, with the distinguished element $\gamma \in H_1(\Theta_{T,Q})$ represented by the unit circle $\{x|x \in \mathbb{R} \subset \mathbb{C}\}/\Gamma_{T,Q}$; so the parameter space of $T,Q$ in section 4 is precisely this moduli space of $j^\star$-structures.

On the other hand, for an $A \in \mathfrak{H} \subset H^2(M;\mathbb{Z})$, $\mathcal{M}_O(M,\gamma_0,A;J,X)$ may be alternatively described as follows.

$$\mathcal{M}_O(M,\gamma_0,A;J,X) = \{(j,\gamma,p,u)\}/\operatorname{Diff}(\mathbb{T}^2),$$

where $(j,\gamma)$ is a $j^\star$ structure, $p \in \mathbb{T}^2$, and $u: \mathbb{T}^2 \to M$ satisfies, in addition to (50), that $u_*\gamma = \gamma_0$.

Similar to $\mathcal{M}_{1,1,A}(M;J,X)$, there is also a $\mathbb{T}^2$ action on $\mathcal{M}_O(M,\gamma_0,A;J,X)$ by translation on $\mathbb{T}^2$, and a "forgetful map"

$$\mathrm{ft}_\gamma : \mathcal{M}_O(M,\gamma_0,A;J,X) \to \mathbf{T}_{1,1}^\gamma, \quad \mathrm{ft}_\gamma(j,\gamma,p,u) = (j,\gamma,p).$$

It follows from the above descriptions that we have the following commutative diagram:

$$\begin{array}{ccc} \mathcal{M}_O(M,\gamma_0,A;J,X) & \xrightarrow{\mathrm{cp}} & \mathcal{M}_{1,1,A}(M;J,X) \\ \downarrow{\mathrm{ft}_\gamma} & & \downarrow{\mathrm{ft}} \\ \mathbf{T}_{1,1}^\gamma & \longrightarrow & \mathcal{M}_{1,1}, \end{array}$$

where the horizontal maps are obtained by deleting $\gamma$ from the quadruple $(j,\gamma,p,u)$ or triple $(j,\gamma,p)$.

**6.1.1 Lemma.** *Suppose that $A \in \mathfrak{H} \subset H_2(M; \mathbb{Z})$ is a primitive class, $X$ is as in (48), and that $(J, X)$ is regular. Then the map cp above is a $\mathbb{T}^2$-equivariant diffeomorphism, and both $\mathcal{M}_O(M, \gamma_0, A; J, X)$ and $\mathcal{M}_{1,1,A}(M; J, X)$ are both compact smooth (nondegenerate) manifolds.*

*Proof.* By Proposition 4.4.2, the regularity of $(J, X)$ guarantees the smoothness and compactness of $\mathcal{M}_O(M, \gamma_0, A; J, X)$; thus we only need to demonstrate that cp is a diffeomorphism. In fact, we shall only show the injectivity of cp, since it is straightforward to see that cp is smooth and surjective.

In other words, we shall show that for any $(j, p, u) \in \mathcal{M}_{1,1,A}(M; J, X)$, $J, \gamma_0 \in H_1(M; \mathbb{Z})$ and $u$ together determine a unique $j^*$-structure $(j, \gamma)$ on the domain of $u$. By primitivity of $A$, $\mathcal{M}_{1,1,A}(M; J, X)/\mathbb{T}^2$ consists of embedded tori in $M$. Thus cp is injective if $u_* : H_1(\mathbb{T}^2) \to H_1(M)$ is injective. Since $u_*\gamma = \gamma_0$, and by assumption $[\gamma_0] \in u_*(H_1(\mathbb{T}^2))$ is primitive, the latter is true if $\underline{e}([u]) \in H_1(M)/\mathbb{Z}[\gamma_0]$ is nontorsion. Since the energy of $u$ $\mathcal{E}(u) = \underline{\mathrm{im}}^*\omega([u]) > 0$, this follows from Lemma 4.1.1 (b). □

Recall that the evaluation map $\mathrm{ev} : \mathcal{M}_{1,1,A}(M; J, X) \to M$ is defined by $\mathrm{ev}(j, p, u) = u(p)$, and the Gromov invariant $\mathrm{Gr}_{A,1,0}(M; J, X)$ is

$$\mathrm{Gr}_{A,1,0}(M; J, X) := [\mathrm{ev}(\mathcal{M}_{1,1,A}(M; J, X))] \cdot B/(A \cdot B) = \chi(\mathcal{M}_{1,1,A}(M; J, X)/\mathbb{T}^2),$$

where $B \in H_{2n-2}(M)$ is such that $A \cdot B \neq 0$. (In the statement of Theorem 2.3.3, $J, X$ is omitted from the notation $\mathrm{Gr}_{A,1,0}(M; J, X)$, since the Gromov invariant is independent of the choice of $J$ and small enough $X$).

Thus, Theorem 2.3.3 follows from the previous lemma and (49).

End of proof for Theorem 2.3.3.

## 6.2 Proof of Corollary D

Suppose $\mathrm{PO}(X; [\gamma_0]) = \emptyset$. Then $X \in \mathcal{X}_{nondeg}$ automatically. Furthermore, by Lemma 4.4.1 one may perturb $X$ slightly in a small neighborhood of $X^{-1}(0)$ into a $X' \in \mathcal{X}_{ad}$, so that $[\theta_{X'}] = [\theta_X]$, and

$$\mathrm{PO}(X'; [\gamma_0]) = \emptyset \tag{51}$$

still. According to Proposition 4.4.2, one may choose a $J$ such that $(J, X')$ is regular. Then $\tau_F(M, \gamma_0; J, X') = 1$ by (51).

Now by assumption $\mathrm{Gr}_{A,1,0}(M; \gamma_0) \neq 0$, and thus by Theorem 2.3.3 and Theorem 2.3.1,

$$I_F^{S^1}(M, \gamma_0; J, X') = \zeta_F(M, \gamma_0; J, X')\tau_F(M, \gamma_0; J, X') = \zeta_F(M, \gamma_0; J, X') \neq 1.$$

Namely, there is a $J|X'$-torus. Since $X'$ can be chosen arbitrarily close to $X$, this implies there is also a $J|X$-torus.

In the special case when $\mathfrak{h}[\theta_X] = [\omega]\big|_{\mathfrak{H}}$, since

$$\underline{\mathrm{im}}^*[\mathcal{Y}_X]\big|_{\mathfrak{H}} = [\omega]\big|_{\mathfrak{H}} - \mathfrak{h}[\theta_X] = 0 = \underline{\mathrm{im}}^*[\mathcal{Y}_{X'}]\big|_{\mathfrak{H}},$$

there is no $J|X$-tori nor $J|X'$-tori. This forces $\operatorname{PO}(X;[\gamma_0]) \neq \emptyset$.

End of proof of Corollary D.

## 6.3 Example: the symplectic mapping torus

We now work out the details of the example of symplectic mapping torus discussed in §2.3.2, 2.3.3, and 2.4.2.

Firstly, from the homotopy sequence of the fibration $\Sigma \to M \to \mathbb{T}^2$, we have
$$\pi_2(M(\Sigma, f)) = \pi_2(\Sigma).$$
On the other hand, observe that
$$c_1(TM(\Sigma, f))\Big|_{\pi_2(\Sigma)} = c_1(T\Sigma) \text{ and } \omega_{M(\Sigma,f)}\Big|_{\pi_2(\Sigma)} = \omega_\Sigma.$$
So $M(\Sigma, f)$ is monotone if $\Sigma$ is.

The homotopy assumption on $\gamma_0$ is clearly satisfied owing to the product structure of $M(\Sigma, f)$.

On the other hand, in this case the fibration $\mathcal{L}M \to \mathcal{C} \to M$ has a section; thus $1 \to \pi_2(M) \to \pi_1(\mathcal{C}) \to \pi_1(M) \to 1$ splits, and (32) now has the more precise form:
$$H_1(\underline{\mathcal{C}}; \mathbb{Z}) = H_1(M)/\mathbb{Z}[\gamma_0] \oplus \pi_2(M) = H_1(\Sigma_f) \oplus \pi_2(\Sigma).$$
With respect to the above expression, it is easy to see that
$$\underline{\operatorname{im}} = (\otimes[\gamma_0]) \oplus (\text{Hurewicz}),$$
the direct sum of the Hurewicz map and the map by tensoring with $[\gamma_0] \in H_1(S^1; \mathbb{Z})$. Now combining this with the fact that
$$c_1(TM) = \pi_1^* c_1(K_f), \tag{52}$$
one obtains the formula for $\mathfrak{H}$ (15). In (52), $K_f$ is the rank-$(2n-2)$ subbundle of $T\Sigma_f$, with fibers consisting of tangent vectors to $\Sigma$. $\pi_1 : \Sigma_f \times S^1 \to \Sigma_f$ denotes the projection.

Next, notice that $\operatorname{Image}(\operatorname{rest}^* \underline{e}^*)$ is precisely the first summand of
$$\operatorname{Hom}(\ker \underline{\psi}_c, \mathbb{R}) = H^1(\Sigma_f; \mathbb{R}) \oplus \operatorname{Hom}(\ker c_1\Big|_{\pi_2(\Sigma)}, \mathbb{R}),$$
and that in this case $\ker \psi_c \simeq \mathfrak{H}$. Thus, any class in $H^1(\Sigma_f; \mathbb{R}) \subset H^1(M; \mathbb{R})$ is $H^2$-induced, and $\mathfrak{h}$ is simply the injection of $H^1(\Sigma_f; \mathbb{R})$ to the first summand above.

Since $M$ is monotone, $[\omega]\Big|_{\ker c_1\big|_{\pi_2(M)}} = 0$. Recalling the definition of $\omega$ in (14), we see that $[\omega]\Big|_{\mathfrak{H}} = \mathfrak{h}[q]$.

With all the above elements in place, the statement in Example 2.4.2 is a straightforward consequence of Corollary D, and the computation of the relevant Gromov invariant in [IP2].

# 7 The invariance proofs

In this section we prove Theorems 2.2.2, 2.2.3, 2.3.1. These theorems are completely analogous to Theorems 1.5, 1.6, and Corollary 1.7 in [L], and not surprisingly, the proofs are straightforward adaptations of [L].

## 7.1 Proofs of Theorems 2.2.2, 2.2.3, 2.3.1

We shall first state a general invariance theorem for all versions of $I_F$, Proposition 7.1.1 below; then show how Theorems 2.2.2, 2.2.3, 2.3.1 follow from this Proposition.

**7.1.1 Proposition.** *If $(J_1, X_1), (J_2, X_2)$ are both regular, and if there is a path $X_\lambda, \lambda \in [1,2]$ in $\mathcal{X}$ connecting $X_1, X_2$, and a $[\mathcal{Y}] \in H^1(\mathcal{C}; \mathbb{R})$ such that $[\mathcal{Y}_{X_\lambda}]\big|_{\ker \psi} = \alpha_\lambda [\mathcal{Y}]\big|_{\ker \psi}$ for some $\alpha_\lambda \geq 0 \ \forall \lambda$, then*

$$I_F(J_1, X_1) = I_F(J_2, X_2).$$

Actually, the last expression in the Proposition needs some clarification. When one of $[\mathcal{Y}_{X_i}]\big|_{\ker \psi}$, $i = 1, 2$ is trivial and the other is not, the two sides of the equation do not really take values in the same monoid. Say $[\mathcal{Y}_{X_1}]\big|_{\ker \psi} = 0$ and $[\mathcal{Y}_{X_2}]\big|_{\ker \psi} \neq 0$. In this case we take the equality to mean

$$i_{-[\mathcal{Y}_{X_2}]} I_F(J_1, X_1) = I_F(J_2, X_2).$$

The proof of this Proposition for the two versions of Floer theories under discussion in this paper is postponed to §7.2.

*Proof of Theorem 2.2.2.* Part (a) is a direct consequence of the Proposition, since $[\mathcal{Y}_{X+\chi H_i}] = [\mathcal{Y}_X]$. On the other hand, part (b) follows from the next, more general theorem. □

We need some preliminaries to state the next theorem.

The monotonicity of $L'$ means that there is a constant $\alpha > 0$ such that $\omega\big|_{\pi_2(M,L')} = \alpha c_1\big|_{\pi_2(M,L')}$. This implies via the exact sequence (19) that

$$\psi_\omega - \alpha \psi_\mu = e_0^* \theta_{L,L'} \tag{53}$$

for some $\theta_{L,L'} \in H^1(L; \mathbb{R})$.

**7.1.2 Theorem.** *Under the assumptions of Theorem 2.2.2(a), suppose additionally that the class $\theta_{L,L'} \in H^1(L; \mathbb{R})$ defined above satisfies $\theta_{L,L'} = \iota_L^* \rho$ for some $\rho \in H^1(M; \mathbb{R})$. (Recall that $\iota_L : L \hookrightarrow M$ is the embedding). Then $I_F(L, L')$ is invariant under symplectic isotopies in the sense of (11).*

*Proof.* This is an analog of Corollary 11.2 in [L]. The idea is that if one can find an $X \in \mathcal{X}$ such that $[\mathcal{Y}_X]\big|_{\ker \psi} = 0$, then all other $X' \in \mathcal{X}$ can be connected to $X$ via a path satisfying the condition of Proposition 7.1.1, and one would then have invariance under symplectic isotopies by Proposition 7.1.1. Such an $X$ exists under our assumptions: any (path of) symplectic vector fields with flux $\rho$ will do, because of (53) and (26). □

*Proof of Theorem 2.2.3.* In this case $L' = \phi(L)$ for a symplectomorphism $\phi$ connected to $\phi_0 = \text{Id}$ via a symplectic isotopy $\{\phi_t | t \in [0,1]\}$. According to §3.1.3,

$$I_F^{\gamma_\phi}(L, \phi(L); J, \chi_H) = (\Phi)_* I_F^{\gamma_p}(L, L; J', X + \chi_{H'}), \tag{54}$$

where $X = \{X_t\}$ is the path of symplectic vector fields generating $\phi_t$.

On the other hand, in this case the fibration $e_0$ in (18) has a section, and the associated homotopy splits to yield

$$\pi_1(\Omega) = \pi_1(L) \oplus \pi_2(M, L),$$

and hence the decomposition

$$H_1(\Omega) = H_1(L) \oplus \pi_2(M, L). \tag{55}$$

Under this decomposition,

$$\psi_\mu = 0 \oplus N_{L,\mu}; \; \psi_\omega = \iota_L^* \theta_\phi \oplus N_{L,\omega}, \tag{56}$$

$\theta_\phi$ being the flux of the symplectic isotopy $\{\phi_t\}$. Thus the conditions of Theorem 7.1.2 are satisfied when $L$ is monotone, and we have invariance under symplectic isotopies. This combines with (54) to establish the first assertion of Theorem 2.2.3.

Next, by the invariance under symplectic isotopy we may compute $I_F^{\Phi \cdot \gamma_0}(L, L'; J, \chi_H)$ via $I_F^{\gamma_0}(L, L; J, \chi_h)$ for a $t$-independent $J$ and a small $t$-independent $h \in \mathcal{H}$. The latter may be identified with the classical invariant under the additional assumption that $\pi_2(M, L) = 0$: Under this assumption, (55) simplifies to say that $e_{0*} : H_1(\Omega; \mathbb{Z}) \to H_1(L; \mathbb{Z})$ is an isomorphism, and (56) says that $\ker \psi = H_1(\Omega; \mathbb{Z}) = H_1(L; \mathbb{Z})$, and $[\mathcal{Y}_{\chi_h}] = 0$. Hence $\mathcal{M}_O^{\gamma_p}(L, L; J, \chi_h)$ is empty, and the energy of all finite energy flows in $\mathcal{M}_P^{\gamma_p}(L, L; J, \chi_h)$ are *uniformly* bounded by $\max(h) - \min(h)$. We can now apply [F89a] Theorem 3 and [Oh96b] Proposition 4.1, which shows that $I_F^{\gamma_p}(L, L; J, \chi_h)$ equals the Reidemeister torsion of the Morse complex of $h$, which in turn equals $\tau(L)$. □

*Proof of Theorem 2.3.1.* Again, the strategy is to find an $X$ such that $[\mathcal{Y}_X] = -\underline{\text{im}}^* \omega + \underline{e}^* \theta_X$ is trivial over $\ker \psi_c$, and then apply Proposition 7.1.1. Such an $X$ exists by Lemma 4.1.1 (a): it suffices to take any $X$ with $[\theta_X] = \theta_\omega$. □

## 7.2 Proof of Proposition 7.1.1

This Proposition is proved for a different version of symplectic Floer theory (Floer theory of symplectomorphisms) in [L]. In view of the length of [L], we shall not reproduce arguments or definitions from [L]; instead, we assume familiarity with [L] and shall only indicate how it should be modified.

Roughly, the strategy of proof in [L] is as follows. When two regular pairs $(J_1, X_1)$, $(J_2, X_2)$ can be connected by a path $(J_\lambda, X_\lambda)$, $\lambda \in [1, 2]$ satisfying the conditions described in the Proposition, one may perturb this path into an "admissible" one ([L] Theorem 7.2), for which all the relevant parametrized moduli spaces are smooth and compact except at the points of bifurcation. There are two types of possible bifurcations: "handleslides", where there is a flow line from a nondegenerate critical point to another with the same index (possibly itself), and "birth-deaths", where there is a degenerate critical point which is standard in the sense of [L] Definition 7.1. At most one bifurcation occurs at a $\lambda$ for an admissible bifurcation. The smoothness and compactness results of the parametrized moduli spaces show that $\tau_F$ and $\zeta_F$ only change at bifurcations. The main part of the proof consists of the bifurcation analysis—analyzing how the moduli spaces change at the bifurcation to see how $\zeta_F$ and $\tau_F$ change. This involves explicit descriptions of the behavior of flow lines near a degenerate critical point, and a series of gluing theorems for flow lines.

### 7.2.1 The Lagrangian intersection version.

The proof in [L] may be transcribed almost verbatim to this version, noting the following changes:

(a) The proof of transversality and compactness results for parametrized moduli spaces (in [L] §7.1) should be modified in the manner of §3.2.2.

(b) In [L] section 10.3 (constructing nonequivariant perturbations for type II handleslides), the map $\pi$ is replaced in our case by $\pi : H_1(\Omega) \to H_2(M, L \cup L')$, and again $\pi(h) \neq 0$. The exact sequence (10.21) there is replaced by

$$H_2(M, L') \xrightarrow{\pi_2} H_2(M, L \cup L') \xrightarrow{\pi_1} H_1(L);$$

so when $\iota_{L*}\pi_1\pi(h)$ is nontorsion, we are again in the simpler case A of [L], where we just need to replace $H$ by a generic nonequivariant $\tilde{H}$ on the covering space. Otherwise, we need to choose a $\nu \in H^2(M, L \cup L')$ with $\nu(h) \neq 0$, and consider nonlocal perturbations of the form

$$\wp_{\chi P}([z, \mu]) := \chi(\int_{[0,1]\times[0,1]} \mu^*\boldsymbol{\nu})\nabla P(z),$$

where $\boldsymbol{\nu}$ is a 2-form representing $\nu$. This replaces (10.22) in [L], and we may go through the rest of the subsection in [L] with these substitutions.

(c) When computing the sign of an element of $\hat{\mathcal{M}}_O^0$, it is convenient to interpret it as a spectral flow. (See [L] §10.2). In this version of Floer theory, this is given as follows:

The sign of an $\hat{u} \in \hat{\mathcal{M}}_O^0$ is given by the spectral flow from $\mathcal{D}_{(1/T,u)}$ to $1 \oplus D_{\hat{\gamma}_0}$, where $\hat{\gamma}_0$ is the constant loop at $\gamma_0$, and

$$\mathcal{D}_{(1/T,u)} : \mathbb{R} \oplus L_1^p(u^*TM) \longrightarrow \mathbb{R} \oplus L^p(u^*TM),$$

$$\mathcal{D}_{(1/T,u)} = \tilde{D}_{(1/T,u)} + d_u^* = \begin{pmatrix} 0 & \Pi_{\partial_s u} \\ T\partial_s u & D_u \end{pmatrix},$$

where $\Pi_{\partial_s u}$ denotes $L^2$-orthogonal projection to $\partial_s u$, and $d_u^*$ is the formal $L^2$-adjoint of $d_u(\lambda) = (0, \lambda\partial_s u)$. (Cf. [L] lemma 5.5).

On the other hand, notice that we need an interpretation of the index of critical points as spectral flows in the proof of [L] Lemma 10.2 (signs for Type II handleslides). In this version, this is supplied by (27).

### 7.2.2 The $S^1$-equivariant version.

Here is a list of the main changes to [L] in this version.

**(a)** CLASSIFICATION OF GENERIC BIFURCATION. ([L] §7.1)

Because in this version of Floer theory we achieve transversality by perturbing the almost complex structure $J$ instead of the symplectic vector field, the content of [L] §10.1 needs to be modified along the following sketch:

In [L] Definition 7.1 (defining "good" minimally degenerate critical points), condition 4 (a) there should be replaced by the following

*Condition.* There is a small neighborhood $U_y$ of the image of the minimally degenerate critical point $y$ in $M$, such that over $U_y$, $J_{\lambda'}$ is constant in $\lambda'$ in a small neighborhood of $\lambda$.

Then, in the proof of [L] Theorem 7.2, replace the small perturbations in $H_\lambda$ or $H^\Lambda$ used to obtain transversality in Step 2 condition (3), Step 3 and Step 4 by small perturbations in $J_\lambda$ that vanishes in $U_y$ at a death-birth bifurcation $(\lambda, y)$, or by extensions of such perturbations in the parametrized version. Such perturbations are sufficient for achieving transversality for $U_y$ is sufficiently small; see the proof of Proposition 4.4.2.

**(b)** MODIFYING THE GLUING THEOREMS.

Because of the $S^1$-invariance of $\tilde{\mathcal{A}}_X$, the gluing theorems in [L] sections 7, 8, 9 need only minor modifications.

We have already briefly touched on the main differences between the $S^1$-equivariant version here and the non-equivariant version in [L]: One is that the products in all statements of gluing theorems should be replaced in the $S^1$-equivariant version by fiber products in the manner of (42).

Another difference is that in this version, type II handleslides include flow lines starting from a nondegenerate critical critical point $y$, ending at a *different* critical point in the same critical $S^1$-orbit. Gluing such a flow with itself will result in "twisted orbits", requiring us to include holomorphic tori of all conformal structures in the definition of $\mathcal{M}_O$ (cf. §4.1).

We now make more specific comments on the necessary modifications. In [L], we divide the standard gluing procedure into four steps: (1) constructing approximate solutions (pregluings) and estimating the error; (2) proving that the relevant deformation operator has a uniformly bounded right inverse; (3) checking that the nonlinear part of the equation satisfy the required bilinear estimate; (4) the uniqueness proof, i.e. showing that all solutions near the pregluings are those obtained from the gluing procedure.

*In Step (1):* given [L], the only extra work needed is to estimate the $\mathfrak{f}$ components of the errors. These are sufficiently small; in fact they are often trivial due to the following reasons:

- Our construction of pregluing typically involves reparameterization of given solutions $u$ of the flow equation; namely we often set the pregluing
  $$w(s) := u(\gamma(s)) \quad \text{for } \gamma : \mathbb{R} \to \mathbb{R}.$$
  This contributes only $L$-components in the error, and is thus trivial in $\mathfrak{f}$ direction. (Cf. e.g. [L] §8.3).

- In a parametrized version of gluing theorem such as [L] Theorem 8.1, there are additional terms in the error due to changes in $X$ or $J$. However, $\partial_\lambda X_\lambda$ is orthogonal to the $\mathfrak{f}$ direction due to the $\gamma_0$-exactness of the symplectic vector fields $X_\lambda$. On the other hand, $\partial_\lambda J_\lambda$ vanishes near $\lambda = 0$ by the assumption that the minimally degenerate critical point $y$ at $\lambda = 0$ is good.

- In [L] Theorem 9.1, the pregluing involves a cutoff construction near $y$, resulting in an additional error term involving $\Pi^\perp_{\mathbf{e}_y} \eta_{\lambda-}$, where $\eta_{\lambda\pm}$ are the local coordinates of the pair of new critical points $y_{\lambda\pm}$. Extending the work of [L] in the manner of §4.2.3, it is easy to see that the $\mathfrak{f}$-component, $(\eta_{\lambda-})_\mathfrak{f}$, may be estimated similarly to $\Pi^\perp_{\ker A_y} \eta_{\lambda-}$:
  $$\|(\eta_{\lambda-})_\mathfrak{f}\|_{2,2,t} \leq C\gamma^{-1}\|(\eta_{\lambda-})_L\|_{2,2,t} \leq C'|\lambda|.$$

*In Step (2):* In the standard case (gluing theorems involving only nondegenerate critical points), it follows basically from the argument of Lemma 4.5.2, plus the fact that the operators to be glued approach asymptotically constant operators exponentially, as seen in §4.2.3.

When minimally degenerate critical points are involved, we again need to augment the work in [L] by estimating the $\mathfrak{f}$-component. Take [L] Theorem 8.1 for example. Let $\hat{W}'^{;+}_\chi = \mathbb{R}_\alpha \oplus W'^{;+}_\chi$, $\hat{E}'^{;+}_\chi, V_\chi$ be respectively the obvious analog of $\hat{W}'_\chi, \hat{E}_\chi, X_\chi$ in [L] §8.4 in this $S^1$-equivariant version, and adapt the definitions of $W_\chi, L_\chi$ in [L]

Definition 8.7 similarly to the adaptations of $W_u$, $L_u$ in §4.3.3. Recall that we want to show that there is no sequence $\{(\alpha, \xi) = (\alpha_\lambda, \xi_\lambda) \in \hat{W}_\chi^{\cdot,+} | \lambda \to 0\}$ with unit $\hat{W}_\chi^{\cdot,+}$-norm such that $\|\hat{E}_\chi^{\cdot,+}(\alpha,\xi)\|_{L_\chi} \to 0$ as $\lambda \to 0$. Suppose $(\alpha,\xi) = (\alpha, \xi_\mathfrak{f} + \xi_\mathfrak{n})$ is in the kernel of $\hat{E}_\chi^{\cdot,+}$. By the same reason as (40), we have $(E_\chi^{\cdot,+}\xi_\mathfrak{f})_\mathfrak{n} = 0$; thus

$$(\hat{E}_\chi^{\cdot,+}(\alpha,\xi_\mathfrak{n}))_\mathfrak{n} \to 0; \tag{57}$$
$$(E_\chi^{\cdot,+}\xi_\mathfrak{f})_\mathfrak{f} + (E_\chi^{\cdot,+}\xi_\mathfrak{n})_\mathfrak{f} + \alpha(V_\chi)_\mathfrak{f} \to 0 \quad \text{in } L_\chi\text{-norm as } \lambda \to 0. \tag{58}$$

Now (57) implies via the arguments in [L] section 8 that

$$(\alpha, \xi_\mathfrak{n}) \to 0 \text{ in } \hat{W}_\chi^{\cdot,+} \text{ norm as } \lambda \to 0. \tag{59}$$

On the other hand, writing $\xi_\mathfrak{f} = \underline{\xi}_\mathfrak{f}(s)\dot{u}$, (58) above takes the form

$$\rho(s)\underline{\xi}_\mathfrak{f}'(s) + \nu(s) = 0,$$

where $\rho(s) \geq C\,\forall s$. Because $(V_\chi)_\mathfrak{f} = 0$ by the $\gamma_0$-exact assumption on the symplectic vector fields, $\nu(s)$ may be estimated by the same argument of (41):

$$\|\sigma_\chi \nu\|_p \leq C\|\xi_\mathfrak{n}\|_{p,1} \to 0 \quad \text{as } \lambda \to 0$$

by (59). Integrating, one finds $\|\xi_\mathfrak{f}\|_{\hat{W}_\chi^{\cdot,+}} \to 0$ as well.

*The last two steps* are essentially the same as [L]. The description of the asymptotic behaviors of the flows required in Step (4) is obtained by modifying the arguments in [L] in the manner described in §4.2.3.

(c) THE FINITE COVERING TRICK. ([L] §10.3)

To apply the finite covering trick of [H], we need to construct nonequivariant perturbations on some finite cyclic coverings of $\underline{C}$. The required coverings are finite cyclic coverings obtained as quotients of the infinite cyclic covering of $\underline{C}$ determined by a primitive cohomology class $[\nu^\mathcal{C}] \in H^1(\underline{C}, \mathbb{Z})$, which is in turn required to satisfy $[\nu^\mathcal{C}](h) = 1$, where $[h] = [u]/d \in H_1(\underline{C}; \mathbb{Z})$ is primitive, $[u] \in H_1(\underline{C}; \mathbb{Z})$ is the homology class of a Type II handleslide $u$ (i.e. a flow line from a nondegenerate critical point back to itself). $d \in \mathbb{Z}^+$ is the divisibility of $[u]$. In [L], we divide the possibilities into the two cases, Case A and B; Case B is the more complicated case when nonlocal perturbations are needed. In the present version of Floer theory, Case B is the case when $\underline{e}[u]$ is torsion. Namely, there is a $k \in \mathbb{Z}^+$ such that $\underline{\text{im}}(k[u]) \in \ker c_1\big|_{\pi_2(M)}$. However, for such class $k[u]$,

$$[\mathcal{Y}_X](k[u]) = \underline{\text{im}}^*\omega(k[u]) - \underline{e}^*\theta_X(k[u]) = 0$$

when $M$ is monotone, contradicting the fact that $u$ has positive energy.

Thus, the monotonicity assumption guarantees that we are in Case A: in this case $\underline{e}[u]$ is nontorsion; furthermore since $\underline{e}$ is surjective, $\underline{e}^*$ maps a primitive class in $H^1(M;\mathbb{Z})/\mathbb{Z}[\gamma_0]$ to a primitive class in $H^1(\underline{\mathcal{C}};\mathbb{Z})$. Choose a $\gamma_0$-exact primitive class $[\theta] \in H^1(M;\mathbb{Z})$ such that $[\theta](\underline{e}[u]) \neq 0$, then we may take $[\nu^{\mathcal{C}}] = \underline{e}^*[\theta]$. A finite covering of $\underline{\mathcal{C}} = \underline{\Omega}(M;\gamma_0)$ satisfying the requirement mentioned above may therefore be obtained as $\underline{\Omega}(\tilde{M};\tilde{\gamma}_0)$, where $\tilde{M}$ is a finite cyclic covering of $M$ determined by $[\theta]$, and $\tilde{\gamma}_0$ is a lift of $\gamma_0$ in $\tilde{M}$. ($\gamma_0$ lifts to a circle in $\tilde{M}$ by the $\gamma_0$-exactness of $[\theta]$). The required nonequivariant perturbations may simply be taken to be perturbations by Hamiltonian or symplectic vector fields on $\tilde{M}$.

**(d) SIGNS IN THE GLUING THEOREMS.** ([L] §10.1)

To modify [L], the basic principle is to replace $E_u$ there by $E_u^{-+}$, and $D_u$ there by $D_{Q,u}: \mathbb{R}_q \oplus L_1^p(u^*TM) \longrightarrow \mathbb{R} \oplus L^p(u^*TM)$,

$$D_{Q,u} := \begin{pmatrix} 0 & \Pi_{\dot{u}} \\ -1/T\dot{u} & D_u \end{pmatrix}.$$

To illustrate this, we now sketch how the proof of Lemma 10.2 in [L] should be modified. First, it follows from our definition of orientation that the sign of an element $(1/T, Q, u) \in \hat{\mathcal{M}}_O^0$ is the spectral flow from $1 \oplus 1 \oplus D_0$ to $\mathcal{D}_{(1/T,Q,u)}$, where $D_0: L_1^p \to L^p$ is a complex linear differential operator with the same symbol as the Cauchy-Riemann operator.

Thus, to compute the sign of $w$, we need to compute the spectral flow from $1 \oplus 1 \oplus D_0$ to $\mathcal{D}_{1/T,Q,w}$. This is divided into three steps, (a), (b), (c) below.

*Step (a):* From $1 \oplus 1 \oplus D_0$ to $1 \oplus 1 \oplus (D_y - \delta)$. The spectral flow of this path is by definition $\mathrm{ind}(y) \mod 2$.

*Step (b):* From

$$\begin{pmatrix} 1 & 0 & 0 \\ 0 & 1 & 0 \\ 0 & 0 & D_y - \delta \end{pmatrix} \text{ to } \begin{pmatrix} 1 & 0 & 0 \\ 0 & 0 & \Pi_{\dot{y}} \\ 0 & \dot{y} & D_y \end{pmatrix};$$

the spectral flow of this path is $0 \mod 2$ by the argument of [L].

*Step (c):* From

$$\begin{pmatrix} 1 & 0 & 0 \\ 0 & 0 & \Pi_{\dot{y}} \\ 0 & \dot{y} & D_y \end{pmatrix} \text{ to } \mathcal{D}_{1/T,Q,w}.$$

The spectral flow of this is $-\mathrm{sign}(\lambda)\mathrm{sign}(u)$ again by following [L], substituting [L] Lemma 10.3 by the next lemma.

Before introducing the next lemma, we need to first describe the map $\#_{T,Q}$ and the operator $\acute{D}_{\#_{T,Q}}$ obtained by "$Q$-twisted self-splicing".

In the notation of §4.5: Let $E^{-+} \in \Sigma_P(x, [w]; x', [w'])$, where $(x, [w]), (x', [w']) = R_Q(x, [w])$ are in the same $\mathbb{R}$-orbit in $\tilde{\mathcal{P}}$. Recall that for sufficiently large $T$, $E^{-+}$ agrees respectively with $\partial_s + \mathbb{A}_{(x,[w])}$ and $\partial_s + \mathbb{A}_{(x',[w'])}$ as $s < -T/4$ and $s > T/4$. The latter two operators are in turn related by rotation along $S_t^1$ by $Q$ (cf. the discussion in the beginning of §4.5.1). Starting from the trivial $\mathbb{R}^{2n}$ bundle over $\mathbb{R}_s \times S_t^1$, define the $\mathbb{R}^{2n}$-bundle $V_{T,Q}$ over $\Theta_{T,Q}$ by identifying the fiber over $(-T/2, t) \in \mathbb{R} \times S^1$ with $(T/2, t + Q)$. It is clear from the above discussion that $E$ formally defines a "spliced" operator $D_{\#_{T,Q}} : L_1^p(\Theta_{T,Q}, V_{T,Q}) \to L_1^p(\Theta_{T,Q}, V_{T,Q})$. We need however an extension to compare with $E^{-+}$.

Note that there is a linear map
$$\#_{T,Q} : L_{1:(\delta,\delta)}^{-+,p}(S^1 \times \mathbb{R}; \mathbb{R}^{2n}; x, [w]; x', [w']) \to \mathbb{R}_q \oplus L_1^p(\Theta_{T,Q}, V_{T,Q})$$

defined as follows. Writing any $\bar{\xi} \in L_{1:(\delta,\delta)}^{-+,p}(\mathbb{R} \times S^1; \mathbb{R}^{2n}; x, [w]; x', [w'])$ as
$$\bar{\xi}(s) = q\beta(s)\Phi_{x',[w']}\dot{x}' + C\left(\beta(-s)\Phi_{x,[w]}\dot{x} + \beta(s)\Phi_{x',[w']}\dot{x}'\right) + \xi(s),$$

where $C, q \in \mathbb{R}$ and $\xi \in L_{1:(\delta,\delta)}^p(\mathbb{R} \times S^1, \mathbb{R}^{2n})$,

$$\#_{T,Q}(\bar{\xi})(\underline{s}) := \left(q, C(\beta(-\underline{s})\Phi_{x,[w]}\dot{x} + \beta(\underline{s})\Phi_{x',[w']}\dot{x}')\right.$$
$$\left. + \sum_{s=\underline{s} \mod T} \beta(s + (T+1)/2)(1 - \beta(s - (T-1)/2))\xi(s)\right)$$

for $\underline{s} \in [-T/2, T/2]$, and the circles $\{-T/2\} \times S_t^1$ and $\{T/2\} \times S_t^1$ are identified by $R_Q$ as discussed.

The required extension $\acute{D}_{\#_{T,Q}} : \mathbb{R}_q \oplus L_1^p(\Theta_{T,Q}, V_{T,Q}) \to L^p(\Theta_{T,Q}, V_{T,Q})$ is defined by
$$\acute{D}_{\#_{T,Q}}(q, \eta) = q\beta'(s)\Phi_{x',[w']}\dot{x}' + D_{\#_{T,Q}}\eta.$$

**Lemma.** *Let $E^{-+} \in \Sigma_P(x, [w]; x', [w'])$, where $(x, [w]), (x', [w']) = R_Q(x, [w])$ as above, and let $[\ker_g E^{-+}; \operatorname{coker}_g E^{-+}]$ be a K-model for $E^{-+}$. Then for sufficiently large $T$, $[\#_{T,Q} \ker_g E^{-+}; \#_{T,Q} \operatorname{coker}_g E^{-+}]$ forms a K-model of $\acute{D}_{\#_{T,Q}}$ defined above.*

Note that the K-model above might not be co-oriented with the orientation of $\det \acute{D}_{\#_{T,Q}}$ chosen in §4.5, but this does not matter for our purpose.

The proof of this Lemma follows the standard argument in Lemma 4.5.2.

**Acknowledgements.** We thank H. Hofer, E. Ionel, F. Lalonde, D. McDuff, Y.-G. Oh, T. Parker, among others, for helpful conversations, and M. Hutchings for pointing out an error in §2.1.1 in a previous version of this article. We acknowledge Harvard University, IPAM, and MSRI for hospitality and support during the (re)-writing of this manuscript.